\documentclass[smallextended]{svjour3}   %

\usepackage[T1]{fontenc}
\usepackage[utf8]{inputenc}
\usepackage{amsmath,amssymb}
\usepackage{booktabs}
\usepackage{placeins}   %
\usepackage{graphicx}
\usepackage{pgfplots}
\pgfplotsset{compat=1.16}
\usepackage{algorithm}
\usepackage{algpseudocode}
\usepackage{natbib}
\usepackage{xcolor}
\usepackage{hyperref}

\hypersetup{
  colorlinks=true,
  linkcolor=blue!60!black,
  citecolor=blue!60!black,
  urlcolor=blue!60!black
}

\spnewtheorem{observation}{Observation}{\bfseries}{\itshape}
\def\subclassname{{\bfseries Mathematics Subject Classification (2020)}\enspace}
\newcommand{\R}{\mathbb{R}}
\newcommand{\Z}{\mathbb{Z}}
\newcommand{\Q}{\mathbb{Q}}
\newcommand{\atoms}{\mathcal{A}}
\newcommand{\redatoms}{\mathcal{F}}
\newcommand{\rows}{\mathcal{R}}
\newcommand{\ints}{\mathcal{I}}
\newcommand{\lbs}{\mathcal{L}}
\newcommand{\ubs}{\mathcal{U}}
\newcommand{\supp}{V}
\newcommand{\card}[1]{\lvert #1 \rvert}
\newcommand{\nruns}[1]{#1 \ifnum#1=1 run\else runs\fi}
\newcommand{\ninst}[1]{#1 \ifnum#1=1 instance\else instances\fi}

\newcommand{\valSeeds}{5}

\newcommand{\valStoppedRuns}{53}
\newcommand{\valSolverRuns}{1{,}760}

\newcommand{\valLbCutoff}{595}
\newcommand{\valLbOptimal}{838}
\newcommand{\valLbRctwo}{13}
\newcommand{\valLbTimelimit}{314}
\newcommand{\valReplayMaxShort}{655}
\newcommand{\valReplayMaxLong}{791}
\newcommand{\valAblNInst}{10}
\newcommand{\valFailTotal}{87}
\newcommand{\valFailSegv}{49}
\newcommand{\valFailOom}{38}

\newcommand{\valFailSegvInst}{\texttt{neos-4344877-stour} and \texttt{neos-ythan}}
\newcommand{\valFailOomInst}{\texttt{enlight11}, \texttt{enlight4}, \texttt{enlight9} and \texttt{g2\_enlight\_hard}}
\newcommand{\valPlannedRuns}{5{,}580}
\newcommand{\valAllRuns}{5{,}493}
\newcommand{\valCrossGroups}{42}

\newcommand{\valStoppedTopN}{48}
\newcommand{\valStoppedTopInst}{\texttt{neos859080}, \texttt{enlight11}, \texttt{misc05inf}, \texttt{hole10}}
\newcommand{\valAbxGClMt}{40}

\newcommand{\valAbxGClFrac}{50}

\newcommand{\valAbxGClFracShort}{59}

\newcommand{\valAbxGClFull}{58}

\newcommand{\valAbxGDen}{100}

\newcommand{\valAblGClFrac}{5}

\newcommand{\valAblGClFracShort}{16}

\newcommand{\valAblGUbFracShort}{27}

\newcommand{\valAblGClFracShortTree}{16}

\newcommand{\valAblGUbFracShortTree}{27}
\newcommand{\valAblGClInit}{18}

\newcommand{\valAblGUbInit}{50}
\newcommand{\valAblGClFull}{18}

\newcommand{\valAblGUbFull}{50}
\newcommand{\valAblGDen}{50}

\newcommand{\valPairGInitUbW}{4}
\newcommand{\valPairGInitUbT}{23}
\newcommand{\valPairGInitUbL}{0}
\newcommand{\valPairGInitUbN}{27}

\newcommand{\valPairGHfixUbW}{9}
\newcommand{\valPairGHfixUbT}{35}
\newcommand{\valPairGHfixUbL}{6}
\newcommand{\valPairGHfixUbN}{50}

\newcommand{\valPairGShortW}{41}

\newcommand{\valPairGShortN}{50}

\newcommand{\valPairGInitW}{10}

\newcommand{\valPairGInitL}{9}

\newcommand{\valPairGHfixW}{8}
\newcommand{\valPairGHfixT}{31}
\newcommand{\valPairGHfixL}{11}

\newcommand{\valPairGMtIntW}{40}

\newcommand{\valPairGMtIntN}{50}

\newcommand{\valPairGFullMtW}{45}

\newcommand{\valPairGFullMtN}{50}

\newcommand{\valPairGShortMtW}{45}
\newcommand{\valPairGShortMtT}{5}
\newcommand{\valPairGShortMtL}{0}

\newcommand{\valInstGFullMtW}{9}

\newcommand{\valInstGFracMtW}{7}

\newcommand{\valInstGFracMtL}{3}
\newcommand{\valPtsRatioGMed}{4.4}

\newcommand{\valAbxSClMt}{40}

\newcommand{\valAbxSUbInt}{99}
\newcommand{\valAbxSClFrac}{44}

\newcommand{\valAbxSClFracShort}{53}

\newcommand{\valAbxSClFull}{51}

\newcommand{\valAblSUbInt}{39}
\newcommand{\valAblSClFrac}{0}

\newcommand{\valAblSClFracShort}{9}

\newcommand{\valAblSUbFracShort}{41}

\newcommand{\valAblSClFracShortTree}{10}

\newcommand{\valAblSUbFracShortTree}{45}
\newcommand{\valAblSClInit}{9}

\newcommand{\valAblSUbInit}{46}
\newcommand{\valAblSClFull}{10}

\newcommand{\valAblSUbFull}{50}
\newcommand{\valAblSDen}{50}

\newcommand{\valPairSInitUbW}{18}
\newcommand{\valPairSInitUbT}{17}
\newcommand{\valPairSInitUbL}{8}
\newcommand{\valPairSInitUbN}{43}

\newcommand{\valPairSHfixUbW}{15}
\newcommand{\valPairSHfixUbT}{20}
\newcommand{\valPairSHfixUbL}{11}
\newcommand{\valPairSHfixUbN}{46}

\newcommand{\valPairSShortW}{47}

\newcommand{\valPairSShortN}{49}

\newcommand{\valPairSInitW}{9}

\newcommand{\valPairSInitL}{24}
\newcommand{\valPairSInitN}{46}

\newcommand{\valPairSHfixW}{16}
\newcommand{\valPairSHfixT}{19}
\newcommand{\valPairSHfixL}{11}

\newcommand{\valPairSMtIntW}{35}

\newcommand{\valPairSMtIntN}{48}

\newcommand{\valPairSFullMtW}{50}

\newcommand{\valPairSFullMtN}{50}

\newcommand{\valPairSShortMtW}{49}
\newcommand{\valPairSShortMtT}{0}
\newcommand{\valPairSShortMtL}{0}
\newcommand{\valPairSShortMtN}{49}

\newcommand{\valInstSFullMtW}{10}

\newcommand{\valInstSFracMtW}{9}

\newcommand{\valInstSFracMtL}{1}
\newcommand{\valPtsRatioSMed}{3.3}

\newcommand{\valInstGFracMtLose}{\texttt{mod008inf}, \texttt{grid4\_31} and \texttt{vertexcover-grid10x19}}
\newcommand{\valInstSFracMtLose}{\texttt{vertexcover-grid10x19}}
\newcommand{\valInstGFracMtOthers}{ (by at most two atoms on the other ones)}
\newcommand{\valInstShortSentence}{With Gurobi and the instance as unit, the median difference over the seeds (\texttt{frac+short} minus \texttt{frac}) is positive on 8 of the 10 instances, zero on 1 and negative on 1 (\texttt{g503inf}, by one atom); with SCIP, it is positive on all the 10 instances.}
\newcommand{\valAblGUbOthersMax}{27}
\newcommand{\valAbxN}{20}
\newcommand{\valAbxRuns}{1{,}968}
\newcommand{\valAbxPlanned}{2{,}000}
\newcommand{\valAbxFailed}{32}
\newcommand{\valAbxFracGW}{12}

\newcommand{\valAbxFracSW}{11}

\newcommand{\valAbxShortGW}{10}

\newcommand{\valAbxShortSW}{9}

\newcommand{\valAbxFullMtGW}{12}

\newcommand{\valAbxFullMtSW}{11}

\newcommand{\valAbxLoopMtGW}{12}

\newcommand{\valAbxLoopMtSW}{11}

\newcommand{\valAbxIntMtGW}{1}

\newcommand{\valAbxIntMtGL}{7}
\newcommand{\valAbxIntMtSW}{4}

\newcommand{\valAbxIntMtSL}{6}
\newcommand{\valAbxFracMtGW}{9}
\newcommand{\valAbxFracMtGT}{9}
\newcommand{\valAbxFracMtGL}{2}
\newcommand{\valAbxFracMtSW}{10}
\newcommand{\valAbxFracMtST}{9}
\newcommand{\valAbxFracMtSL}{1}
\newcommand{\valAbxShortLoopGW}{4}
\newcommand{\valAbxShortLoopGT}{16}
\newcommand{\valAbxShortLoopGL}{0}
\newcommand{\valAbxShortLoopSW}{7}
\newcommand{\valAbxShortLoopST}{11}
\newcommand{\valAbxShortLoopSL}{1}
\newcommand{\valAbxShortLoopEndGW}{4}
\newcommand{\valAbxShortLoopEndGT}{14}
\newcommand{\valAbxShortLoopEndGL}{2}
\newcommand{\valAbxShortLoopEndSW}{3}
\newcommand{\valAbxShortLoopEndST}{14}
\newcommand{\valAbxShortLoopEndSL}{2}
\newcommand{\valAbxMcsMtGW}{12}

\newcommand{\valAbxMcsMtSW}{11}

\newcommand{\valAbxLoopMcsGW}{2}

\newcommand{\valAbxLoopMcsSW}{3}

\newcommand{\valAbxSegv}{30}
\newcommand{\valAbxOom}{2}
\newcommand{\valAbxBothMin}{6}
\newcommand{\valAbxBothMax}{10}
\newcommand{\valAbxStopped}{53}
\newcommand{\valAbxOver}{59}
\newcommand{\valAbxGUbOthersMax}{60}
\newcommand{\valAbxMedFracG}{4.5}
\newcommand{\valAblMedFracG}{51}

\newcommand{\valAbxMedFracS}{0.75}
\newcommand{\valAblMedFracS}{19}

\newcommand{\valAbxMedShortG}{0.5}
\newcommand{\valAblMedShortG}{6}

\newcommand{\valAbxMedShortS}{0}
\newcommand{\valAblMedShortS}{13}
\newcommand{\valAbxZeroShortSW}{9}
\newcommand{\valAbxZeroShortST}{8}
\newcommand{\valAbxZeroShortSL}{2}
\newcommand{\valAbxMedFullMtG}{3.5}
\newcommand{\valAblMedFullMtG}{10.5}

\newcommand{\valAbxMedFullMtS}{6}
\newcommand{\valAblMedFullMtS}{40.5}

\newcommand{\valBenchNTwo}{30}
\newcommand{\valBenchGPairs}{150}

\newcommand{\valBenchGCl}{79}
\newcommand{\valBenchGInstAll}{15}
\newcommand{\valBenchGInstAny}{16}
\newcommand{\valBenchGMinOpen}{48}

\newcommand{\valBenchGOpenUb}{70}
\newcommand{\valBenchGStopped}{8}
\newcommand{\valBenchGTimeMed}{804}
\newcommand{\valBenchSPairs}{150}

\newcommand{\valBenchSCl}{71}
\newcommand{\valBenchSInstAll}{13}
\newcommand{\valBenchSInstAny}{15}
\newcommand{\valBenchSMinOpen}{41}

\newcommand{\valBenchSOpenUb}{70}
\newcommand{\valBenchSStopped}{11}
\newcommand{\valBenchSTimeMed}{1276}
\newcommand{\valBenchInstAny}{17}

\newcommand{\valReplayLoseSentence}{ The replay can lose a match: in 3 runs with SCIP on \texttt{g503inf} the branch-and-cut tree of Phase~2 was solved to optimality, but the covering problem of the replay reached its time limit with a smaller bound, and the run is counted as unmatched.}
\newcommand{\valSpreadOpenN}{31}
\newcommand{\valSpreadWideN}{5}
\newcommand{\valSpreadInst}{\texttt{enlight11}}
\newcommand{\valSpreadSolver}{Gurobi}
\newcommand{\valSpreadMin}{96}
\newcommand{\valSpreadMax}{231}
\newcommand{\valFdBankRuns}{240}
\newcommand{\valFdBankPoints}{48{,}410}

\newcommand{\valFdBankSentence}{No point violates a bound or an integrality requirement.}

\newcommand{\valIisTrunc}{33}
\newcommand{\valIisTimeMed}{8.8}
\newcommand{\valIisTimeMax}{1386}
\newcommand{\valIisZeroRows}{15}

\newcommand{\valIisNoSub}{16}
\newcommand{\valIisRuns}{300}

\newcommand{\valFdNeosC}{480}
\newcommand{\valFdNeosRows}{17}

\newcommand{\valFdNeosRemText}{33 atoms in all with removable domains (Table~\ref{tab:v2-b30})}
\newcommand{\valFdPzeroC}{201}
\newcommand{\valFdPzeroRows}{2}
\newcommand{\valFdPzeroTot}{203}

\newcommand{\valFdPzeroRemText}{9 atoms in all with removable domains}

\newcommand{\valFdN}{12}
\newcommand{\valNoiseGInst}{7}
\newcommand{\valNoiseGRuns}{35}
\newcommand{\valNoiseGLbDiff}{6}

\newcommand{\valNoiseSRuns}{35}
\newcommand{\valNoiseSLbDiff}{18}

\newcommand{\valNoiseSentence}{With Gurobi, the final lower bound of the ablation is the larger one in 4 of the 6 pairs of runs where the two differ, and the same pairs are matched in the two experiments; with SCIP, the final lower bound of the ablation is the larger one in 7 of the 18 pairs of runs where the two differ, and two pairs are matched in one of the two experiments only (on \texttt{bell3binf} and \texttt{mod008inf}). Differences of two matched pairs with SCIP are therefore within the noise of the time limits, and SCIP is the noisier of the two solvers.}
\newcommand{\valFdTimeGN}{10}
\newcommand{\valFdTimeGMt}{11.4}
\newcommand{\valFdTimeGOurs}{16.7}
\newcommand{\valFdTimeGMtFaster}{8}
\newcommand{\valFdTimeSN}{8}
\newcommand{\valFdTimeSMt}{30.4}
\newcommand{\valFdTimeSOurs}{30.2}
\newcommand{\valFdTimeSMtFaster}{4}
\newcommand{\valFdMatchSentence}{With Gurobi, the two methods match the bounds in almost the same number of pairs: 55 of the 60 for the multi-tree loop and 54 for our method. With SCIP, the counts are 48 and 54; the difference is on \texttt{misc05inf} and \texttt{neos-942886}, where the multi-tree loop matches the bounds in 0 and 3 of the 5 seeds, respectively, and our method in all of them.}
\newcommand{\valRedMinPctTwo}{26}
\newcommand{\valRedMaxPctTwo}{66}
\newcommand{\valPreGRed}{33}
\newcommand{\valPreGAll}{26}
\newcommand{\valPreGDen}{60}
\newcommand{\valPreSRed}{21}
\newcommand{\valPreSAll}{14}

\newcommand{\valPreAllCountSentence}{ On all atoms, a subsystem that is not optimal can keep a bound of a variable together with its integrality requirement, which already enforces it; our code counts this bound and the official checker does not. This happens in one run of this configuration, whose bounds do not match, and it does not affect the counts above.}
\newcommand{\valComN}{7}
\newcommand{\valComDen}{35}
\newcommand{\valComFdMtG}{35}
\newcommand{\valComFdOursG}{34}
\newcommand{\valComRemMtG}{5}
\newcommand{\valComRemOursG}{18}
\newcommand{\valComFdMtS}{30}
\newcommand{\valComFdOursS}{34}
\newcommand{\valComRemMtS}{0}
\newcommand{\valComRemOursS}{10}
\newcommand{\valAblUnivMax}{571}
\newcommand{\valMtMcsPctG}{98}
\newcommand{\valMtMcsPctS}{97}
\newcommand{\valCompInst}{16}
\newcommand{\valCompRuns}{140}
\newcommand{\valCompIntKeptMed}{22}
\newcommand{\valCompIntKeptLtHalf}{11}
\newcommand{\valCompIntKeptMin}{0.8}
\newcommand{\valCompIntKeptMinInst}{\texttt{grid4\_31}}
\newcommand{\valCompIntKeptMax}{100}
\newcommand{\valCompIntKeptMaxInst}{\texttt{mathologer4a}}
\newcommand{\valCompIntShareMed}{28}
\newcommand{\valCompRowShareMed}{57}
\newcommand{\valIntRuns}{140}

\newcommand{\valIntNeededMed}{12}
\newcommand{\valIntNeededOne}{5}
\newcommand{\valIntNeededLess}{7}
\newcommand{\valIntNeededOneLess}{2}
\newcommand{\valIntNeededAll}{8}

\newcommand{\valIntUntestedSentence}{the time limit of 300\,s of this filter was never reached}

\newcommand{\valSepObjPairs}{50}
\newcommand{\valSepObjFailed}{18}

\newcommand{\valSepObjLam}{10^{-4}}
\newcommand{\valSepObjNoShortGMin}{11}
\newcommand{\valSepObjNoShortGMax}{21}
\newcommand{\valSepObjNoShortSMin}{6}
\newcommand{\valSepObjNoShortSMax}{7}
\newcommand{\valSepObjShortGMin}{1}
\newcommand{\valSepObjShortGMax}{14}
\newcommand{\valSepObjShortSMin}{1}
\newcommand{\valSepObjShortSMax}{5}

\newcommand{\valNbFailed}{8}

\newcommand{\valNbRadius}{20}
\newcommand{\valNbLbGW}{5}
\newcommand{\valNbLbGT}{36}
\newcommand{\valNbLbGL}{9}
\newcommand{\valNbLbGN}{50}
\newcommand{\valNbLbSW}{13}
\newcommand{\valNbLbST}{27}
\newcommand{\valNbLbSL}{7}
\newcommand{\valNbLbSN}{47}
\newcommand{\valNbProxGW}{10}
\newcommand{\valNbProxGT}{32}
\newcommand{\valNbProxGL}{8}

\newcommand{\valNbProxSW}{12}
\newcommand{\valNbProxST}{26}
\newcommand{\valNbProxSL}{9}
\newcommand{\valNbProxSN}{47}
\newcommand{\valNbNoneGW}{5}
\newcommand{\valNbNoneGT}{37}
\newcommand{\valNbNoneGL}{8}

\newcommand{\valNbNoneSW}{12}
\newcommand{\valNbNoneST}{28}
\newcommand{\valNbNoneSL}{7}
\newcommand{\valNbNoneSN}{47}
\newcommand{\valNbNoIncSentence}{ A run that reaches the time limit without a subsystem is counted as larger; there are 2 such runs with Gurobi (one of local branching and one of proximity search), on \texttt{enlight11}.}
\newcommand{\valNbLbGLbW}{8}
\newcommand{\valNbLbGLbT}{33}
\newcommand{\valNbLbGLbL}{9}
\newcommand{\valNbLbSLbW}{13}
\newcommand{\valNbLbSLbT}{26}
\newcommand{\valNbLbSLbL}{8}
\newcommand{\valNbProxGLbW}{6}
\newcommand{\valNbProxGLbT}{30}
\newcommand{\valNbProxGLbL}{14}
\newcommand{\valNbProxSLbW}{9}
\newcommand{\valNbProxSLbT}{27}
\newcommand{\valNbProxSLbL}{11}
\newcommand{\valNbNoneGLbW}{6}
\newcommand{\valNbNoneGLbT}{36}
\newcommand{\valNbNoneGLbL}{8}
\newcommand{\valNbNoneSLbW}{12}
\newcommand{\valNbNoneSLbT}{24}
\newcommand{\valNbNoneSLbL}{11}
\newcommand{\valNbHardGPts}{5{,}621}
\newcommand{\valNbHardSPts}{1{,}003}
\newcommand{\valNbHardGImpr}{10}
\newcommand{\valNbLbGPts}{4{,}186}
\newcommand{\valNbLbSPts}{321}
\newcommand{\valNbLbGImpr}{6}
\newcommand{\valNbProxGPts}{2{,}564}
\newcommand{\valNbProxSPts}{339}
\newcommand{\valNbProxGImpr}{21}
\newcommand{\valMsfBudget}{1{,}200}
\newcommand{\valOobMissing}{9}
\newcommand{\valOobMissingFull}{2}
\newcommand{\valOobMissingInst}{\texttt{enlight9}, \texttt{g2\_enlight\_hard}}
\newcommand{\valOobN}{12}
\newcommand{\valOobNG}{8}
\newcommand{\valOobNH}{4}

\newcommand{\valOobLbGW}{9}

\newcommand{\valOobLbSW}{10}

\newcommand{\valOobShGW}{9}

\newcommand{\valOobShSW}{10}

\newcommand{\valOobFmGW}{5}
\newcommand{\valOobFmGT}{3}
\newcommand{\valOobFmGL}{4}
\newcommand{\valOobFmSW}{9}
\newcommand{\valOobFmST}{3}
\newcommand{\valOobFmSL}{0}
\newcommand{\valOobPairs}{118}
\newcommand{\valOobMsfW}{32}
\newcommand{\valOobMsfT}{20}
\newcommand{\valOobMsfL}{36}
\newcommand{\valOobMsfN}{88}
\newcommand{\valOobIisW}{39}
\newcommand{\valOobIisT}{18}
\newcommand{\valOobIisL}{31}

\newcommand{\valOobNone}{30}
\newcommand{\valOobStopped}{27}
\newcommand{\valOobNoSub}{33}
\newcommand{\valOobWholeRuns}{10}
\newcommand{\valReplayAboveFull}{67}
\newcommand{\valReplayAboveFullN}{100}
\newcommand{\valReplayAboveShort}{71}
\newcommand{\valReplayAboveShortN}{99}

\newcommand{\valLbFromExact}{14}

\newcommand{\valTolRuns}{409}

\newcommand{\valMatchedRuns}{601}
\newcommand{\valMatchedExact}{535}

\newcommand{\valMatchedExactSentence}{In all the other 66 runs the second covering problem was solved to optimality, so that the smaller $\text{LB}_{\mathrm{exact}}$ comes from the longer cuts of the exact supports and not from a time limit.}
\newcommand{\valInstMatched}{17}
\newcommand{\valInstMatchedExact}{16}
\newcommand{\valLbGapMax}{15}
\newcommand{\valLbGapInst}{\texttt{misc05inf}}

\newcommand{\valStopSentence}{in 41 the check did not finish within the limit $\bar T$ either, in one the checker did not confirm an answer of infeasibility, and in 11 the cut of the feasible point found did not exclude the candidate.}
\newcommand{\valFinCalls}{197}
\newcommand{\valFinRuns}{114}
\newcommand{\valFinUndecided}{69}
\newcommand{\valOverRuns}{78}
\newcommand{\valOverMed}{185}
\newcommand{\valOverMax}{599}
\newcommand{\valOverMatched}{0}

\newcommand{\valOverSensSentence}{We repeated the counts of wins, ties and losses on the final lower bound of Table~\ref{tab:v2-pairs} and of the comparisons of Appendix~\ref{app:abl} (24 in all) without the pairs in which one of the two runs is of this kind (at most 7 pairs per count): the sign of the difference between wins and losses changes in one of them only, \texttt{full} against \texttt{frac+short+tree} with Gurobi (9--28--13 on all the pairs, 9--28--8 without them), a comparison that is close to a tie in both cases.}

\newcommand{\valMtRuleSentence}{the check with tighter tolerances that follows an answer not confirmed by the checker keeps its short limit there, and a point found after the time limit is not repaired. No run of \texttt{multitree} in the ablation and in the fixed-domain comparison was stopped for either reason; outside the benchmark the first one stopped 5 runs, all on \texttt{g2\_ej} with SCIP, where the other configurations are stopped by the checker as well (Appendix~\ref{app:details}, Table~\ref{tab:v2-oob}).}
\newcommand{\valRecRuns}{5{,}493}
\newcommand{\valRecPoints}{5{,}982{,}756}
\newcommand{\valRecCmp}{1{,}939{,}005{,}520}

\newcommand{\valRecTiesSentence}{the rounded distance is never equal to $\delta$}
\newcommand{\valDilSubsystems}{2{,}301}

\newcommand{\valDilRuns}{3{,}963}

\newcommand{\valDilOutcomeSentence}{Both solvers declared 1{,}979 of them infeasible, although 117 of these, on \texttt{bell3binf}, are feasible in exact arithmetic (Example~\ref{ex:bell3b}). For 33 of them, 22 of which were reported by runs with matched bounds, both solvers (32 of them) or one solver declare the relaxed system feasible and none declares it infeasible. A point returned satisfies all atoms with exact comparisons within $10^{-5}(1 + 10^{-4})$, but not within $10^{-5}$: it violates the relaxed system by at most the feasibility tolerance $10^{-9}$ of the solves, i.e., it lies on the border of the relaxed system. For these subsystems, on \texttt{g\_neos-1425699} (27 of them), \texttt{misc05inf}, \texttt{mod008inf} and \texttt{stein45inf}, the tolerance reading has no margin, and we count them as undecided. The other 289 subsystems remain undecided because a solve reached its time limit of 300\,s (both solves in 123 cases). They are on \texttt{enlight11}, \texttt{enlight4}, \texttt{enlight9}, \texttt{g2\_enlight\_hard}, \texttt{g2\_pigeon-08}, \texttt{g\_noswot}, \texttt{g\_p0201}, \texttt{g\_supportcase14}, \texttt{g\_supportcase16}, \texttt{hole10}, \texttt{hole9}, \texttt{neos-2626858-aoos}, \texttt{neos-4344877-stour}, \texttt{neos-942886}, \texttt{neos-ythan} and \texttt{neos859080}.}
\newcommand{\valCertN}{12}
\newcommand{\valCertTot}{16}
\newcommand{\valCertA}{14}
\newcommand{\valCertB}{14}

\newcommand{\valBtbRows}{68}
\newcommand{\valBtbNnz}{290}
\newcommand{\valBtbChanged}{42}
\newcommand{\valBtbFree}{66}
\newcommand{\valBtbRankB}{66}
\newcommand{\valBtbRankD}{65}
\newcommand{\valBtbDecNodes}{3}
\newcommand{\valBtbUb}{118}

\newcommand{\valCertOpenSentence}{The other parts reached the node or time limit also with the wider limit: (A) on \texttt{enlight4}; (B) on \texttt{bell5inf} and \texttt{grid4\_31}.}
\newcommand{\valBtbRestSentence}{}
\newcommand{\valBtbAbsRun}{1.1\times 10^{24}}
\newcommand{\valBtbRunsPhrase}{all the 10}
\newcommand{\valBelRef}{117}
\newcommand{\valBelDistinct}{190}
\newcommand{\valBelRefRuns}{124}
\newcommand{\valBelRefMatched}{123}
\newcommand{\valBelRefMatchedSub}{116}
\newcommand{\valBelMatchedDistinct}{118}

\newcommand{\valBtbAbsMin}{6.1\times 10^{23}}
\newcommand{\valBtbAbsMax}{7.9\times 10^{24}}
\newcommand{\valFlugLbExact}{4}
\newcommand{\valFlugUbMin}{23}
\newcommand{\valFlugUbMax}{25}
\newcommand{\valFlugAtoms}{37}
\newcommand{\valFlugOptTol}{23}
\newcommand{\valFlugClosedTol}{10}
\newcommand{\valFlugRunsTol}{10}
\newcommand{\valFlugPoneLp}{23}
\newcommand{\valReplayGMedMt}{1.3}
\newcommand{\valReplaySMedMt}{0.9}

\newcommand{\valReplayGMedFull}{191}
\newcommand{\valReplaySMedFull}{65}
\newcommand{\valReplayAblMax}{655}

\newcommand{\valEndImputed}{3}
\newcommand{\valBoneN}{10}
\newcommand{\valBoneMcsSentence}{ We also ran \texttt{multitree} itself after the Phase~1 of \texttt{frac+short}, i.e., with one minimum-cardinality correction per iteration and no shortening: it is ahead of \texttt{multitree} on 8 instances with Gurobi and on 10 with SCIP, and never behind, and against it the final lower bound of the loop with our cuts is larger, equal and smaller on 5, 4 and 1 instances with Gurobi, and on 6, 3 and 1 with SCIP.}

\newcommand{\valBoneMissing}{1}

\newcommand{\valBoneTreeFinGW}{5}
\newcommand{\valBoneTreeFinGT}{5}
\newcommand{\valBoneTreeFinGL}{0}
\newcommand{\valBoneTreeFinSW}{8}
\newcommand{\valBoneTreeFinST}{1}
\newcommand{\valBoneTreeFinSL}{1}
\newcommand{\valBoneTreeEndGW}{2}
\newcommand{\valBoneTreeEndGT}{4}
\newcommand{\valBoneTreeEndGL}{4}
\newcommand{\valBoneTreeEndSW}{4}
\newcommand{\valBoneTreeEndST}{1}
\newcommand{\valBoneTreeEndSL}{5}
\newcommand{\valBoneCutsFinGW}{9}

\newcommand{\valBoneCutsFinSW}{10}

\newcommand{\valBoneTreeClosedG}{16}
\newcommand{\valBoneTreeSubG}{27}
\newcommand{\valBoneTreeRunsG}{50}
\newcommand{\valBoneTreeClosedS}{9}
\newcommand{\valBoneTreeSubS}{41}
\newcommand{\valBoneTreeRunsS}{49}
\newcommand{\valBoneMtClosedG}{15}

\newcommand{\valBoneMtClosedS}{0}

\title{Lower Bounds for Minimum-Cardinality Infeasible Subsystems in Mixed-Integer Linear Programs}

\author{Matteo Fischetti}

\institute{M. Fischetti \at
  Department of Information Engineering, University of Padova, Via Gradenigo 6/A, 35131 Padova, Italy\\
  \email{matteo.fischetti@unipd.it}}

\date{Received: date / Accepted: date}

\begin{document}
\maketitle

\begin{abstract}
Given an infeasible Mixed-Integer Linear Program (MILP), we look for an infeasible subsystem of minimum cardinality. In the setting of the 2027 Land-Doig MIP Competition (MIPcc27), the elements that can be relaxed---which we call atoms---are the linear constraints, the variable bounds and the integrality requirements, each at unit cost. The problem has an implicit MILP hitting-set formulation: for every point violating a set of atoms, every infeasible subsystem must keep at least one of them. Our focus is on the lower bound on the minimum cardinality, which tells how far a given subsystem can be from a minimum one. A dominance argument removes all bound atoms of integer variables. The hitting-set MILP is solved by a branch-and-cut scheme using fractional separation followed by a cut-shortening step. We tested the method with Gurobi and SCIP on the \valBenchNTwo{} public instances of the competition with at most 1{,}000 atoms; lower and upper bounds match on about half of them. Exact certificates confirm both bounds for one run on \valCertN{} of the \valCertTot{} instances where the lower bound valid in exact arithmetic also reaches the upper bound, and on one instance they show that a subsystem accepted by the checker is feasible in exact arithmetic. An ablation shows the effectiveness of fractional separation and cut shortening.
\end{abstract}

\section{Introduction}\label{sec:intro}

Diagnosing an infeasible Mixed-Integer Linear Program (MILP) is a classical task in optimization practice: before a model can be repaired, one has to find out which of its constraints are in conflict \citep{GuieuChinneck1999,AmaldiPfetschTrotter2003,Chinneck2008}. The standard tool is an \emph{Irreducible Infeasible Subsystem} (IIS), i.e., an infeasible subsystem that becomes feasible as soon as any of its elements is removed. MILP solvers typically offer the computation of an IIS, e.g., by a deletion filter.

Unfortunately, irreducibility says nothing about size. The same model can have IISs of very different cardinality, and a deletion filter stops at the first one it meets. A small subsystem is what a modeler wants to read, see e.g.\ \citet{Chinneck1997} for linear programs, hence one is interested in an infeasible subsystem of \emph{minimum} cardinality, which we call a Min-ISS. Computing a Min-ISS is NP-hard already for linear programs \citep{AmaldiPfetschTrotter2003}, and the 2027 Land-Doig MIP Competition (MIPcc27) \citep{MIPcc27} has chosen it as its topic, with an official checker and a public set of instances.

The competition also fixes what a subsystem is made of and how its size is measured, and this is the setting we address. In its definition, which follows \citet{GuieuChinneck1999}, the variable bounds and the integrality requirements are elements of the model that can be kept or relaxed, besides the linear constraints. Each of them has unit cost, and we call all of them \emph{atoms}. The closest work on MILPs, the X-MILP framework of \citet{LeraLeri2025}, where a smallest IIS is computed by an adaptation of FORQES \citep{IgnatievPrevitiLiffitonMarquesSilva2015}, takes instead the linear constraints as the only elements. The difference matters in two ways. On the one hand, a subsystem that can drop the domain of a variable gives a sharper diagnosis, as it tells the variables whose domain takes part in the conflict from the other ones. On the instances where our bounds match (\texttt{bell3binf} excepted, Example~\ref{ex:bell3b}), the subsystems we find keep the integer domain of a median of \valCompIntKeptMed\% of the integer variables of the model, and on \valIntNeededOne{} of these instances a single variable has to be integer once the bounds are kept (Section~\ref{sec:bench30}). On the other hand, the problem is larger and harder: on the instances where we ran both settings, our implementation of a classical implicit hitting-set loop matches lower and upper bounds in almost all runs when only the constraints can be relaxed, and in very few of them when bounds and integrality can be relaxed too (Section~\ref{sec:form-fixdom}).

In this paper we investigate how to compute lower (dual) bounds for the Min-ISS problem on general MILPs in this setting, within an algorithm that also looks for a subsystem whose cardinality matches them. Our starting point is the implicit hitting-set view, well established for unsatisfiable Boolean formulas \citep{Liffiton2009,IgnatievPrevitiLiffitonMarquesSilva2015}: every point $x$ violates a set $\supp(x)$ of elements of the model, and an infeasible subsystem must keep at least one element of every such set. The resulting master is a set-covering problem with exponentially many constraints, which are generated on the fly by a separation oracle. This is the scheme of Benders' decomposition \citep{Benders1962} with combinatorial cuts \citep{CodatoFischetti2006}, and in its original form the master is solved again to integer optimality at every iteration. This is the multi-tree loop that we use as a reference (Section~\ref{sec:setup}). The modern alternative is to separate the cuts within a branch-and-cut on the master, also at its fractional points \citep{McDanielDevine1977}, and at first glance this is all that is needed. However, a direct implementation stalls on MILPs for three reasons. First, integer separation returns long sets $\supp(x)$ of the atoms violated by a point $x$, so the lower bound grows slowly. Second, the LP relaxation of the master needs many cuts before its value moves, and each cut separation costs a MILP solve. Third, a search driven by the lower bound may run for a long time without a small infeasible subsystem, i.e., a good upper bound, in hand.

Our algorithm addresses the three points in this order. A dominance argument (Observation~\ref{obs:dom}) removes the bound atoms of integer variables, which are between \valRedMinPctTwo\% and \valRedMaxPctTwo\% of all atoms on the \valFdN{} instances where we measured it. The LP relaxation of the master is solved at the root node by fractional separation at a smoothed point, and the first points returned by each separation are moved, by a sequence of feasibility checks, to points that violate fewer atoms (\emph{cut shortening}). The integer master is then solved by branch-and-cut with lazy constraints, starting from an initial subsystem built from the LP solution, with fractional separation also inside the tree, and with two hard-fixing matheuristics that search for smaller subsystems inside the incumbent, \texttt{HardFix(U)}, and inside the incumbent plus the support of the LP solution, \texttt{HardFix(W)}. An oracle call that cannot be decided, for instance within its time limit, stops the search instead of pruning, so that the bounds we report never depend on an undecided oracle answer.

We stress that our focus is on dual bounds. What we look for is a lower bound that holds for every infeasible subsystem of the model and that, when it reaches the cardinality of the subsystem found, shows that this subsystem is of minimum cardinality. What is exact in our bounds, and what relies instead on the floating-point arithmetic of the MILP solver, is stated in Section~\ref{sec:semantics}. The initial subsystem and the two matheuristics only play the role that primal heuristics have in any branch-and-cut, namely, to give the search a good incumbent early. When a small subsystem is all that is wanted, with no bound on its quality, one can resort to heuristics such as a deletion filter. This purely heuristic context is not addressed in this paper, and we make no claim on our method as a heuristic.

Note that none of the components above is new by itself. The implicit hitting-set master, the separation of its cuts within a branch-and-cut, the reduction of a violation set towards a minimal correction, and the smoothing of the separation point all have close relatives in the literature reviewed in Section~\ref{sec:related}. To the best of our knowledge, what is new is twofold. On the one hand, these components are adapted to general MILPs in a setting where variable bounds and integrality requirements can be relaxed at unit cost. This includes the dominance argument of Observation~\ref{obs:dom}, and a statement of what the resulting bounds certify when the oracle is a floating-point MILP solver. On the other hand, we report a computational study that measures the contribution of each component, and that separates the effect of the cuts from that of the search scheme.

We evaluate the algorithm with the commercial MILP solver Gurobi \citep{Gurobi13} and with the open-source solver SCIP \citep{SCIP10}, \valSeeds{} seeds per configuration, on the \valBenchNTwo{} instances of the MIPcc27 public set with at most 1{,}000 atoms after preprocessing. An ablation on \valAblNInst{} of them, repeated on the other \valAbxN{}, isolates the effect of each component, and a multi-tree implicit hitting-set loop in the spirit of \citet{IgnatievPrevitiLiffitonMarquesSilva2015} and \citet{LeraLeri2025} is our baseline. It turns out that a branch-and-cut that uses only integer separation gives worse lower bounds than the multi-tree loop. Fractional separation makes the lower bound move, and cut shortening makes it reach the upper bound. When the bounds do not match, a verified subsystem comes, with Gurobi, from the initial subsystem rather than from the hard-fixing matheuristics, whereas with SCIP the search finds one by itself in most runs. The cuts, rather than the search scheme, explain most of these results. In a further experiment on the instances of the ablation we gave the multi-tree loop the same Phase~1 and the same cut shortening of a branch-and-cut without the other components: its lower bounds improve on almost all instances, and those of the branch-and-cut are larger when the final bound is taken, but not when the bound at the end of the search is. The branch-and-cut, however, matches the bounds in more runs with SCIP and in about as many with Gurobi, and it can end with a verified subsystem also when the bounds do not match. Without our cuts, our implementation of the multi-tree loop is about as good as our method when only the constraints can be relaxed. In the setting of the competition, where bounds and integrality can be relaxed too, it is clearly behind.

The paper is organized as follows. Section~\ref{sec:related} reviews related work. Section~\ref{sec:formulation} defines the problem, the hitting-set master and the preprocessing step, and states what the computed bounds certify. Section~\ref{sec:alg} describes the algorithm, and Section~\ref{sec:setup} the computational setup. Results are reported in Section~\ref{sec:results}. Some conclusions are finally drawn in Section~\ref{sec:concl}. The appendices are meant for the reader who wants to check our numbers: Appendix~\ref{app:params} lists the parameters of the algorithm, Appendix~\ref{app:abl} gives the per-instance results and the secondary comparisons of the ablation, Appendix~\ref{app:validity} the checks on the validity of the bounds, including the exact certificates, and Appendix~\ref{app:details} the remaining details of the experiments.

\section{Related work}\label{sec:related}

\paragraph{Infeasibility analysis in LP and MILP.}
The computation of an IIS of a linear system goes back to the filtering algorithms of \citet{ChinneckDravnieks1991}, see also \citet{GleesonRyan1990}, and was extended to MILPs by \citet{GuieuChinneck1999}. Generalized filtering schemes are given by \citet{AtlihanSchrage2008}, and the book of \citet{Chinneck2008} covers the whole area. \citet{WitzigBertholdHeinz2021} study infeasibility analysis inside a MILP solver, i.e., the derivation of conflicts during branch-and-bound. As to complexity, finding an IIS of minimum cardinality is NP-hard and hard to approximate already for linear systems \citep{AmaldiPfetschTrotter2003}, and it is strongly NP-hard even for flow networks \citep{JoormannOrlinPfetsch2016}.

A different but closely related problem asks for the smallest set of constraints whose removal makes the system feasible, i.e., the Minimum IIS Cover, equivalent to the Maximum Feasible Subsystem (MaxFS) problem. It is NP-hard even for totally unimodular matrices \citep{Sankaran1993}, and its approximability is studied by \citet{AmaldiKann1995,AmaldiKann1998}. \citet{ParkerRyan1996} formulate it as a set-covering problem whose constraints are the IISs, generated only when violated, and \citet{Pfetsch2008} gives a branch-and-cut algorithm for MaxFS based on the separation of IIS inequalities. Our master~\eqref{eq:master-full} is the hitting-set dual of theirs. In the IIS-cover master the covering constraints are IISs and a solution is a set of constraints to remove. In ours, the cuts are violation sets of points, i.e., sets of constraints whose removal makes the system feasible, and a solution is an infeasible subsystem.

\paragraph{Hitting-set duality in SAT and constraint programming.}
The duality between minimal conflicts and minimal corrections is the basis of model-based diagnosis \citep{Reiter1987}, and it has been used to enumerate the minimal unsatisfiable subsets (MUSes) of a constraint set \citep{BaileyStuckey2005,LiffitonSakallah2008,LiffitonPrevitiMalikMarquesSilva2016}. Minimal correction subsets are computed by the algorithms of \citet{MarquesSilvaHerasJanotaPrevitiBelov2013} and used to speed up MUS computation by \citet{BacchusKatsirelos2015}, while QuickXplain \citep{Junker2004} finds one conflict by divide and conquer. A smallest MUS is computed by branch and bound in \citet{Liffiton2009} and by an implicit hitting-set loop in FORQES \citep{IgnatievPrevitiLiffitonMarquesSilva2015}, which alternates a minimum hitting set of the corrections found so far with a satisfiability test of that hitting set. A satisfiable set is grown to a maximal satisfiable subset, whose complement becomes a new set to hit, and the paper also computes the minimum hitting set incrementally, by one call to an incremental MaxSAT solver, starts from a collection of disjoint corrections, and gives an early upper bound by enumerating MUSes within a time limit. OCUS \citep{GambaBogaertsGuns2023} extends the scheme to cost-optimal unsatisfiable subsets for explanation. The implicit hitting-set framework has been formalized by \citet{ChandrasekaranKarpMorenoCentenoVempala2011} and \citet{MorenoCentenoKarp2013}, and is the basis of the MaxSAT solvers of \citet{DaviesBacchus2011,DaviesBacchus2013,DaviesBacchus2013b} and \citet{SaikkoBergJarvisalo2016}, which use a MILP solver for the hitting-set problem. \citet{BacchusHyttinenJarvisaloSaikko2017} add reduced-cost fixing from the LP relaxation of the hitting-set master. Implicit hitting sets have also been used for weighted CSPs \citep{DelisleBacchus2013} and for problems beyond NP \citep{SaikkoWallnerJarvisalo2016}.

The closest work to ours is X-MILP \citep{LeraLeri2025}, a framework that explains the answer to a user query on the solution of a MILP by an IIS of the infeasible system obtained when the query is added as a constraint. In X-MILP the IIS is computed by the commercial MILP solver CPLEX \citep{CPLEX22} and is not of minimum cardinality. To assess the difference, the authors also compute a smallest IIS by an adaptation of FORQES, in which a maximal feasible subsystem is computed at every iteration, and they report that the smallest IIS has more than two constraints less than the one of CPLEX in 7\% of their instances only, with a computing time that is larger by four orders of magnitude. In the code distributed with that paper, the hitting-set problem is solved from scratch at every iteration by an ILP, the hitting set is tested by a constraint-programming model, and the maximal feasible subsystem is computed by a big-M MILP that minimizes the number of relaxed constraints. X-MILP considers the linear constraints as the only elements (its experiments are on binary programs from auctions and scheduling), which corresponds to the fixed-domain variant of Section~\ref{sec:formulation}. Our \texttt{multitree} configuration is our implementation of the same scheme on general MILPs (Section~\ref{sec:setup}). The comparisons of Section~\ref{sec:results} are with this implementation, not with the codes of FORQES and X-MILP.

\paragraph{Cuts from infeasible subproblems.}
From the MILP side, our master and its oracle form a Benders' decomposition \citep{Benders1962} in which the subproblem is a feasibility check. In the original scheme the master is solved again before every call of the subproblem, which is what the implicit hitting-set loops above do. In modern implementations the cuts are instead separated within the branch-and-cut on the master, at integer and at fractional points. Cuts from the fractional points of the master go back to \citet{McDanielDevine1977}. Our branch-and-cut algorithm is of the second kind, and the \texttt{multitree} baseline of the first. The cuts of our master are close relatives of combinatorial Benders' cuts \citep{CodatoFischetti2006}, which forbid the assignments of a master that lead to an infeasible subproblem, and of logic-based Benders' cuts \citep{HookerOttosson2003}. The choice of a short cut among the many valid ones, which we address by cut shortening, is the topic of \citet{FischettiSalvagninZanette2010} in the Benders setting. A set-covering model is also at the heart of the search for small branching backdoors in \citet{FischettiMonaci2013}. Our smoothing of the separation point is in the spirit of the in--out approach of \citet{BenAmeurNeto2007}, used for disjunctive cuts by \citet{FischettiSalvagnin2010}.

\paragraph{Exact MILP.}
Exact rational MILP solvers \citep{CookKochSteffyWolter2013,EiflerGleixner2023}, with safe cutting planes \citep{EiflerGleixner2024} and verifiable certificates \citep{CheungGleixnerSteffy2017}, would be the natural tool to turn the matched bounds of Section~\ref{sec:semantics} into proofs. See also \citet{GleixnerSteffy2020} for LP with limited-precision oracles, and \citet{Szeider2026}, who builds verifiable certificates of optimality for integer programs from the floating-point output of a black-box solver. In our setting the certificate must also tie the MPS data to the subsystem and to the supports of the stored points. We use them as a reference for what our bounds do \emph{not} certify. For one run on each of \valCertTot{} instances, we attempt to certify both bounds with a small exact branching procedure of our own (Section~\ref{sec:validity}).

\section{Problem formulation}\label{sec:formulation}

Consider an infeasible MILP defined over $x \in \R^n$:
\begin{equation}\label{eq:milp}
  \begin{aligned}
  & a_i^\top x \lesseqqgtr_i b_i \;\; (i \in \rows),\\
  & x_j \ge l_j \;\; (j \in \lbs), \qquad
  x_j \le u_j \;\; (j \in \ubs), \qquad
  x_j \in \Z \;\; (j \in \ints),
  \end{aligned}
\end{equation}
where $\rows = \{1,\dots,m\}$ indexes the linear constraints, each with its sense ${\lesseqqgtr_i} \in \{\le, \ge, =\}$, $\ints \subseteq \{1,\dots,n\}$ the integer variables, and $\lbs, \ubs$ the finite lower and upper bounds. All data $(a_i,b_i,l_j,u_j)$ are taken as the exact rationals represented by the IEEE~754 binary64 numbers read from the input MPS file.\footnote{Binary64 is the 64-bit double-precision format of the IEEE~754 standard, in which solvers store real numbers. A finite binary64 number is $s \cdot 2^{e}$ with integers $s$ and $e$, $|s| < 2^{53}$, hence a rational whose denominator is a power of two. A decimal written in the MPS file, such as $0.9$, is in general not of this form, and it is read as the nearest binary64 number, with a relative error of at most $2^{-53} \approx 1.1 \times 10^{-16}$. Our data are these rounded values, taken exactly. Examples~\ref{ex:flugpl} and~\ref{ex:bell3b} show that the difference from the decimals can change the answer.}

Let $\atoms := \rows \cup \ints \cup \lbs \cup \ubs$ be the set of \emph{atoms}, i.e., of the elements of the model that can be relaxed. The four subsets are kept disjoint by labeling their elements $\mathrm{row}(i)$, $\mathrm{int}(j)$, $\mathrm{lb}(j)$ and $\mathrm{ub}(j)$. Following the removable-domain convention of MIPcc27, keeping $\mathrm{int}(j)$ in the model enforces the whole integer domain
\begin{equation}\label{eq:dom-sem}
  x_j \in \Z \;\wedge\; (x_j \ge l_j \text{ if } j \in \lbs) \;\wedge\; (x_j \le u_j \text{ if } j \in \ubs),
\end{equation}
whereas relaxing $\mathrm{int}(j)$ makes $x_j$ continuous, and its bounds $\mathrm{lb}(j)$ and $\mathrm{ub}(j)$ can then be kept or relaxed individually, each at unit cost. If some $j \in \ints \cap \lbs \cap \ubs$ has $\lceil l_j \rceil > \lfloor u_j \rfloor$, then $\{\mathrm{int}(j)\}$ is an optimal subsystem of size one. We therefore assume that every integer domain contains at least one integer. In the \emph{fixed-domain} variant, considered for instance by \citet{LeraLeri2025}, all bounds and integrality requirements are always kept and only the rows in $\rows$ can be relaxed.

For any $x \in \R^n$, its violation support $\supp(x) \subseteq \atoms$ is the set of atoms violated by $x$, where $\mathrm{int}(j)$ counts as violated whenever $x_j$ does not satisfy~\eqref{eq:dom-sem}. By definition, for every $j \in \ints$ and every $x \in \R^n$,
\begin{equation}\label{eq:dom-imp}
  \mathrm{lb}(j) \in \supp(x) \implies \mathrm{int}(j) \in \supp(x),
  \qquad
  \mathrm{ub}(j) \in \supp(x) \implies \mathrm{int}(j) \in \supp(x).
\end{equation}
A subsystem $K \subseteq \atoms$ is infeasible if and only if every $x \in \R^n$ violates at least one atom of $K$. Hence the Min-ISS problem is the set-covering program with infinitely many constraints
\begin{equation}\label{eq:master-full}
  z^*_0 := \min \left\{ \sum_{a \in \atoms} z_a \;:\; \sum_{a \in \supp(x)} z_a \ge 1 \;\;\forall x \in \R^n,\;\; z \in \{0,1\}^{\atoms} \right\},
\end{equation}
where $z_a = 1$ if atom $a$ is kept in the model, and $z_a = 0$ if it is relaxed. We call the covering constraint of a point $x$ its \emph{cut}, to distinguish it from the rows $\mathrm{row}(i)$ of the MILP.
By~\eqref{eq:dom-imp}, the variables $z_{\mathrm{lb}(j)}$ and $z_{\mathrm{ub}(j)}$ with $j \in \ints$ are dominated by $z_{\mathrm{int}(j)}$, at equal cost. Indeed, any weight on $z_{\mathrm{lb}(j)}$ or $z_{\mathrm{ub}(j)}$, binary or fractional, can be moved onto $z_{\mathrm{int}(j)}$ without increasing the objective and without uncovering any cut. The same equal-cost dominance is used in set covering and in MaxSAT hitting-set solvers \citep{DaviesBacchus2011}. We therefore work on the reduced atom set
\begin{equation}\label{eq:red-univ}
  \redatoms := \rows \;\cup\; \ints \;\cup\; (\lbs \setminus \ints) \;\cup\; (\ubs \setminus \ints),
\end{equation}
with projected supports $\supp_{\redatoms}(x) := \supp(x) \cap \redatoms$.

\begin{observation}\label{obs:dom}
For any subset of points $X \subseteq \R^n$, restricting~\eqref{eq:master-full} to the cuts of $X$, to the variables $z_a$ with $a \in \redatoms$, and to the projected supports $\supp_{\redatoms}(x)$ preserves both the integer optimum and the value of the linear relaxation.
\end{observation}

Note that the reduction also removes the corresponding binary variables from the separation MILP of Section~\ref{sec:srp}. Its effect on the size of the atom set and on the bounds is measured in Section~\ref{sec:form-fixdom}.

\subsection{What the bounds certify}\label{sec:semantics}

The lower and upper bounds that we report are obtained in two different ways, and it is worth stating precisely what each of them guarantees.

Every point $x$ found by the algorithm is a \emph{witness}: it is computed to satisfy a given set $K$ of atoms, and it is stored together with the atoms it violates. The solver returns it in floating point. After rounding its integer variables, we check whether it satisfies all the atoms of $K$ in exact arithmetic, against the binary64 data of the MPS file. If it does not, we try to turn it into a rational vector that does (Section~\ref{sec:grow}). In both cases we call the point an \emph{exact witness}. Otherwise the binary64 point is kept as a \emph{tolerance witness}, provided that it violates no atom of $K$ by more than $\delta = 10^{-9}$. A point that fails this test as well is stored with the atoms it violates, but it says nothing on $K$ (Section~\ref{sec:safe-bc}). Violations are absolute and unscaled: the residual of a row, the distance from a bound, and the distance from the nearest integer. For a row the comparison with $\delta$ is rigorous, i.e., the residual is computed together with a bound on its rounding error, and in rational arithmetic when this bound does not allow a decision. The distance from the nearest integer is computed without rounding error. The distance from a bound is the difference of two binary64 numbers, which is rounded, and in the runs it is compared with $\delta$ as it is. The outcome can differ from that of the exact comparison only when the rounded distance is equal to $\delta$. As a control, we evaluated again with the exact comparison the \valRecPoints{} binary64 points stored in all the \valRecRuns{} runs of Section~\ref{sec:results}: in the \valRecCmp{} comparisons with a finite bound, \valRecTiesSentence{}, and no row changes.

For $\varepsilon \ge 0$, let $z^*_\varepsilon$ be the smallest cardinality of a subsystem that no point satisfies within the absolute tolerance $\varepsilon$, so that $z^*_{\mathrm{exact}} := z^*_0$ is the value of~\eqref{eq:master-full}, the Min-ISS value in exact arithmetic. A subsystem with no point within $\varepsilon$ has no point within a smaller tolerance either, hence $z^*_\varepsilon$ is nondecreasing in $\varepsilon$. For a stored point $x$, let $\supp^0(x)$ be the set of the atoms of $\redatoms$ that $x$ violates in exact arithmetic, and $\supp^\delta(x) \subseteq \supp^0(x)$ the set of those that it violates by more than $\delta$. For a point stored as a rational vector we set $\supp^\delta(x) := \supp^0(x)$. During the search, the cut of a point in the covering master is $\supp^\delta(x)$. Every subsystem that no point satisfies within $\delta$ contains an atom of $\supp^\delta(x)$, because $x$ satisfies all the other atoms within $\delta$. In the same way, every subsystem that is infeasible in exact arithmetic contains an atom of $\supp^0(x)$. The reduction to $\redatoms$ applies to both kinds of cuts, as a bound of an integer variable $x_j$ that is violated by more than $\delta$ puts $\mathrm{int}(j)$ in $\supp^\delta(x)$. After the search, the two supports are recomputed for every stored point. The set-covering problem with one cut $\supp^\delta(x)$ per point is then solved again by Gurobi or SCIP, in up to three stages, each with its own time limit (Appendix~\ref{app:details}). If the search has an upper bound UB and was completed, i.e., Phase~2 ended by itself and without undecided oracle calls, the covering problem is first solved with the extra constraint $\sum_a z_a \le \text{UB}-1$. Its infeasibility gives the bound UB directly. Otherwise, the covering problem is solved to optimality, and at the time limit the bound is the dual bound of the solver rounded up. In that case the MaxSAT solver RC2 \citep{IgnatievMorgadoMarquesSilva2018,IgnatievMorgadoMarquesSilva2019} is also tried, and the larger of the two values is kept. We call the result $\text{LB}_\delta$. By the argument above, $\text{LB}_\delta \le z^*_\delta$. Of the \valSolverRuns{} runs of the benchmark, of the ablation and of the fixed-domain comparison (Section~\ref{sec:validity}), \valLbCutoff{} end in the first stage, \valLbOptimal{} with an optimal covering solve, \valLbRctwo{} with RC2, and \valLbTimelimit{} with the dual bound at the time limit. The same covering problem with the cuts $\supp^0(x)$, solved by the last two stages only, gives a second bound $\text{LB}_{\mathrm{exact}} \le z^*_{\mathrm{exact}}$, valid in exact arithmetic. Since $z^*_{\mathrm{exact}} \le z^*_\delta$, the value $\text{LB}_{\mathrm{exact}}$ is a valid lower bound on $z^*_\delta$ as well, and the \emph{replayed lower bound} that we report, called the \emph{final lower bound} in Section~\ref{sec:results}, is $\text{LB} := \max\{\text{LB}_\delta, \text{LB}_{\mathrm{exact}}\} \le z^*_\delta$. Every cut $\supp^0(x)$ contains the cut $\supp^\delta(x)$ of the same point, hence the optimal value of the second covering problem cannot exceed that of the first one, and $\text{LB}_{\mathrm{exact}} > \text{LB}_\delta$ is possible only when the first problem ends at its time limit (\valLbFromExact{} runs), or when its first stage gives $\text{LB}_\delta = \text{UB}$ and $\text{LB}_{\mathrm{exact}} > \text{UB}$, which never occurred. The 0/1 data of the two covering problems make their solution numerically benign, although the solve itself is done in floating point and is not accompanied by a proof certificate. For one run on each of \valCertTot{} instances, Section~\ref{sec:validity} attempts such a certificate in exact arithmetic.

The upper bound UB is the cardinality of a subsystem accepted as infeasible by the official MIPcc27 checker \citep{MIPcc27}, with its default feasibility tolerance of $10^{-5}$. The checker decides infeasibility by a floating-point solve of the subsystem with the same solver used in the run, so UB inherits the reliability of that solver. We write $z^*_{\mathrm{tol}} := z^*_\varepsilon$ with $\varepsilon = 10^{-5}$. As $\delta \le \varepsilon$, we have $z^*_\delta \le z^*_{\mathrm{tol}}$. What an accepted subsystem says on these values depends on how the verdict of the checker is read, and two readings have to be kept apart. In the \emph{exact reading}, the accepted subsystem has no point in exact arithmetic; then $z^*_{\mathrm{exact}} \le \text{UB}$, and nothing follows for a positive tolerance, as a subsystem can be infeasible in exact arithmetic and have a point within $\varepsilon$ (Example~\ref{ex:flugpl}). In the \emph{tolerance reading}, the accepted subsystem has no point within $\varepsilon$; then $z^*_{\mathrm{tol}} \le \text{UB}$. The tolerance reading is an assumption on the checker, which a correct verdict of its solver does not imply. The exact reading is an assumption as well, and Example~\ref{ex:bell3b} shows that it can fail. Appendix~\ref{app:validity} reports a check of both readings on its runs. Under the tolerance reading,
\begin{equation}\label{eq:sandwich}
  \text{LB}_{\mathrm{exact}} \;\le\; z^*_{\mathrm{exact}} \;\le\; z^*_\delta \;\le\; z^*_{\mathrm{tol}} \;\le\; \text{UB}
  \qquad\text{and}\qquad \text{LB}_{\mathrm{exact}} \;\le\; \text{LB} \;\le\; z^*_\delta ,
\end{equation}
whereas the exact reading only gives $\text{LB}_{\mathrm{exact}} \le z^*_{\mathrm{exact}} \le \text{UB}$. We use the words \emph{matched bounds} for the equality LB $=$ UB of the two computed values, and we do not claim more. When the bounds match, under the tolerance reading the optimal values $z^*_\delta$ and $z^*_{\mathrm{tol}}$ coincide with the matched value, up to the floating-point covering solve. If in addition $\text{LB}_{\mathrm{exact}} = \text{UB}$, the exact reading is enough to conclude that the matched value is $z^*_{\mathrm{exact}}$. This is the stronger of the two statements, and Appendix~\ref{app:validity} says in how many runs it holds. The converse does not hold: as the first stage is not used for $\text{LB}_{\mathrm{exact}}$, in a run with matched bounds a smaller $\text{LB}_{\mathrm{exact}}$ can also be the dual bound of a covering solve that ended at its time limit. Two remarks are in order. First, the verdict of the checker comes from a solver with its own scaling and presolve, and is not literally a statement on the points within an absolute tolerance. Should the checker accept a subsystem that has a point within $\delta$, LB could exceed UB; we compared the two bounds in every run, and this never occurred (Section~\ref{sec:validity}). This is a consistency check, and not a proof of the tolerance reading: an accepted subsystem with a point goes undetected as long as LB does not exceed its cardinality, as in Example~\ref{ex:bell3b}. Second, $z^*_{\mathrm{exact}}$ and $z^*_{\mathrm{tol}}$ can differ: a small subsystem can be infeasible in exact arithmetic only because of the binary64 representation of a coefficient, and feasible within the tolerance. On such an instance exact witnesses alone cannot match the two bounds, and the gap is due to the two definitions of infeasibility rather than to the search. This is not a remote possibility, as the following example shows.

\begin{example}\label{ex:flugpl}
In the benchmark instance \texttt{flugplinf}, the equality row \texttt{ANZ4} reads $0.9\,\mathrm{STM3} + \mathrm{ANM3} - \mathrm{STM4} = 0$, where the three variables are integer, with $57 \le \mathrm{STM3}, \mathrm{STM4} \le 75$ and $0 \le \mathrm{ANM3} \le 18$. The binary64 number closest to $0.9$ is $q = 8106479329266893/2^{53}$, whose numerator is odd, hence $q\,\mathrm{STM3}$ is integer only if $\mathrm{STM3}$ is a multiple of $2^{53}$. In exact arithmetic, the row and the integer domains of its three variables are therefore an infeasible subsystem of 4 atoms. On the other hand, the integer point $(60, 3, 57)$ violates the row by $3/2^{51} \approx 1.3 \times 10^{-15}$ only, and the checker finds this subsystem feasible with both solvers. The same holds for the rows \texttt{ANZ3}, \texttt{ANZ5} and \texttt{ANZ6} of the model. We ran the configuration \texttt{full} on this instance in two ways. With exact witnesses only ($\delta = 0$), the lower bound is \valFlugLbExact{} in every run, so that $z^*_{\mathrm{exact}} = \valFlugLbExact$ up to the floating-point covering solve, whereas the subsystems accepted by the checker have between \valFlugUbMin{} and \valFlugUbMax{} of the \valFlugAtoms{} atoms of the model, and the bounds cannot match. With tolerance witnesses, the bounds match at \valFlugOptTol{} in \valFlugClosedTol{} of the \valFlugRunsTol{} runs. The lower bound alone gives $z^*_\delta \ge \valFlugOptTol > \valFlugLbExact = z^*_{\mathrm{exact}}$, so the two notions of infeasibility differ on this instance whatever the reading of the checker. Under the tolerance reading, $z^*_\delta = z^*_{\mathrm{tol}} = \valFlugOptTol$. The value $z^*_{\mathrm{exact}} = \valFlugLbExact$ also holds without the covering solve. Ten of the points stored by a run of an earlier version of our code with Gurobi and seed~0 satisfy, in exact arithmetic with the binary64 data, every set of at most three of the 37 atoms of $\redatoms$, so that by Observation~\ref{obs:dom} no subsystem of fewer than four atoms is infeasible; the ten points and the script that checks them are available from the author upon request.
\end{example}

The binary64 representation can also act the other way, and turn an infeasible subsystem into a feasible one.

\begin{example}\label{ex:bell3b}
In the benchmark instance \texttt{bell3binf}, the run of Table~\ref{tab:cert} has matched bounds at \valBtbUb{}, and its subsystem of \valBtbUb{} atoms is accepted by the checker with both solvers. Its \valBtbRows{} rows are inequalities of type $\le$, with \valBtbNnz{} nonzero coefficients, and \valBtbFree{} of its continuous variables have no finite bound. With the decimals written in the MPS file, this subsystem is infeasible: our certifier of Appendix~\ref{app:validity} proves it with a tree of \valBtbDecNodes{} nodes, as a combination of the rows cancels all the variables with no finite bound and leaves a contradiction. However, \valBtbChanged{} of the coefficients, such as $1.0585$, are not binary64 numbers, and after rounding the cancellation fails in the last bits: the rows restricted to the variables with no finite bound have rank \valBtbRankD{} with the decimals and \valBtbRankB{} with the binary64 data. These variables can then make up for the contradiction, at very large values. We computed a rational point that satisfies all the \valBtbUb{} atoms in exact arithmetic with the binary64 data, with integer values for the integer variables and continuous values up to $\valBtbAbsRun$ in absolute value. With the data of this paper, the subsystem is therefore feasible in exact arithmetic, and within any tolerance. The same construction gives an exact point of the subsystem reported in \valBtbRunsPhrase{} runs of \texttt{full} on this instance in Table~\ref{tab:v2-b30}, where the largest absolute value is between $\valBtbAbsMin$ and $\valBtbAbsMax$.\valBtbRestSentence{} In exact arithmetic, these ten points also satisfy \valBelRef{} of the \valBelDistinct{} distinct subsystems of this instance reported as upper bound by the runs of Section~\ref{sec:results}, in \valBelRefRuns{} runs, \valBelRefMatched{} of which with matched bounds. Gurobi and SCIP did not find such points: in the check of Appendix~\ref{app:validity}, they declare all these subsystems infeasible even when every atom is relaxed by the tolerance of the checker. The lower bound is not affected: part (B) of Table~\ref{tab:cert} proves $z^*_{\mathrm{exact}} \ge \valBtbUb$, and what fails on this instance is the upper bound.
\end{example}

\section{The algorithm}\label{sec:alg}

Algorithm~\ref{alg:overview} gives the overall scheme. It keeps a finite subset of points $X \subset \Q^n$, which defines the current covering master, and a verified incumbent subsystem $U$. Phase~1 solves the linear relaxation of the master by fractional separation. Phase~2 solves the integer master by branch-and-cut, with the two matheuristics of Section~\ref{sec:hardfix}. Every point found anywhere is added to $X$, and the final lower bound is computed from $X$ only, as explained in Section~\ref{sec:semantics}. The parameters are collected in Table~\ref{tab:params}.

\begin{algorithm}[t]
\caption{Branch-and-cut for Min-ISS (overview)}\label{alg:overview}
\begin{algorithmic}[1]
\Require infeasible MILP~\eqref{eq:milp}, time budget $T$
\State compute $\redatoms$ (Observation~\ref{obs:dom}); initialize the set of points $X \subset \Q^n$ by disjoint seeding (Section~\ref{sec:srp})\label{alg:init}
\Repeat \Comment{Phase~1, at most $0.45\,T$}
  \State solve LP~\eqref{eq:lp-master} on $X$, obtaining $z^*$; update the smoothed point $\bar z$\label{alg:smooth}
  \State solve the separation MILP~\eqref{eq:sep-mip}--\eqref{eq:sep-dom} at $\bar z$;\label{alg:sep} shorten the supports of the first two points (Section~\ref{sec:grow}); add the points to $X$
\Until{the LP is closed, or no progress, or four failed iterations at $z^*$, or the Phase-1 time cap is reached}
\State $U \gets$ initial subsystem from the LP solution, verified by the checker
\State solve the integer master on $X$ by branch-and-cut: \Comment{Phase~2}
\State \quad at an integer solution $z$, with $C := \{a : z_a = 1\}$: if $C$ is feasible, add to $X$ a point $x$ that satisfies $C$, with shortened support (Section~\ref{sec:grow}); if $C$ is infeasible, it becomes the new incumbent $U$; if $C$ is still undecided after the repeated check, stop the search
\State \quad at selected nodes: fractional separation of the LP solution $\hat z$ of the node
\State \quad at a new incumbent: \texttt{HardFix(U)}; on primal stall: \texttt{HardFix(W)}
\State recompute the supports $\supp^\delta$ and $\supp^0$ of the points of $X$ and solve the two covering problems (Section~\ref{sec:semantics}) $\to$ LB; \Return LB, $U$
\end{algorithmic}
\end{algorithm}

\subsection{Phase 1: fractional separation}\label{sec:srp}

Given $X$, we solve the covering LP on $\redatoms$ without the redundant upper bounds $z_a \le 1$:
\begin{equation}\label{eq:lp-master}
  L(X) := \min \left\{ \sum_{a \in \redatoms} z_a \;:\; \sum_{a \in \supp_{\redatoms}(x)} z_a \ge 1 \;\;\forall x \in X,\;\; z \ge 0 \right\}.
\end{equation}
Here and in the rest of this section, the cut $\supp_{\redatoms}(x)$ of a point $x$ is the set $\supp^\delta(x)$ of Section~\ref{sec:semantics}. Omitting the redundant upper bounds $z_a \le 1$ guarantees that the dual multipliers $\lambda_x \ge 0$ give dual loads $\ell_a := \sum_{x \in X : a \in \supp_{\redatoms}(x)} \lambda_x \le 1$ for every atom. The dual load $\ell_a$ is the price that the cuts of $X$ put on atom $a$, and $1 - \ell_a \ge 0$ is the reduced cost of $z_a$.

At Step~\ref{alg:init} of Algorithm~\ref{alg:overview}, the initial points in $X$ come from two short disjoint-seeding passes: the first one looks for points violating integrality atoms only, the second one for points violating rows only, and every new support is excluded from the next calls so that the supports found are pairwise disjoint.

At Step~\ref{alg:smooth} of Algorithm~\ref{alg:overview}, as a stabilization technique \citep{BenAmeurNeto2007}, the point to be separated at iteration $t$ is the smoothed point $\bar z^{(t)} = (1-\alpha)\,\bar z^{(t-1)} + \alpha\, z^*$ with $\alpha = 0.6$. We go back to $z^*$ itself ($\alpha = 1$) as soon as the smoothed point cannot be separated, and we return to $\alpha = 0.6$ at every iteration whose index is a multiple of 15, provided that $z^*$ itself is in use and the iteration has added a cut. A safeguard completes the loop. If the exact repair of Section~\ref{sec:grow} loses the violation of the first point of a separation, the feasibility of the atoms satisfied by that point is checked again with tight tolerances, and the point so found is used. If this check declares the atoms infeasible, the set is excluded from the next separation MILPs by a no-good constraint (the last 30 such sets are kept). Note that $\bar z$ is a moving average of the LP optima, not a fixed interior point as in the in--out approach of \citet{BenAmeurNeto2007}.

We call \emph{backbone} the set $B$ of atoms $a$ for which $X$ contains a point $x$ with $\supp_{\redatoms}(x) = \{a\}$. Such a point satisfies all atoms but $a$ within $\delta$, so $a$ belongs to every subsystem of $\redatoms$ that has no point within $\delta$, and to every infeasible one if $\supp^0(x) = \{a\}$. In the master, the corresponding cut reads $z_a \ge 1$. Backbone atoms are recognized whenever a point with a singleton support is found, and are also searched for explicitly by a probe, done at the first iteration and then every 40, that checks the feasibility of $\redatoms \setminus \{a\}$ for up to 20 atoms with $z^*_a \ge 0.99$.

At Step~\ref{alg:sep} of Algorithm~\ref{alg:overview}, given $\bar z$, let $S := \{a \in \redatoms : \bar z_a > 0\}$ be its support and $H := B \cup \{a \in S : \bar z_a \ge 0.9995\}$ the \emph{hard set}. The backbone is not implied by the threshold: $\bar z$ is the smoothed point, and for an atom just recognized as backbone $z^*_a = 1$ while $\bar z_a$ can stay below $0.9995$ for several iterations. We separate $\bar z$ by the MILP
\begin{align}
  \min\;\; & \sum_{a \in S \setminus H} \bar z_a\, v_a \label{eq:sep-mip}\\
  \text{s.t.}\;\; & \sum_{a \in S \setminus H} \bar z_a\, v_a \le 0.999, \label{eq:sep-cut}\\
  & a_i^\top x \lesseqqgtr_i b_i, && \mathrm{row}(i) \in H, \label{eq:sep-hard-first}\\
  & x_j \in \Z,\;\; l_j \le x_j \le u_j, && \mathrm{int}(j) \in H, \\
  & x_j \ge l_j, && \mathrm{lb}(j) \in H, \\
  & x_j \le u_j, && \mathrm{ub}(j) \in H, \label{eq:sep-hard-last}\\
  & v_a = 0 \;\Rightarrow\; a_i^\top x \lesseqqgtr_i b_i, && a = \mathrm{row}(i) \in S \setminus H, \label{eq:sep-ind-first}\\
  & v_a = 0 \;\Rightarrow\; x_j = y_j, \;\; y_j \in \Z,\;\; l_j \le y_j \le u_j, && a = \mathrm{int}(j) \in S \setminus H, \label{eq:sep-ind-int}\\
  & v_a = 0 \;\Rightarrow\; x_j \ge l_j, && a = \mathrm{lb}(j) \in S \setminus H, \\
  & v_a = 0 \;\Rightarrow\; x_j \le u_j, && a = \mathrm{ub}(j) \in S \setminus H, \label{eq:sep-ind-last}\\
  & x \in \R^n,\;\; v_a \in \{0,1\}, && a \in S \setminus H. \label{eq:sep-dom}
\end{align}
Constraints \eqref{eq:sep-hard-first}--\eqref{eq:sep-hard-last} enforce all atoms of $H$, i.e., they impose $H \cap \supp_{\redatoms}(x) = \emptyset$. The indicator constraints \eqref{eq:sep-ind-first}--\eqref{eq:sep-ind-last} enforce atom $a$ when $v_a = 0$. In~\eqref{eq:sep-ind-int}, the auxiliary integer variable $y_j$ carries the integer domain~\eqref{eq:dom-sem} of $x_j$ (only finite bounds are imposed on $y_j$). When the separation is run on all atoms (Section~\ref{sec:form-fixdom}), $\mathrm{lb}(j)$ and $\mathrm{ub}(j)$ are both in $H$ and $u_j - l_j \le 10^4$, the indicator~\eqref{eq:sep-ind-int} is replaced by the big-$M$ pair $-(u_j - l_j)\, v_a \le x_j - y_j \le (u_j - l_j)\, v_a$. On the reduced atom set the bounds of the integer variables are not atoms, and the pair is never used. Atoms outside $S$ do not appear at all. Hence, in exact arithmetic, $v_a \ge [\,a \in \supp_{\redatoms}(x)\,]$ for every $a \in S \setminus H$, and any feasible solution gives a point $x$ with $\sum_{a \in \supp_{\redatoms}(x)} \bar z_a \le 0.999 < 1$, i.e., a cut for~\eqref{eq:lp-master} violated by $\bar z$.

The two sets $S$ and $H$ simply remove from the MILP the atoms whose status is known in advance. An atom $a \notin S$ has $\bar z_a = 0$, so violating it costs nothing, and it is relaxed (no variable $v_a$). An atom $a \in H$ has $\bar z_a \ge 0.9995$, so any point violating it gives $\sum_{a \in \supp(x)} \bar z_a \ge 0.9995 > 0.999$, i.e., a cut that $\bar z$ violates by less than the separation threshold; hence $v_a = 0$ can be imposed beforehand, and $a$ is enforced as a hard constraint of the MILP. For a backbone atom $a$ with $\bar z_a < 0.9995$, imposing $v_a = 0$ only loses cuts that are dominated by the cut $z_a \ge 1$ already in the master. We stop each call as soon as the solver has found 11 feasible solutions, each of which gives a cut violated by $\bar z$. The covering LP~\eqref{eq:lp-master} is \emph{closed} when the separation MILP~\eqref{eq:sep-mip}--\eqref{eq:sep-dom}, written at the LP optimum $z^*$ instead of $\bar z$, is infeasible, i.e., when no point gives a cut of weight at most $0.999$ at $z^*$. If every cut of the master is the exact support $\supp^0(x)$ of its point, if LP~\eqref{eq:lp-master} and the separation MILP are solved in exact arithmetic for the same cuts, and if every no-good constraint excludes a set of atoms with no point, then every cut of~\eqref{eq:master-full} has weight larger than $0.999$ at $z^*$, hence $z^*/0.999$ is feasible for the whole covering LP, whose value is between $L(X)$ and $L(X)/0.999$. In our runs these conditions are not guaranteed, and the closure is only a stopping rule, which is not used for the bounds we report. The separation MILP is solved in floating point, with the tolerances of the solver. Moreover, the cut of a point $x$ during the search is $\supp^\delta(x)$, which can be strictly contained in $\supp^0(x)$ (Section~\ref{sec:semantics}). This happens when $x$ is a tolerance witness, i.e., a point that satisfies the atoms of its set $K$ only within $\delta$, and also when $x$ is an exact witness in binary64 that satisfies an atom outside $K$ only within $\delta$. The value $L(X)$ then refers to the covering LP of the subsystems with no point within $\delta$, and it can exceed the optimal value of~\eqref{eq:master-full} itself. On \texttt{flugplinf}, Phase~1 ends with $L(X) = \valFlugPoneLp$ in all the runs of the benchmark, whereas $z^*_{\mathrm{exact}} = \valFlugLbExact$ (Example~\ref{ex:flugpl}). By Observation~\ref{obs:dom}, the MILP has no indicator variables for the bounds of integer variables.

Phase~1 ends at the first of: closed LP; 35 iterations without an improvement of $10^{-3}$ in $L(X)$, after at least 45\,s; four failed iterations at $z^*$; a time cap of 45\% of the time limit. Note that the cap makes the length of Phase~1, and therefore the time at which the bounds of a run match, depend on the budget.

\subsection{Cut shortening}\label{sec:grow}

The points $x$ returned by~\eqref{eq:sep-mip}--\eqref{eq:sep-dom}, or by the integer separation of Phase~2, often violate many atoms, and a long support gives a weak cut. Given such a point $x$, with violated set $V \subseteq \redatoms \setminus K$ for the set $K$ of atoms that $x$ must satisfy, cut shortening greedily tries to move as many atoms of $V$ as possible into the satisfied set $\redatoms \setminus V$. We try to re-impose the atoms of $V$ one block at a time, i.e., to remove them from the violated set. When the satisfied set together with a block is feasible, the block is absorbed: the point $x'$ found by the feasibility check replaces $x$, and $V$ is replaced by the violated set of $x'$, which can lose more atoms than those of the block. In detail:
\begin{enumerate}
  \item if $\card{V} < 6$, the cut is short already and $x$ is kept as it is;
  \item the atoms of $V$ are sorted by decreasing $z^*_a + 0.5\,\ell_a$, where $\ell_a$ is the dual load of Section~\ref{sec:srp}, with the atoms of the incumbent $U$ first in Phase~2. The atoms with large LP value or dual load are thus tried first, so that the remaining violations tend to fall onto atoms of small LP value: the cut gets shorter by construction, and this order aims at making it also more violated;
  \item the sorted atoms are tested in consecutive blocks of $\min(8, \lfloor \card{V}/4 \rfloor)$ atoms (hence between 3 and 8) when $\card{V} \ge 12$, and one by one otherwise. Each check is a pure feasibility MILP on the satisfied set together with the block, with a short time limit, and a block whose check fails is split in two;
  \item the last point $x$, which violates exactly the final set $V$ (the original point if no block was absorbed), is turned into a witness as described in Section~\ref{sec:semantics}, namely by rounding its integer variables and, if the rounded point does not satisfy in exact arithmetic the atoms of the final satisfied set, by a sparse rational projection onto the active rows (skipped above 200 rows). If this repair fails, the floating-point point $x$ is kept as a tolerance witness. The point $x$ is added to $X$, and its cut is $\supp^\delta(x)$.
\end{enumerate}
The feasibility checks run in floating point, and only the last point $x$ is repaired. The time limits per call and per point are given in Table~\ref{tab:params}. To limit its cost, cut shortening is applied to the first two points returned by a separation MILP in Phase~1, and to the first one in Phase~2; the other points of the same call are added with the support they have.

In the terminology of the SAT literature, a \emph{correction subset} is a set of constraints whose removal makes the system feasible, and a \emph{minimal correction subset} (MCS) is one that is minimal by inclusion. The violated set $V$ of any point is a correction subset, and it is minimal if and only if no atom of $V$ can be re-imposed without losing feasibility. Cut shortening is therefore a time-limited version of the ``grow'' procedure used to compute an MCS \citep{MarquesSilvaHerasJanotaPrevitiBelov2013}, which enlarges a satisfiable set of constraints as long as it stays satisfiable. Because of the time limits, the final violated set is short but not necessarily minimal. When it is minimal, no point has a violation set strictly contained in $V$, i.e., the corresponding cut is not dominated by any other cut of~\eqref{eq:master-full}.

\subsection{Phase 2: Branch-and-Cut}\label{sec:safe-bc}

Phase~2 solves the integer covering master on $\redatoms$, with the cuts of $X$ and lazy constraints: in Gurobi through \texttt{LazyConstraints=1} and the \texttt{MIPSOL} callback, in SCIP through a constraint handler (with dual reductions disabled, since the cuts are not all known). At an integer solution $z$ of the master, which keeps the atoms of $C := \{a : z_a = 1\}$, the handler first looks for a violated cut among the points of $X$. If there is none, it checks whether the subsystem $C$ is feasible. On large models ($\card{\redatoms} > 350$ or $\card{\redatoms \setminus C} > 150$) this is a pure feasibility MILP on $C$, on small ones a MILP that minimizes the number of violated atoms outside $C$. A feasible point is shortened and gives a new lazy cut. If $C$ is infeasible, the official checker is called and $C$ becomes the new incumbent.

\emph{Undecided oracle calls and the time limit.} The feasibility check of a candidate $C$ (the \emph{oracle}) can end without a usable answer in three ways: it does not finish within its time limit, even after a second attempt with a longer limit and a different seed; the solver declares the kept atoms infeasible, but the checker does not confirm it and a check with tighter tolerances ($10^{-9}$ on feasibility and integrality) finds no point; or a feasible point is found, but the witness obtained from it (Section~\ref{sec:semantics}) violates an atom of $C$ by more than $\delta$, so that its cut does not exclude $C$. A candidate in this state can be neither accepted nor rejected, and the search cannot go on without an answer. The first case arises in two situations. During the run, the limits of the two attempts are short (Table~\ref{tab:params}), and a harder check is left undecided although most of the time of the run is still available. At the end of the run, the check is cut by the time limit of the run, and the candidate that is lost in this way may be the one that improves the incumbent or matches the bounds. In both situations we let the check finish: a check that does not finish within its limits is repeated once, whatever the time left in the run, with a larger time limit $\bar T$, set to 300\,s in our runs. The checker and the checks with tighter tolerances have the same limit $\bar T$. The price is that a run can end after its time limit. The time limit is tested before every oracle call, and a single call can chain up to four checks with limit $\bar T$ (the repeated check, the checker and two checks with tighter tolerances). We report the actual running times, and Appendix~\ref{app:validity} says how many runs end after their limit, by how much, and what the comparisons become without them. Only if the candidate is still undecided after that does the handler stop the search instead of rejecting $C$. The run then ends with the current incumbent and with the replayed lower bound of the points found so far, so that no reported bound depends on an undecided answer. In the \texttt{multitree} baseline (the multi-tree loop of Section~\ref{sec:related}, whose implementation is described in Section~\ref{sec:setup}), the check of the candidate is repeated in the same way, with two differences: \valMtRuleSentence{} The stopped runs are counted separately in Tables~\ref{tab:v2-abl-G}--\ref{tab:v2-abl-S}, and Appendix~\ref{app:validity} lists them by cause.

\emph{In-tree fractional separation.} At the root and then every 200 nodes (SCIP also at nodes of depth at most 3 among the first 15), we separate the fractional node solution $\hat z$ by~\eqref{eq:sep-mip}--\eqref{eq:sep-dom} with threshold 0.995 in~\eqref{eq:sep-cut}, for at most 20 calls and 12\% of the Phase-2 time. No smoothing is applied inside the tree. Hence an atom fixed to~1 by branching has $\hat z_a = 1$ and belongs to the hard set $H$ because of the threshold. The cut obtained is still globally valid, since its support is the violation set of an actual point.

\subsection{Hard-fixing matheuristics}\label{sec:hardfix}

As in any branch-and-cut, it pays to find a good incumbent early in the search. To this end we use three primal components, all inside the same branch-and-cut code and all adding to $X$ every point they find.
\begin{enumerate}
  \item \emph{Initial subsystem.} At the end of Phase~1, starting from the backbone and from the atoms with $z^*_a > 0.05$, we add atoms until the subsystem is infeasible, and then remove them in blocks, in increasing order of $z^*_a$, as long as the subsystem stays infeasible. The result is verified by the checker and is the first incumbent $U$.
  \item \texttt{HardFix(U)}. At every new incumbent $U$, we solve the Min-ISS problem restricted to $U$, i.e., with $z_a = 0$ fixed for all $a \notin U$, with a short time limit (see Table~\ref{tab:params}). The restricted problem is small, and we address it by a cutting-plane procedure of its own instead of a nested branch-and-cut. The cuts of $X$ are first projected onto $U$. Each atom $a \in U$ is then probed by a short feasibility check of $U \setminus \{a\}$: if the check finds a point, its support gives a new cut and $a$ is needed in $U$, whereas $a$ is removed from $U$ if $U \setminus \{a\}$ is infeasible. A few rounds of fractional separation restricted to $U$ follow. Finally, the covering problem on the projected cuts, with the additional constraint $\sum_{a \in U} z_a \le \card{U} - 1$, is solved from scratch with a time limit of 5\,s every time a cut is added, as in a multi-tree loop. Each of its solutions is checked for feasibility and gives either a smaller subsystem or a new cut. If this covering problem becomes infeasible and no check was left undecided, then no subsystem of $U$ is smaller than $U$, i.e., $U$ is \emph{irreducible} with respect to the removal of atoms; otherwise the run ends with a smaller subsystem inside $U$, or with a lower bound on the size of the smallest one. Note that this statement of irreducibility relies on floating-point solves that are not repeated by the replay of Section~\ref{sec:semantics}. All the points found are valid for the whole problem and are added to $X$.
  \item \texttt{HardFix(W)}. When $U$ is irreducible, \texttt{HardFix(U)} cannot improve it further. If the incumbent has not improved for a while (see Table~\ref{tab:params} for this trigger and the size of $W$), we enlarge $U$ with the atoms outside $U$ with largest LP value $\hat z_a > 0.01$, obtaining $W$. We then address the problem restricted to $W$, with $z_a = 0$ for all $a \notin W$, looking for a subsystem smaller than $U$. The procedure is the one of \texttt{HardFix(U)}, where the probing step is replaced by two attempts to drop from $W$ the 3 and the 6 atoms of $U$ with smallest $\hat z_a + 0.5\,\ell_a$: if what remains is still infeasible, it is reduced by the deletion filter of the initial subsystem. A success gives a new incumbent. In any case, the cuts of the points found are valid for the whole problem and are added to $X$, although they are not guaranteed to cut off the current LP solution, whose weight outside $W$ can be positive.
\end{enumerate}
Fixing $z_a = 0$ outside $W$ makes both the restricted master and its separation MILP small, as the atoms outside $W$ disappear from them. The alternative is to keep all the atoms and to define the neighbourhood of $U$ by a constraint or by the objective function of the master, as in local branching \citep{FischettiLodi2003} and in proximity search \citep{FischettiMonaci2014}. We compare the three options in Section~\ref{sec:abl}.

\section{Computational setup}\label{sec:setup}

The algorithm is implemented in Python~3.12 on top of Gurobi~13.0.3 (through \texttt{gurobipy}) and SCIP~10.0 (through PySCIPOpt~6.3.0) \citep{Gurobi13,SCIP10}. Either solver is used for all MILPs of a run, i.e., for the master, the separation and feasibility MILPs, and the final covering problem, while the covering LP~\eqref{eq:lp-master} is solved by HiGHS through SciPy~1.18 \citep{HuangfuHall2018}. Each solver runs on one thread, and all its other parameters, including its own cutting planes, presolve and primal heuristics, keep their default values, except for the lazy-constraint settings of Section~\ref{sec:safe-bc} and the time and solution limits of Section~\ref{sec:alg}. The runs of the benchmark, of the ablation and of the fixed-domain comparison (Sections~\ref{sec:bench30}--\ref{sec:form-fixdom}, except two side experiments of Section~\ref{sec:abl}, see Appendix~\ref{app:details}) were made on the identical blades of a partition of a departmental cluster of the University of Padova, with an Intel Xeon E3-1220~v2 CPU at 3.10\,GHz (launched in 2012), 4 cores and 14\,GB of memory, allocated exclusively, each running four single-thread processes at a time that share its memory. The configurations of our code compared in a table run on the same type of machine, while the IIS finder used as a reference in Section~\ref{sec:bench30} ran on a second partition; Appendix~\ref{app:details} gives the details, and the machines of the other experiments. Each configuration was run with \valSeeds{} random seeds (0--4) on each instance, and the unit of our statistics is the instance-seed pair. The seeds of the same instance are of course not independent. The search time limit is 600\,s for the ablation and the fixed-domain comparison, and 1{,}200\,s for the benchmark of Section~\ref{sec:bench30}. The replay of Section~\ref{sec:semantics} runs after the search has ended, and its time is not part of these limits.

\paragraph{Instances.}
The MIPcc27 public set (\texttt{smallset}) contains 50 infeasible MILPs. We use all the \valBenchNTwo{} instances with at most 1{,}000 atoms after preprocessing. The other ones have more atoms and are outside the scope of this study. The ablation uses \valAblNInst{} of these instances and the fixed-domain comparison \valFdN{}. Both subsets were chosen during the development of the code, among instances not solved within seconds and of different type, and are not a random sample. The ablation is then repeated on the other \valAbxN{} instances (Section~\ref{sec:abl30}). Note that the parameters of Table~\ref{tab:params} were set during the development of the code on instances of the same set, so the results measure the behavior of the algorithm on the instances it was designed on, not on an independent test set. A check on \valOobN{} instances outside this set is reported in Section~\ref{sec:oob}.

\paragraph{Configurations.}
The ablation compares eight configurations, each defined by its difference from the previous ones:
\begin{itemize}
  \item \texttt{multitree}: our implementation of the multi-tree implicit hitting-set loop of \citet{IgnatievPrevitiLiffitonMarquesSilva2015}, as adapted to MILPs by \citet{LeraLeri2025}. This is the original scheme of \citet{Benders1962}, in which the master is solved again before every call of the subproblem, whereas all the other configurations separate their cuts within a branch-and-cut on the master. At each iteration, the integer master is solved from scratch, with a time limit of $\min(30, 0.25\,\mathrm{rem})$\,s, where $\mathrm{rem}$ is the remaining time, after which its best solution $C$ is used; $C$ is checked for feasibility, and a minimum-cardinality correction outside $C$ is computed by a MILP when $\card{\redatoms \setminus C} \le 600$, with a time limit of $\min(60, 0.15\,\mathrm{rem})$\,s and no other stopping rule, its best solution being used at the time limit (otherwise, or if this MILP finds no solution, the raw feasible point is used). Unlike FORQES, the hitting set is not computed incrementally, and no MUS enumeration gives an early upper bound. It shares the seeding and the replay with our code, and the feasibility checks of the oracle with the differences of Section~\ref{sec:safe-bc} and Appendix~\ref{app:details};
  \item \texttt{int}: branch-and-cut with integer separation only, i.e., with cuts generated only by the oracle calls at integer solutions of the master, after the seeding;
  \item \texttt{frac}: as \texttt{int}, with Phase~1 fractional separation;
  \item \texttt{frac+short} and \texttt{frac+tree}: as \texttt{frac}, with cut shortening, and with in-tree fractional separation, respectively;
  \item \texttt{frac+short+tree}: as \texttt{frac}, with both;
  \item \texttt{+init}: as \texttt{frac+short+tree}, with the initial subsystem of Section~\ref{sec:hardfix};
  \item \texttt{full}: as \texttt{+init}, with \texttt{HardFix(U)} and \texttt{HardFix(W)}.
\end{itemize}
For the benchmark we also report, as a reference for our upper bounds, the IIS finder of each solver (Gurobi \texttt{computeIIS}, SCIP~10 IIS), mapped onto the atoms of $\redatoms$ and followed by a deletion filter in one random order. Both subsystems are verified by the checker, and the filtered one is reported when it passes; Appendix~\ref{app:details} gives the time limits.

\paragraph{Implementation notes.} The following details matter for the interpretation of the results.
\begin{itemize}
  \item The time limits of the separation and feasibility MILPs are in wall-clock seconds, so the sequence of calls, and hence the final bounds, can change from one run to another even with the same seed. This is one source of the variance across seeds.
  \item \texttt{multitree} is our reimplementation, not the original code of \citet{LeraLeri2025}. It uses the same seeding and replay as our method, and the same feasibility checks in its oracle (with one attempt instead of two, see Appendix~\ref{app:details}), but not Phase~1 and cut shortening, so the comparison is between the two search schemes with their own way of producing cuts. A multi-tree loop with our cuts is the object of a separate experiment in Section~\ref{sec:abl}.
  \item All the configurations other than \texttt{multitree}, \texttt{int} included, are branch-and-cut algorithms on the master, with the cuts separated at its integer solutions; \texttt{tree} in a name refers to fractional separation inside the tree. The comparison with the scheme that solves the master from scratch at every iteration is therefore the one between \texttt{multitree} and \texttt{int}, with the caveat on the cuts of the previous item, and, with our cuts on both sides, the one of the experiment just mentioned.
  \item A few runs with the SCIP backend of our code did not end normally (Section~\ref{sec:validity}). A failed run counts as not matched, with the planned \valSeeds{} seeds as denominator.
\end{itemize}

\section{Results}\label{sec:results}

We first report the results of \texttt{full} on the benchmark (Section~\ref{sec:bench30}) and the ablation of its components (Sections~\ref{sec:abl} and~\ref{sec:abl30}), then the fixed-domain mode (Section~\ref{sec:form-fixdom}) and the instances outside the benchmark (Section~\ref{sec:oob}). Section~\ref{sec:validity} closes with the checks on the validity of the bounds, on which all these comparisons rely.

\subsection{Benchmark}\label{sec:bench30}

Table~\ref{tab:v2-b30} reports the configuration \texttt{full} on the \valBenchNTwo{} instances with at most 1{,}000 atoms, with a time limit of 1{,}200\,s. With Gurobi, lower and upper bounds match in every seed on \valBenchGInstAll{} instances, in at least one seed on \valBenchGInstAny{}, and in \valBenchGCl{} of the \valBenchGPairs{} instance-seed pairs. With SCIP, they match in every seed on \valBenchSInstAll{} instances, in at least one seed on \valBenchSInstAny{}, and in \valBenchSCl{} of the \valBenchSPairs{} pairs. Overall, \valBenchInstAny{} instances have matched bounds with at least one solver; Section~\ref{sec:validity} reports what these matches certify. Of the unmatched runs, \valBenchGStopped{} with Gurobi and \valBenchSStopped{} with SCIP stopped on an undecided oracle call, and are marked in the table. The medians hide a wide spread over the seeds on some instances, up to \valSpreadMax{} against \valSpreadMin{} atoms on \valSpreadInst{} with \valSpreadSolver{} (Appendix~\ref{app:details}).

The subsystems with matched bounds also show what the removable-domain setting gives to the modeler. On the \valCompInst{} instances with at least one matched run other than \texttt{bell3binf} (\texttt{flugplinf} included), the subsystems found keep the atom $\mathrm{int}(j)$, i.e., the integer domain, of a median of \valCompIntKeptMed\% of the integer variables of the model, and of less than half of them on \valCompIntKeptLtHalf{} instances. A kept $\mathrm{int}(j)$ also carries the bounds of $x_j$, so we also checked how much integrality itself is needed, by dropping it one variable at a time from each matched subsystem (Appendix~\ref{app:details}). Integrality is needed in every subsystem, but often for few variables: for a single one in every run of \valIntNeededOne{} instances, and for a median of \valIntNeededMed\% of the integer variables of the model. A diagnosis of this kind, made of some rows, of the domain of a minority of the variables and of the integrality of fewer still, is not available when the domains are fixed.

The last two columns of the table give, as a reference for our upper bounds, the subsystem found by the IIS finder of each solver followed by a deletion filter. This is a method of a different kind. It is a heuristic, which gives no lower bound and hence cannot tell how far its subsystem is from the minimum, and it is usually fast, with a median time of \valIisTimeMed\,s and a maximum of \valIisTimeMax\,s (on a different machine, see Appendix~\ref{app:details}), whereas a run of our method goes on until its time limit unless the bounds match. We therefore do not compare the two methods as heuristics, which is outside the scope of this paper: the reader can see in the table that neither subsystem is always the smaller one, and that on some of the instances where our search ends with a gap, such as \texttt{stein45inf}, the filter does better.

\begin{table}[t]
\centering
\scriptsize
\setlength{\tabcolsep}{3pt}
\caption{The full method on the 30 MIPcc27 \texttt{smallset} instances with at most 1{,}000 atoms after preprocessing, 5 seeds (1{,}200\,s of search). Column $|\redatoms|$: number of atoms after Observation~\ref{obs:dom}. For each solver: median [LB, UB] over the seeds, i.e., over the runs that ended and, for UB, that have a verified subsystem, a dash if at most half of them have one (bold single value if every seed matched), number $k$ of seeds with matching bounds, followed, when $k > 0$, by the number of them with $\text{LB}_{\mathrm{exact}} = \text{UB}$ as well (Section~\ref{sec:semantics}), and median search time over the matched seeds; $^{\dagger n}$ = $n$ runs failed; $^{\ast n}$ = $n$ runs stopped on an undecided oracle call. Time is given only when all seeds matched. Matching bounds are results of the checker (Section~\ref{sec:semantics}); $^{\mathrm{c}}$ = the optimum is certified in exact arithmetic (Table~\ref{tab:cert}); $^{\mathrm{r}}$ = subsystems returned by the runs are feasible in exact arithmetic (Example~\ref{ex:bell3b}). Columns IIS + filter: median over the seeds of the cardinality of the subsystem given by the IIS finder of the solver followed by a deletion filter, verified by the checker; a dash means no verified subsystem.}
\label{tab:v2-b30}
\begin{tabular}{@{}lrrrrrrr@{}}
\toprule
 & & \multicolumn{2}{c}{Gurobi} & \multicolumn{2}{c}{SCIP} & \multicolumn{2}{c}{IIS + filter} \\
\cmidrule(lr){3-4}\cmidrule(lr){5-6}\cmidrule(l){7-8}
Instance & $|\redatoms|$ & [LB, UB] ($k$) & s & [LB, UB] ($k$) & s & Gurobi & SCIP \\
\midrule
\texttt{mathologer4a}$^{\mathrm{c}}$ & 5 & \textbf{4} (5; 5) & 0.1 & \textbf{4} (5; 5) & 0.1 & 4 & -- \\
\texttt{stein10i}$^{\mathrm{c}}$ & 24 & \textbf{12} (5; 5) & 1.4 & \textbf{12} (5; 5) & 3.8 & 15 & 15 \\
\texttt{flugplinf} & 37 & \textbf{23} (5; 0) & 7.6 & \textbf{23} (5; 0) & 152.3 & 23 & 24 \\
\texttt{neos-4331170-senio}$^{\mathrm{c}}$ & 37 & \textbf{12} (5; 5) & 4.7 & \textbf{12} (5; 5) & 6.4 & 12 & 12 \\
\texttt{enlight4} & 48 & \textbf{18} (5; 5) & 376.5 & [18, 18] (4; 4)$^{\ast1}$ & -- & 19 & -- \\
\texttt{stein15inf}$^{\mathrm{c}}$ & 52 & \textbf{42} (5; 5) & 36.2 & \textbf{42} (5; 5) & 21.6 & 43 & 43 \\
\texttt{g503inf}$^{\mathrm{c}}$ & 89 & \textbf{30} (5; 3) & 775.5 & [28, 30] (0) & -- & 52 & 54 \\
\texttt{aim-50-2\_0-no-4}$^{\mathrm{c}}$ & 150 & \textbf{27} (5; 5) & 21.1 & \textbf{27} (5; 5) & 67.2 & 27 & 28 \\
\texttt{bell5inf} & 196 & \textbf{94} (5; 5) & 221.8 & \textbf{94} (5; 5) & 479.4 & 106 & 107 \\
\texttt{neos-ythan} & 233 & [73, 129] (0) & -- & [70, 160.5] (0)$^{\dagger1}$ & -- & 179 & 129 \\
\texttt{bell3ainf}$^{\mathrm{c}}$ & 257 & \textbf{117} (5; 5) & 287.6 & \textbf{117} (5; 4) & 330.8 & 124 & 123 \\
\texttt{bell3binf}$^{\mathrm{r}}$ & 257 & \textbf{118} (5; 5) & 574.8 & \textbf{118} (5; 5) & 758.2 & 137 & 126 \\
\texttt{aim-100-1\_6-no-4}$^{\mathrm{c}}$ & 260 & \textbf{59} (5; 5) & 129.2 & \textbf{59} (5; 5) & 150.8 & 62 & 61 \\
\texttt{neos859080} & 324 & [176, 227] (0)$^{\ast5}$ & -- & [174, 250] (0)$^{\ast5}$ & -- & 277 & 322.5 \\
\texttt{mod008inf}$^{\mathrm{c}}$ & 326 & \textbf{323} (5; 3) & 120.9 & \textbf{323} (5; 2) & 718.2 & 324 & 324 \\
\texttt{p0201inf}$^{\mathrm{c}}$ & 335 & \textbf{9} (5; 5) & 1.2 & \textbf{9} (5; 5) & 2.2 & 9 & 9 \\
\texttt{enlight11} & 363 & [89, 199.5] (0)$^{\ast2}$ & -- & [38.5, 325] (0)$^{\dagger3}$ & -- & 120 & -- \\
\texttt{grid4\_31} & 373 & [93, 95] (0) & -- & [94, 95] (2; 2) & -- & 106 & 110 \\
\texttt{stein45inf} & 377 & [22, 241] (0) & -- & [21, 229] (0) & -- & 211 & 211 \\
\texttt{misc05inf} & 431 & [97, 122] (0)$^{\ast1}$ & -- & [88, 121] (0) & -- & 146 & 128 \\
\texttt{domset-hexgrid\_03\_052} & 469 & [158, 188] (0) & -- & [158, 227] (0) & -- & 239 & 245 \\
\texttt{hole9} & 505 & [436, 480] (0) & -- & [434, 490] (0) & -- & 481 & 487 \\
\texttt{aim-200-1\_6-no-2}$^{\mathrm{c}}$ & 520 & \textbf{101} (5; 5) & 265.3 & \textbf{101} (5; 5) & 780.9 & 106 & 105 \\
\texttt{vertexcover} & 571 & [197, 200] (0) & -- & [196, 200] (0) & -- & 200 & 200 \\
\texttt{hole10} & 671 & [584, 643] (0) & -- & [582, 658] (0)$^{\ast5}$ & -- & 644 & 652 \\
\texttt{neos-942886}$^{\mathrm{c}}$ & 839 & [33, 33] (4; 4) & -- & [24, 42] (0) & -- & 72 & 36 \\
\texttt{jnh215} & 900 & [65, 104] (0) & -- & [55, 120] (0) & -- & 275 & 274 \\
\texttt{jnh10} & 950 & [59, 67] (0) & -- & [51, 89] (0) & -- & 253 & 259 \\
\texttt{jnh15} & 950 & [57, 141] (0) & -- & [46, 141] (0) & -- & 265 & 308 \\
\texttt{neos-4344877-stour} & 990 & [70, 125] (0) & -- & --$^{\dagger5}$ & -- & 105 & 105 \\
\bottomrule
\end{tabular}
\end{table}

\subsection{Ablation}\label{sec:abl}

Table~\ref{tab:v2-pairs} summarizes the ablation, i.e., the eight configurations of Section~\ref{sec:setup} on the \valAblNInst{} instances of the ablation, through paired comparisons of two configurations on the same instance and seed; the bounds of every configuration on every instance are in the upper blocks of Tables~\ref{tab:v2-abl-G} and~\ref{tab:v2-abl-S} in Appendix~\ref{app:abl}. Two runs that should be identical can differ, because the time limits are in wall-clock seconds: \texttt{full}, run twice with the same seeds on \valNoiseGInst{} instances, gives a different final lower bound (the replayed bound LB of Section~\ref{sec:semantics}) in \valNoiseGLbDiff{} of the \valNoiseGRuns{} pairs of runs with Gurobi and in \valNoiseSLbDiff{} of \valNoiseSRuns{} with SCIP (Appendix~\ref{app:abl}). Small differences, in particular with SCIP, are within this noise.

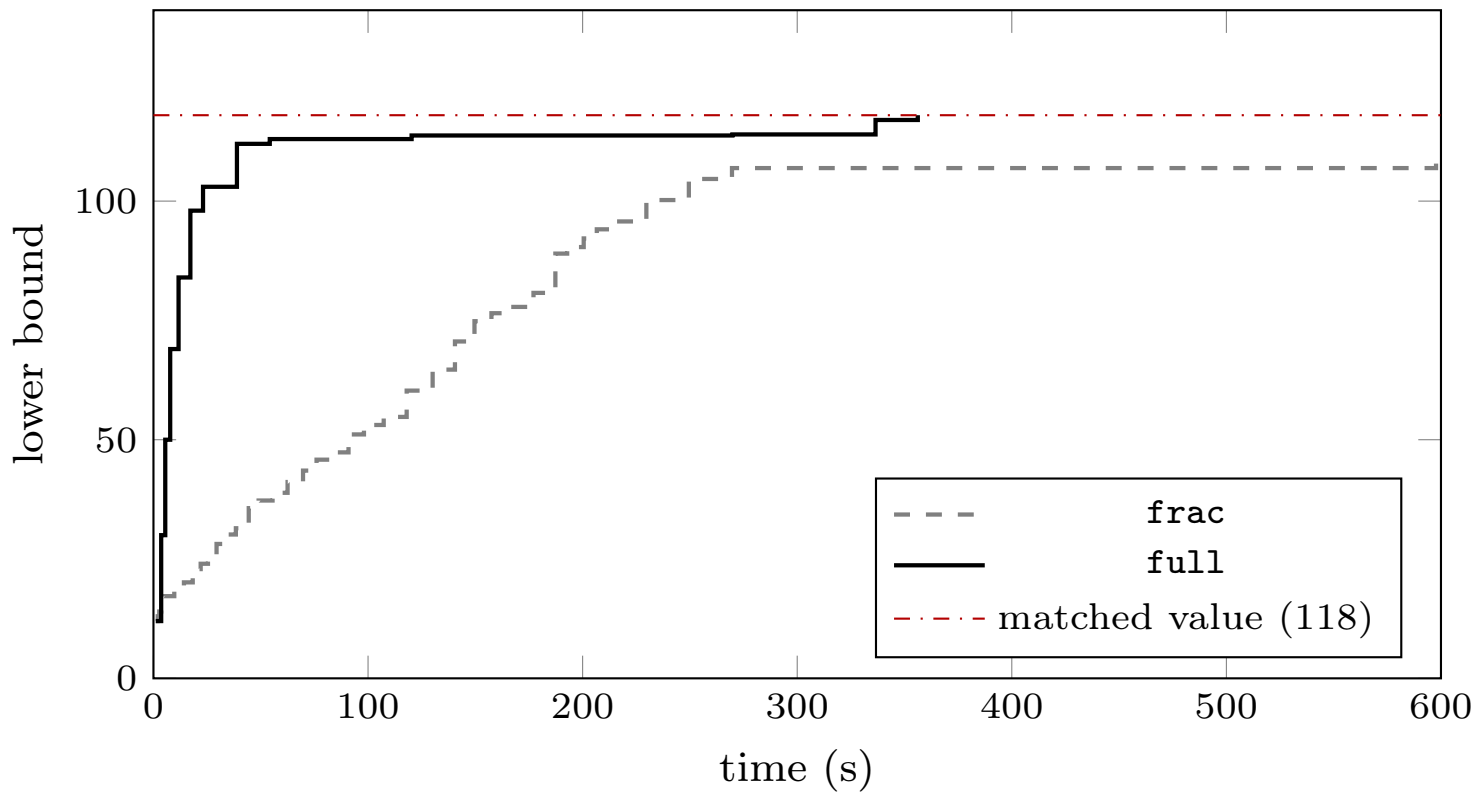
\begin{figure}[t]
\centering
\begin{tikzpicture}
\begin{axis}[width=0.85\textwidth, height=6cm, xlabel={time (s)}, ylabel={lower bound}, ymin=0, ymax=140, xmin=0, xmax=600, legend pos=south east, legend style={font=\scriptsize}, tick label style={font=\scriptsize}]
\addplot[dashed, gray, thick, const plot] coordinates {(1.1,12.00) (2.0,13.00) (2.7,14.00) (3.6,15.57) (4.5,17.22) (9.7,18.75) (14.2,20.11) (18.3,22.89) (22.0,24.01) (25.4,25.65) (29.4,28.16) (33.6,30.14) (38.4,31.45) (44.4,35.66) (49.0,37.26) (55.5,38.88) (62.5,41.14) (69.7,43.53) (76.0,45.82) (83.9,47.35) (90.9,51.11) (98.1,53.09) (107.3,54.78) (118.0,60.29) (130.1,64.69) (140.5,70.59) (149.6,74.84) (157.5,76.51) (166.4,77.83) (177.1,80.77) (187.3,89.00) (192.8,90.37) (200.5,92.18) (206.6,94.09) (214.2,95.74) (229.7,100.22) (249.5,104.64) (269.6,106.92) (597.7,110.00)}; \addlegendentry{\texttt{frac}}
\addplot[solid, black, thick, const plot] coordinates {(1.1,12.00) (3.6,30.00) (5.5,50.00) (7.8,69.00) (11.7,84.00) (17.2,98.00) (23.1,103.00) (38.9,112.00) (54.2,113.00) (120.3,113.75) (269.8,113.97) (269.8,113.97) (336.5,117.00) (356.1,117.00) (356.2,118.00)}; \addlegendentry{\texttt{full}}
\addplot[dashdotted, red!70!black, thin, domain=0:600] {118}; \addlegendentry{matched value (118)}
\end{axis}
\end{tikzpicture}
\caption{Lower bound during the search on \texttt{bell3binf} with Gurobi, seed 0, for \texttt{frac} and \texttt{full}: value of the covering LP in Phase~1, bound of the branch-and-cut tree in Phase~2. Horizontal line: the value 118 at which the bounds of \texttt{full} match; on this instance it is only a lower bound (Example~\ref{ex:bell3b}). The final replayed lower bound is computed after the search and is not shown.}
\label{fig:traj}
\end{figure}

Figure~\ref{fig:traj} gives a first view of the main effect, on a single run on \texttt{bell3binf} where the two configurations of the figure start from the same seeding. The figure is about the lower bound, since on this instance most subsystems accepted by the checker are feasible in exact arithmetic (Example~\ref{ex:bell3b}). Cut shortening makes the lower bound of \texttt{full} grow much faster than that of \texttt{frac}, which does not reach the value \valBtbUb{} of the matched bounds within the time limit. We omit \texttt{frac+short+tree}, whose curve almost coincides with that of \texttt{full}.

\begin{table}[t]
\centering
\scriptsize
\setlength{\tabcolsep}{3pt}
\caption{Paired comparisons of the ablation (Section~\ref{sec:abl}), configuration A against configuration B, on the instance-seed pairs where both runs ended. W--T--L: pairs where the lower bound of A is larger, equal, smaller, for the final lower bound (after the replay of Section~\ref{sec:semantics}) and for the bound at the end of the search (the dual bound of the tree; for \texttt{multitree}, the bound of its last master solved to optimality; the final bound for a run that records no bound during the search). $\Delta$LB: median over the instance-seed pairs of the final lower bound of A minus that of B (Appendix~\ref{app:abl} gives also the median over the instances of the medians over the seeds). Time: median over the pairs of search plus replay time, in seconds, of A and of B. Matched: pairs with matched bounds for A only and for B only; in parentheses those on \texttt{bell3binf}, where a match is only a lower bound (Example~\ref{ex:bell3b}).}
\label{tab:v2-pairs}
\begin{tabular}{@{}llccrrc@{}}
\toprule
A & B & W--T--L final & W--T--L end & $\Delta$LB & time A / B & matched A / B \\
\midrule
\multicolumn{7}{@{}l}{\emph{Gurobi}} \\
\texttt{multitree} & \texttt{int} & 40--0--10 & 45--1--4 & 42.5 & 599 / 759 & 5 / 0 \\
\texttt{frac} & \texttt{int} & 45--2--3 & 46--2--2 & 51 & 795 / 759 & 5 / 0 \\
\texttt{frac} & \texttt{multitree} & 35--0--15 & 31--4--15 & 4 & 795 / 599 & 5 / 5 \\
\texttt{frac+short} & \texttt{frac} & 41--6--3 & 38--9--3 & 6 & 660 / 795 & 11 (5) / 0 \\
\texttt{frac+tree} & \texttt{frac} & 17--28--5 & 11--33--6 & 0 & 808 / 795 & 0 / 0 \\
\texttt{frac+short+tree} & \texttt{frac+short} & 11--28--11 & 6--34--10 & 0 & 653 / 660 & 0 / 0 \\
\texttt{full} & \texttt{multitree} & 45--5--0 & 45--5--0 & 11 & 693 / 599 & 13 (5) / 0 \\
\midrule
\multicolumn{7}{@{}l}{\emph{SCIP}} \\
\texttt{multitree} & \texttt{int} & 35--0--13 & 38--0--10 & 10 & 599 / 813 & 0 / 0 \\
\texttt{frac} & \texttt{int} & 46--0--2 & 48--0--0 & 19 & 982 / 813 & 0 / 0 \\
\texttt{frac} & \texttt{multitree} & 43--0--6 & 43--0--6 & 17 & 994 / 599 & 0 / 0 \\
\texttt{frac+short} & \texttt{frac} & 47--2--0 & 46--3--0 & 13 & 649 / 994 & 9 (1) / 0 \\
\texttt{frac+tree} & \texttt{frac} & 28--12--8 & 32--11--5 & 1.5 & 983 / 982 & 4 / 0 \\
\texttt{frac+short+tree} & \texttt{frac+short} & 18--26--4 & 21--22--5 & 0 & 669 / 651 & 1 / 0 \\
\texttt{full} & \texttt{multitree} & 50--0--0 & 50--0--0 & 38 & 657 / 599 & 10 (3) / 0 \\
\bottomrule
\end{tabular}
\end{table}

The configuration \texttt{int} matches the bounds on no instance-seed pair, and \texttt{multitree} only on \texttt{mod008inf} with Gurobi, in every seed. With SCIP, it matches none. In our runs, \texttt{multitree} has an upper bound only on the pairs where it matches the bounds. These are the two configurations without fractional separation, and they differ in the search scheme and in how each cut is produced. The multi-tree loop is clearly the better one of the two: its final lower bound is larger than that of \texttt{int} in \valPairGMtIntW{} of the \valPairGMtIntN{} pairs with Gurobi and in \valPairSMtIntW{} of \valPairSMtIntN{} with SCIP. In other words, moving the separation inside the branch-and-cut does not help by itself, at least when the cuts are taken as they come from the oracle, while the multi-tree loop computes a minimum-cardinality correction at every iteration. Fractional separation alone gives a larger lower bound still: on most pairs, \texttt{frac} beats both \texttt{int} and \texttt{multitree}, with either solver, although \texttt{multitree} stays ahead on a few instances (Appendix~\ref{app:abl}).

Since the seeds of an instance are not independent, we also take the instance as unit, i.e., we count instances instead of instance-seed pairs: for each instance we take the median over its seeds of the difference between the two lower bounds, and we count the instances where this median is positive or negative. Here and below, in a comparison of a configuration A against a configuration B the difference is the lower bound of A minus that of B, so that a positive median means that A is ahead, and we write the counts of the comparison as three numbers W--T--L (win, tie, loss): the pairs, or the instances, where A is better than B, equal to B, and worse than B, i.e., with a larger lower bound or a smaller subsystem. Against \texttt{multitree}, the configuration \texttt{full} has a larger final lower bound in \valPairGFullMtW{} of the \valPairGFullMtN{} pairs with Gurobi and in \valPairSFullMtW{} of \valPairSFullMtN{} with SCIP. With the instance as unit, the median is positive, i.e., \texttt{full} is ahead, on \valInstGFullMtW{} of the \valAblNInst{} instances with Gurobi and on \valInstSFullMtW{} with SCIP, and it is never negative. This advantage is already there with \texttt{frac+short}, whose counts against \texttt{multitree} on the final lower bound are \valPairGShortMtW{}--\valPairGShortMtT{}--\valPairGShortMtL{} with Gurobi and \valPairSShortMtW{}--\valPairSShortMtT{}--\valPairSShortMtL{} of \valPairSShortMtN{} with SCIP: it comes from fractional separation and cut shortening, not from the components added later.

Cut shortening is what makes the lower bound reach the upper bound. Against \texttt{frac}, \texttt{frac+short} improves the lower bound in \valPairGShortW{} of the \valPairGShortN{} pairs with Gurobi and in \valPairSShortW{} of \valPairSShortN{} with SCIP, and the matched pairs go from \valAblGClFrac{} to \valAblGClFracShort{} with Gurobi and from \valAblSClFrac{} to \valAblSClFracShort{} with SCIP. \valInstShortSentence{} The shorter cuts also reduce the work of the master: the configuration \texttt{frac} collects a median of \valPtsRatioGMed{} times as many points as \texttt{frac+short} with Gurobi (\valPtsRatioSMed{} with SCIP). In-tree separation has a much smaller effect, both alone and on top of cut shortening: combined with cut shortening, it gives \valAblGClFracShortTree{} matched pairs with Gurobi against \valAblGClFracShort{} of \texttt{frac+short}, and \valAblSClFracShortTree{} against \valAblSClFracShort{} with SCIP. On the lower bound, its wins and losses against \texttt{frac+short} are balanced with Gurobi, and with SCIP it helps by much less than cut shortening itself.

All the counts of wins, ties and losses above are on the final lower bound, which is computed by the replay of Section~\ref{sec:semantics} after the search has ended, outside its time limit. The three main conclusions (\texttt{full} ahead of \texttt{multitree}, \texttt{frac+short} ahead of \texttt{frac}, and \texttt{frac} ahead of \texttt{int}) hold as well with the bound that each configuration has at the end of its 600\,s of search (Table~\ref{tab:v2-pairs}); Appendix~\ref{app:abl} gives the times of the replay.

The comparisons with \texttt{multitree} do not separate the cuts from the search scheme, because \texttt{multitree} has neither Phase~1 nor cut shortening. We therefore ran, within the ablation, a multi-tree loop with our cuts: after the same Phase~1 of \texttt{frac+short}, the master is solved from scratch at every iteration, and its solution is handled as an integer solution of the master in Section~\ref{sec:safe-bc}, i.e., with the same call to the oracle, the same subset of points $X$ and the same shortening of the first cut. With the instance as unit, the counts of \texttt{frac+short} against this loop are \valBoneTreeFinGW{}--\valBoneTreeFinGT{}--\valBoneTreeFinGL{} with Gurobi and \valBoneTreeFinSW{}--\valBoneTreeFinST{}--\valBoneTreeFinSL{} with SCIP on the final lower bound, but \valBoneTreeEndGW{}--\valBoneTreeEndGT{}--\valBoneTreeEndGL{} and \valBoneTreeEndSW{}--\valBoneTreeEndST{}--\valBoneTreeEndSL{} on the bound at the end of the search, which for the loop is the value of the last master solved to optimality, as for \texttt{multitree}, i.e., within the time limit neither scheme is ahead of the other. The difference comes from the replay, which often raises the final bound of the branch-and-cut above the bound of its tree (Appendix~\ref{app:abl}). On the other hand, the multi-tree loop with our cuts is ahead of \texttt{multitree} on \valBoneCutsFinGW{} of the \valBoneN{} instances with Gurobi and on \valBoneCutsFinSW{} with SCIP, and it is never behind, whichever of the two bounds is taken. Most of the distance between \texttt{multitree} and our method is therefore due to the cuts, i.e., to Phase~1 and to cut shortening, and not to the branch-and-cut.

What the branch-and-cut adds is of a different kind. It matches the bounds in \valBoneTreeClosedG{} of its \valBoneTreeRunsG{} runs with Gurobi and in \valBoneTreeClosedS{} of \valBoneTreeRunsS{} with SCIP, against \valBoneMtClosedG{} and \valBoneMtClosedS{} for the multi-tree loop with our cuts, and it ends with a verified subsystem in \valBoneTreeSubG{} and \valBoneTreeSubS{} runs, whereas the multi-tree loop has one only when the bounds match.

The objective function~\eqref{eq:sep-mip} only measures the violation of the cut, not its length. We tested three alternatives from the literature on cut separation, which also take into account the number of violated atoms \citep{FischettiLodi2007,AndreelloCapraraFischetti2007,FischettiSalvagninZanette2010}. Without cut shortening, they improve the final lower bound: against the violation, their wins minus losses in \texttt{frac+tree}, on the instance-seed pairs where both runs ended (at most \valSepObjPairs{}), are between \valSepObjNoShortGMin{} and \valSepObjNoShortGMax{} with Gurobi, and between \valSepObjNoShortSMin{} and \valSepObjNoShortSMax{} with SCIP. With cut shortening, i.e., in \texttt{frac+short+tree}, they do not: their losses minus wins are between \valSepObjShortGMin{} and \valSepObjShortGMax{} with Gurobi, and between \valSepObjShortSMin{} and \valSepObjShortSMax{} with SCIP. A better objective thus replaces cut shortening in part, but does not add to it, and we kept the violation; Appendix~\ref{app:abl} describes the three objectives.

The last two configurations add the primal components one at a time. With Gurobi, the initial subsystem of \texttt{+init} is what gives an upper bound when the bounds do not match: \texttt{+init} ends with a verified subsystem in \valAblGUbInit{} of the \valAblGDen{} runs, against \valAblGUbFracShortTree{} for \texttt{frac+short+tree}. With SCIP, the search finds a verified subsystem by itself in most runs, already with integer separation only (\valAblSUbInt{} of the \valAblSDen{} runs of \texttt{int}), and the initial subsystem adds little (\valAblSUbInit{} runs against \valAblSUbFracShortTree{}). The hard-fixing matheuristics of \texttt{full} add little on top of \texttt{+init}: with Gurobi, the upper bound of \texttt{full} is smaller than that of \texttt{+init} in \valPairGHfixUbW{} pairs and larger in \valPairGHfixUbL{} of the \valPairGHfixUbN{} pairs where both have one, and with SCIP we read no effect, given the noise of SCIP. On the lower bound, the primal components are not free: with SCIP the final lower bound of \texttt{+init} is smaller than that of \texttt{frac+short+tree} in \valPairSInitL{} pairs and larger in \valPairSInitW{} of \valPairSInitN{}, and with Gurobi that of \texttt{full} is smaller than that of \texttt{+init} in \valPairGHfixL{} pairs and larger in \valPairGHfixW{}. Replacing \texttt{HardFix(W)} by local branching \citep{FischettiLodi2003} or by proximity search \citep{FischettiMonaci2014}, or dropping it, makes no clear difference either. Appendix~\ref{app:abl} gives the counts of these comparisons.

The ablation thus points to two uses of the method. When only the lower bound is wanted, \texttt{frac+short}, with or without in-tree separation, is enough: it has the two components that explain our lower bounds, and none of the primal components, which take time from the search and on some pairs give a smaller lower bound. When a subsystem is wanted as well, as in the competition, the configuration to run is \texttt{full}, whose primal components give a verified subsystem in all its \valAblGUbFull{} runs of the ablation with Gurobi and in \valAblSUbFull{} with SCIP, against \valAblGUbFracShort{} and \valAblSUbFracShort{} for \texttt{frac+short}. This is why we use \texttt{full} for the benchmark of Section~\ref{sec:bench30}.

\subsection{The ablation on the other benchmark instances}\label{sec:abl30}

The \valAblNInst{} instances of the ablation were chosen during the development of the code. To see how much the conclusions of Section~\ref{sec:abl} depend on this choice, we ran the same eight configurations, and the multi-tree loop with our cuts of that section, with its variant of Appendix~\ref{app:abl}, on the other \valAbxN{} instances of the benchmark of Section~\ref{sec:bench30}, with the same seeds, time limit, parameters and machines; the lower blocks of Tables~\ref{tab:v2-abl-G} and~\ref{tab:v2-abl-S} report the eight configurations. These instances are not new to the code, since they belong to the set on which it was developed and some of them were used in tests of single components, but they entered neither the ablation nor the choice of its instances.

The conclusions to test and the rule to read the runs were written before the runs. We test four conclusions on the final lower bound: that fractional separation raises it (\texttt{frac} against \texttt{int}), that cut shortening raises it further (\texttt{frac+short} against \texttt{frac}), that \texttt{full} is ahead of \texttt{multitree}, and that the multi-tree loop with our cuts is ahead of \texttt{multitree}. With the instance as unit, a conclusion is confirmed if, with both solvers, the median over the seeds of the difference between the two lower bounds (the configuration expected to be ahead minus the other one) is positive on more instances than those where it is negative, and on at least twice as many.

All four conclusions are confirmed. The median is positive, i.e., the first configuration of the comparison is ahead, on \valAbxFracGW{} of the \valAbxN{} instances with Gurobi and on \valAbxFracSW{} with SCIP for \texttt{frac} against \texttt{int}, on \valAbxShortGW{} and \valAbxShortSW{} for \texttt{frac+short} against \texttt{frac}, on \valAbxFullMtGW{} and \valAbxFullMtSW{} for \texttt{full} against \texttt{multitree}, and on \valAbxLoopMtGW{} and \valAbxLoopMtSW{} for the multi-tree loop with our cuts against \texttt{multitree}. It is zero on all the other instances with a pair of runs, and it is never negative. The effects are smaller than on the instances of Section~\ref{sec:abl}, since on many of these instances both configurations match the bounds in every seed. \valAbxFailed{} of the \valAbxPlanned{} runs failed, all with SCIP, and the pairs with a failed run are not counted. On \texttt{neos-4344877-stour} this removes the instance from the three comparisons whose configuration under test uses cut shortening; if it is counted there as an instance with a negative median, the four conclusions are still confirmed. The other comparisons of Section~\ref{sec:abl} behave similarly, with smaller margins; in particular, \texttt{frac+short} and the multi-tree loop with our cuts are again even at the end of the search. Appendix~\ref{app:abl} gives the sizes of the effects, the failed runs, and the same rule applied to $\text{LB}_{\mathrm{exact}}$, which confirms the four conclusions again.

\subsection{Fixed-domain mode and preprocessing}\label{sec:form-fixdom}

In the fixed-domain mode only the linear constraints can be relaxed, while bounds and integrality are always kept. This is the setting of the X-MILP framework of \citet{LeraLeri2025}. The problem is much smaller in this mode, and it is the natural ground for the multi-tree loop, which we run here as described in Section~\ref{sec:setup}, with the same replay as our code. Table~\ref{tab:v2-fd} compares it with our method, i.e., \texttt{full} run in the fixed-domain mode, on \valFdN{} instances. \valFdMatchSentence{} For the times, we take the median search time over the seeds, and the shifted geometric mean (shift 1\,s) over the instances with matched bounds in every seed for both methods. With Gurobi, this mean is \valFdTimeGMt\,s for the multi-tree loop and \valFdTimeGOurs\,s for our method on \valFdTimeGN{} instances, and the multi-tree loop is faster on \valFdTimeGMtFaster{} of them. With SCIP, the means are \valFdTimeSMt\,s and \valFdTimeSOurs\,s on \valFdTimeSN{} instances, and the multi-tree loop is faster on \valFdTimeSMtFaster{}. These times should be read with the length of our Phase~1 in mind, which depends on the time limit. In any case, they show no advantage of our method. In other words, in the fixed-domain mode our implementation of the multi-tree loop is as good as our method with Gurobi, and faster. With SCIP it takes about the same time, and matches the bounds somewhat less often.

The picture changes when bounds and integrality can be relaxed, and not because of a different choice of instances. On the \valComN{} instances that are both in this experiment and in the ablation, with the same time limit, the multi-tree loop matches the bounds in \valComFdMtG{} of the \valComDen{} pairs with Gurobi and in \valComFdMtS{} with SCIP in the fixed-domain mode, against \valComFdOursG{} and \valComFdOursS{} for our method. In the removable-domain mode, it matches them in \valComRemMtG{} and \valComRemMtS{} pairs, against \valComRemOursG{} and \valComRemOursS{} for \texttt{full}. The removable-domain problem is thus harder for both methods, and it is where the two methods separate. The last column of the table shows what this mode gives in terms of cardinality: relaxing bounds and integrality gives much smaller subsystems, by an order of magnitude on \texttt{neos-942886} and \texttt{p0201inf}, mostly because the fixed-domain mode always keeps all the domain atoms, whose number is the constant $\kappa$ of Table~\ref{tab:v2-fd}. The two diagnoses are thus of a different kind: the fixed-domain one names few rows and leaves every domain in the conflict, while the removable-domain one also says which domains take part in it.

The same experiment measures the effect of Observation~\ref{obs:dom}. On these instances the dominance argument removes between \valRedMinPctTwo\% and \valRedMaxPctTwo\% of the atoms. The configuration \texttt{full} matches the bounds in \valPreGRed{} of the \valPreGDen{} pairs with Gurobi when run on the reduced atom set, and in \valPreGAll{} when run on all atoms. With SCIP, the counts are \valPreSRed{} and \valPreSAll{}.\valPreAllCountSentence{}

\begin{table}[t]
\centering
\scriptsize
\setlength{\tabcolsep}{3pt}
\caption{Fixed-domain mode (only rows can be relaxed) on 12 instances, 5 seeds (600\,s of search). Values are cardinalities of the whole subsystem, i.e., kept rows plus the number $\kappa$ of always-kept domain atoms, so that the upper bounds of this table are upper bounds on the removable-domain optimum of the last column. Each cell gives the median over the seeds of the final lower bound and of the verified upper bound, [LB, UB], and in parentheses the number of seeds in which the two bounds match; a bold single value means that every seed matched. Our method is \texttt{full} in fixed-domain mode; last column: \texttt{full} in removable-domain mode with Gurobi. $^{\ast n}$ = $n$ runs stopped on an undecided oracle call; $^{\mathrm{r}}$ = an instance on which subsystems accepted by the checker can be feasible in exact arithmetic (Example~\ref{ex:bell3b}); the fixed-domain subsystems of this table are not among those we checked.}
\label{tab:v2-fd}
\setlength{\tabcolsep}{2pt}
\begin{tabular}{@{}lrrrrrr@{}}
\toprule
 & & \multicolumn{2}{c}{multi-tree} & \multicolumn{2}{c}{our method} & removable \\
\cmidrule(lr){3-4}\cmidrule(lr){5-6}
Instance & $\kappa$ & Gurobi & SCIP & Gurobi & SCIP & domains \\
\midrule
\texttt{aim-200-1\_6-no-2} & 200 & \textbf{280} (5) & \textbf{280} (5) & \textbf{280} (5) & \textbf{280} (5) & \textbf{101} (5) \\
\texttt{bell3ainf} & 133 & \textbf{191} (5) & \textbf{191} (5) & \textbf{191} (5) & \textbf{191} (5) & \textbf{117} (5) \\
\texttt{bell3binf}$^{\mathrm{r}}$ & 133 & \textbf{195} (5) & \textbf{195} (5) & \textbf{195} (5) & \textbf{195} (5) & \textbf{118} (5) \\
\texttt{bell5inf} & 104 & \textbf{143} (5) & \textbf{143} (5) & \textbf{143} (5) & \textbf{143} (5) & \textbf{94} (5) \\
\texttt{g503inf} & 48 & \textbf{60} (5) & \textbf{60} (5) & \textbf{60} (5) & \textbf{60} (5) & [30, 30] (3) \\
\texttt{misc05inf} & 130 & \textbf{177} (5) & [174, --] (0) & \textbf{177} (5) & \textbf{177} (5) & [94, 128] (0)$^{\ast2}$ \\
\texttt{mod008inf} & 319 & \textbf{324} (5) & \textbf{324} (5) & \textbf{324} (5) & \textbf{324} (5) & \textbf{323} (5) \\
\texttt{neos-942886} & 480 & \textbf{497} (5) & [497, 497] (3) & \textbf{497} (5) & \textbf{497} (5) & [32, 34] (0) \\
\texttt{neos859080} & 160 & [294, --] (0) & [293, --] (0) & [296, 317] (0) & [294, 319] (0)$^{\ast1}$ & [176, 225] (0)$^{\ast3}$ \\
\texttt{grid4\_31} & 124 & \textbf{218} (5) & \textbf{218} (5) & [218, 218] (4) & \textbf{218} (5) & [92, 95] (0) \\
\texttt{vertexcover} & 190 & \textbf{381} (5) & \textbf{381} (5) & \textbf{381} (5) & [381, 381] (4) & [197, 200] (0) \\
\texttt{p0201inf} & 201 & \textbf{203} (5) & \textbf{203} (5) & \textbf{203} (5) & \textbf{203} (5) & \textbf{9} (5) \\
\midrule
matched & & 55/60 & 48/60 & 54/60 & 54/60 & \\
\bottomrule
\end{tabular}
\end{table}

\FloatBarrier
\subsection{Instances outside the benchmark}\label{sec:oob}

All the results above are on instances of the set on which the code was developed. As a check, we ran the same code, without any change, on \valOobN{} infeasible instances with at most 1{,}000 atoms that are not in the public set: \valOobNH{} infeasible instances of MIPLIB~2017 \citep{GleixnerEtAl2021}, and \valOobNG{} feasible instances of the same library made infeasible by a constraint on the value of the objective function (Appendix~\ref{app:details} describes them, and Table~\ref{tab:v2-oob} gives the results). We ran \texttt{multitree}, \texttt{frac}, \texttt{frac+short} and \texttt{full} as in Section~\ref{sec:abl}, on the second partition of the cluster, and \texttt{full} with the 1{,}200\,s of Section~\ref{sec:bench30}, and the criteria to read the experiment were written before the runs. Most of these instances had appeared in earlier experiments made with previous versions of the code, so this is a check outside the benchmark rather than a test on unseen families. The two conclusions of the ablation on the lower bound that we had chosen to test are confirmed. With the instance as unit, the lower bound of \texttt{full} is larger than that of \texttt{multitree} on \valOobLbGW{} of the \valOobN{} instances with Gurobi and on \valOobLbSW{} with SCIP, that of \texttt{frac+short} is larger than that of \texttt{frac} on \valOobShGW{} and on \valOobShSW{} instances, respectively, and neither is ever smaller.

The upper bounds, instead, are weaker than on the benchmark: lower and upper bounds match on two instances only, \valOobNoSub{} of the \valOobPairs{} runs of \texttt{full} with 1{,}200\,s that ended normally have no verified subsystem, and by the criterion fixed before the runs, which is a descriptive rule, our method is equivalent as to the size of the subsystem to the deletion filter, both in a single pass and in a multi-start version (Appendix~\ref{app:details} gives the rule, the multi-start filter and the counts). What our method adds on these instances is the lower bound. Without cut shortening, on the other hand, \texttt{frac} is not ahead of \texttt{multitree} on these instances with Gurobi, while it is with SCIP.

\subsection{Validity of the bounds}\label{sec:validity}

All the comparisons above rely on the validity of the bounds, which we check here; Appendix~\ref{app:validity} gives the details. In the \valAllRuns{} runs of our code in Sections~\ref{sec:bench30}--\ref{sec:oob}, no lower bound exceeds the upper bound of the same run, and every reported upper bound was accepted by the official checker. In the \valCrossGroups{} groups of runs of Sections~\ref{sec:bench30}--\ref{sec:form-fixdom} on the same instance and in the same mode (removable or fixed domains), all runs with matched bounds report the same value. These checks are necessary conditions for valid bounds, not sufficient ones. Appendix~\ref{app:validity} also compares the two readings of the verdict of the checker of Section~\ref{sec:semantics}, and solves again every reported subsystem after relaxing it by the tolerance of the checker. Beyond these checks, for one run on each of the \valCertTot{} instances of the benchmark where \texttt{full} has a run with $\text{LB}_{\mathrm{exact}} = \text{UB}$, we attempted to certify both bounds in exact arithmetic with the binary64 data, independently of the checker and of the two solvers. Both are certified for \valCertN{} of these runs (Table~\ref{tab:cert} in Appendix~\ref{app:validity}), so that on these instances $z^*_{\mathrm{exact}} = \text{UB}$, up to the correctness of our own checker of the certificates, which uses rational arithmetic only. On \texttt{flugplinf}, which is not among them, Example~\ref{ex:flugpl} proves $z^*_{\mathrm{exact}} = \valFlugLbExact$ directly. On the other ones, one part reached its node or time limit, except on \texttt{bell3binf}, where the subsystem accepted by the official checker has a point in exact arithmetic (Example~\ref{ex:bell3b}), so that there the matched value is only a lower bound.

Of the \valPlannedRuns{} planned runs of our code, \valFailTotal{} did not end normally, all of them with the SCIP backend, and the search was stopped on an undecided candidate (Section~\ref{sec:safe-bc}) in \valStoppedRuns{} of the \valSolverRuns{} runs of the benchmark, of the ablation and of the fixed-domain comparison, in \valAbxStopped{} of the \valAbxRuns{} runs of Section~\ref{sec:abl30}, and in \valOobStopped{} of the \valOobPairs{} runs of \texttt{full} with 1{,}200\,s of Section~\ref{sec:oob}. A stopped run still reports valid bounds, in the sense of Section~\ref{sec:semantics}. Table~\ref{tab:campaigns} in Appendix~\ref{app:validity} gives these counts by experiment.

\section{Conclusions}\label{sec:concl}

We have proposed a branch-and-cut algorithm for the minimum-cardinality infeasible subsystem problem on general MILPs, where linear constraints, variable bounds and integrality requirements can all be relaxed. The algorithm works on an implicit hitting-set master: a dominance argument first removes the bound atoms of integer variables, the covering LP is solved by fractional separation with cut shortening, and the integer master is solved by branch-and-cut, starting from an initial subsystem, with in-tree separation and two hard-fixing matheuristics. The lower bound is obtained by a replay, in which every stored point is evaluated again on the original data, in exact arithmetic and within a tolerance of $10^{-9}$, and the two covering problems are solved again. The upper bound is the cardinality of a subsystem accepted by the official checker of the competition. What the two bounds certify is stated in Section~\ref{sec:semantics}. Our aim has been to compute lower bounds that hold for every infeasible subsystem ($\text{LB}_{\mathrm{exact}}$ for infeasibility in exact arithmetic, LB for infeasibility within the tolerance $\delta$) and, whenever possible, to match them with the cardinality of a verified subsystem: the primal components have the role that primal heuristics have in any branch-and-cut, and we have not addressed the purely heuristic context, where a small subsystem is wanted and no bound on its quality is required.

The computational study, made with two solvers and \valSeeds{} seeds per configuration, has shown that the components of the algorithm do not contribute in the same way. Branch-and-cut with integer separation only has given worse lower bounds than our implementation of the classical loop, which solves the master from scratch at every iteration. Fractional separation makes the lower bound move, cut shortening makes it reach the upper bound, and in-tree separation adds little on top of them. The effects of fractional separation and of cut shortening on the lower bound, and the advantage over the classical loop, are confirmed on the other \valAbxN{} instances of the benchmark. When the classical loop is given the same Phase~1 and the same cut shortening, its lower bounds improve on almost all the instances of the ablation. The corresponding branch-and-cut is then ahead on the final lower bound, which is computed by the replay, but not on the bound that the two schemes have at the end of the search. Our lower bounds are thus explained by the cuts much more than by the search scheme. The branch-and-cut has however matched the bounds in more runs with SCIP (in about as many with Gurobi) and, unlike the loop, it can end with a verified subsystem when the bounds do not match. In the ablation, almost every run of the configurations with the initial subsystem ends with a verified subsystem, which with Gurobi is due to the initial subsystem, while the hard-fixing matheuristics make it smaller in about as many runs as they make it larger. For the lower bound alone, \texttt{frac+short} is therefore enough, and \texttt{full} is the variant to use when a subsystem is wanted as well. On the instances of the MIPcc27 competition with at most 1{,}000 atoms, lower and upper bounds have matched on about half of them. Exact certificates have confirmed both bounds for one run on \valCertN{} of the \valCertTot{} instances where a run has $\text{LB}_{\mathrm{exact}} = \text{UB}$, and they have shown on \texttt{bell3binf} subsystems accepted by the checker that are feasible in exact arithmetic: how the verdict of the checker is read is not a formality. Where the bounds have matched (\texttt{bell3binf} excepted), the subsystems found have mostly kept the integer domain of less than half of the integer variables, and on several instances a single variable has to be integer once the bounds are kept. This is the kind of diagnosis that the setting of the competition makes possible. In the classical setting where only the linear constraints can be relaxed and the bounds and integrality requirements of the variables stay fixed (the fixed-domain setting of Section~\ref{sec:form-fixdom}), our implementation of the multi-tree loop has matched the bounds as often as our method when Gurobi is used, in less time, and somewhat less often when SCIP is used. When bounds and integrality can be relaxed too, the problem is harder for both, and our method has given a larger lower bound than the multi-tree loop in almost all paired runs.

The limits of the method are also clear from our runs. We did not address the instances of the public set with more than 1{,}000 atoms, and outside the set on which the code was developed we only have a check on \valOobN{} instances, most of which were not new to earlier versions of the code. There, the advantage over the multi-tree loop and the effect of cut shortening on the lower bound have been confirmed, but our method has ended without a verified subsystem in more than a quarter of its runs, and its subsystems have turned out to be neither smaller nor larger than those of a deletion filter. When the bounds do not match, the subsystem in hand comes with a gap, and a clever heuristic could give a better upper bound in the same time. Moreover, a number of runs are stopped on a candidate that cannot be decided, the SCIP backend of our code ended with a segmentation fault on two instances, and with SCIP some runs on the instances of the \texttt{enlight} family exceeded the memory of the blade. Natural directions for future work are a Phase-1 exit rule that does not depend on the time limit, a way to continue the search after an undecided oracle call without compromising the bounds, and stronger valid inequalities for the covering master.

\FloatBarrier
\appendix

\section{Parameters of the algorithm}\label{app:params}

\begin{table}[t]
\centering
\scriptsize
\caption{Parameters of the algorithm, as set in the code used for all experiments. Where the two backends differ, both values are given (Gurobi / SCIP). $\mathrm{rem}$ = remaining time.}
\label{tab:params}
\begin{tabular}{@{}p{1.9cm}p{2.5cm}p{6.6cm}@{}}
\toprule
Component & Parameter & Value \\
\midrule
Phase~1 & time cap; smoothing & 45\% of the time limit (12\% without fractional separation); $\alpha = 0.6$, restored every 15 iterations at $\alpha = 1$ \\
 & stop on plateau & 35 iterations without $+10^{-3}$, after at least 45\,s \\
Separation MILP~\eqref{eq:sep-mip}--\eqref{eq:sep-dom} & hard set; threshold & $\bar z_a \ge 0.9995$; 0.999 (Phase~1, matheuristics), 0.995 (in tree) \\
 & solutions per call & 11 (Phase~1), 5 (in tree), 6 (matheuristics) \\
 & time limit per call & 8\,s if $\card{S} \le 150$, else 14\,s (Phase~1); 2.0\,s / 2.5--1.5\,s (in tree) \\
Cut shortening & activation; blocks & $\card{V} \ge 6$; blocks of $\min(8, \lfloor \card{V}/4 \rfloor)$ when $\card{V} \ge 12$ \\
 & time per check; per point & 0.12--0.18\,s; 0.30--0.65\,s \\
 & exact repair & one per point; projection only if at most 200 active rows; tolerance witness with $\delta = 10^{-9}$ if it fails \\
Phase~2 oracle & feasibility probe & if $\card{\redatoms} > 350$ or $\card{\redatoms\setminus C} > 150$; $\min(8, \mathrm{rem}-2)$\,s, retry $\min(20, 0.15\,\mathrm{rem})$\,s \\
 & repeated check & once, $\bar T = 300$\,s, whatever the remaining time; same limit for the checker and for the checks with tolerances $10^{-9}$ \\
In-tree separation & nodes; caps & root, every 200 nodes (SCIP: also depth $\le 3$ among the first 15); 20 calls, 12\% of Phase~2 \\
\texttt{HardFix(U)} & budget & $\min(25, \max(6, 0.09\,\mathrm{rem}))$\,s / $\min(28, \max(6, 0.09\,\mathrm{rem}))$\,s \\
\texttt{HardFix(W)} & trigger; size; budget & 45\,s without improvement and 55\,s since last call; $\min(80,\max(25,\card{U}/2))$ atoms added to $U$; $\min(25, \max(8, 0.08\,\mathrm{rem}))$\,s \\
Final bounds & replay; checker & supports at tolerance $\delta$ and at tolerance 0, covering solves within the replay time limit, RC2 fallback; tolerance $10^{-5}$ \\
\bottomrule
\end{tabular}
\end{table}

\FloatBarrier
\section{Detailed results of the ablation}\label{app:abl}

Tables~\ref{tab:v2-abl-G} and~\ref{tab:v2-abl-S} give the bounds of every configuration of Section~\ref{sec:abl} on every instance of the benchmark, and Table~\ref{tab:v2-replay} the time of the replay of the lower bound; the paired comparisons drawn from them are in Table~\ref{tab:v2-pairs}.

\begin{table}[t]
\centering
\scriptsize
\setlength{\tabcolsep}{3pt}
\caption{Ablation with Gurobi: in the upper block, the 10 instances of the ablation subset (Section~\ref{sec:setup}); in the lower block, the other 20 instances of the benchmark (Section~\ref{sec:abl30}). 5 seeds per configuration (600\,s of search, four 1-thread runs per exclusive blade). Each cell gives the median over the seeds of the final lower bound and of the verified upper bound, [LB, UB], and in parentheses the number of seeds in which the two bounds match; a bold single value means that every seed matched. A dash means no verified upper bound in most seeds; an upper bound equal to the number of atoms (the whole model) is not counted. Configurations as in Section~\ref{sec:setup}. $^{\dagger n}$ = $n$ runs failed (counted as not matched); $^{\ast n}$ = $n$ runs stopped on an undecided oracle call. Last rows of each block: matched instance-seed pairs; runs stopped on an undecided oracle call; runs whose search ended more than 5\,s after the time limit (Section~\ref{sec:safe-bc}); failed runs.}
\label{tab:v2-abl-G}
\resizebox{\textwidth}{!}{%
\begin{tabular}{@{}lrrrrrrrr@{}}
\toprule
Instance & \texttt{multitree} & \texttt{int} & \texttt{frac} & \texttt{frac+short} & \texttt{frac+tree} & \texttt{frac+short+tree} & \texttt{+init} & \texttt{full} \\
\midrule
\texttt{bell3binf} & [104, --] (0) & [14, --] (0) & [112, --] (0) & \textbf{118} (5) & [113, --] (0) & \textbf{118} (5) & \textbf{118} (5) & \textbf{118} (5) \\
\texttt{bell5inf} & [88, --] (0) & [15, --] (0) & \textbf{94} (5) & \textbf{94} (5) & \textbf{94} (5) & \textbf{94} (5) & \textbf{94} (5) & \textbf{94} (5) \\
\texttt{enlight11} & [31, --] (0) & [35, --] (0) & [34, --] (0) & [90, --] (0) & [34, --] (0) & [89, --] (0) & [85, 246] (0)$^{\ast2}$ & [82, 241] (0)$^{\ast5}$ \\
\texttt{g503inf} & [18, --] (0) & [24, --] (0) & [29, 30] (0) & [28, 30] (1) & [29, 30] (0) & [29, 30] (1) & [30, 30] (3) & [30, 30] (3) \\
\texttt{hole9} & [427, --] (0) & [422, --] (0) & [433, --] (0) & [434, --] (0) & [433, --] (0) & [435, --] (0) & [435, 479] (0) & [435, 480] (0) \\
\texttt{misc05inf} & [61, --] (0) & [24, --] (0) & [78, --] (0) & [95, --] (0) & [78, --] (0) & [95, --] (0)$^{\ast1}$ & [95, 128] (0) & [94, 131] (0)$^{\ast1}$ \\
\texttt{mod008inf} & \textbf{323} (5) & [27, --] (0) & [115, --] (0) & \textbf{323} (5) & [122, --] (0) & \textbf{323} (5) & \textbf{323} (5) & \textbf{323} (5) \\
\texttt{neos-ythan} & [63, --] (0) & [35, --] (0) & [66, --] (0) & [72, --] (0) & [66, --] (0) & [72, --] (0) & [72, 145] (0) & [72, 151] (0) \\
\texttt{grid4\_31} & [74, --] (0) & [24, --] (0) & [71, --] (0) & [93, 107] (0) & [72, --] (0) & [93, 114.5] (0) & [92, 95] (0) & [92, 95] (0) \\
\texttt{vertexcover} & [196, --] (0) & [93, --] (0) & [195, --] (0) & [197, 200] (0) & [195, --] (0) & [197, 200] (0) & [197, 200] (0) & [197, 200] (0) \\
\midrule
matched & 5/50 & 0/50 & 5/50 & 16/50 & 5/50 & 16/50 & 18/50 & 18/50 \\
stopped & 0 & 0 & 0 & 0 & 0 & 1 & 2 & 6 \\
beyond limit & 0 & 0 & 0 & 0 & 0 & 0 & 2 & 5 \\
failed & 0 & 0 & 0 & 0 & 0 & 0 & 0 & 0 \\
\midrule
\texttt{aim-100-1\_6-no-4} & \textbf{59} (5) & [52, --] (0) & \textbf{59} (5) & \textbf{59} (5) & \textbf{59} (5) & \textbf{59} (5) & \textbf{59} (5) & \textbf{59} (5) \\
\texttt{aim-200-1\_6-no-2} & [97, --] (0) & [67, --] (0) & [93, --] (0) & \textbf{101} (5) & [93, --] (0) & \textbf{101} (5) & \textbf{101} (5) & \textbf{101} (5) \\
\texttt{aim-50-2\_0-no-4} & \textbf{27} (5) & \textbf{27} (5) & \textbf{27} (5) & \textbf{27} (5) & \textbf{27} (5) & \textbf{27} (5) & \textbf{27} (5) & \textbf{27} (5) \\
\texttt{bell3ainf} & [114, --] (0) & [9, --] (0) & \textbf{117} (5) & \textbf{117} (5) & \textbf{117} (5) & \textbf{117} (5) & \textbf{117} (5) & \textbf{117} (5) \\
\texttt{enlight4} & [16, --] (0) & \textbf{18} (5) & \textbf{18} (5) & \textbf{18} (5) & \textbf{18} (5) & [18, 18] (4)$^{\ast1}$ & \textbf{18} (5) & \textbf{18} (5) \\
\texttt{flugplinf} & \textbf{23} (5) & \textbf{23} (5) & \textbf{23} (5) & \textbf{23} (5) & \textbf{23} (5) & \textbf{23} (5) & \textbf{23} (5) & \textbf{23} (5) \\
\texttt{hole10} & [575, --] (0) & [548, --] (0) & [581, --] (0) & [582, --] (0) & [581, --] (0) & [582, --] (0) & [582, 644] (0)$^{\ast1}$ & [583, 644] (0)$^{\ast1}$ \\
\texttt{jnh10} & [30, --] (0) & [30, --] (0) & [50, --] (0) & [55, --] (0) & [50, --] (0) & [55, --] (0) & [55, 85] (0) & [54, 96] (0) \\
\texttt{jnh15} & [28, --] (0) & [28, --] (0) & [46, --] (0) & [51, --] (0) & [46, --] (0) & [52, --] (0) & [51, 154] (0) & [51, 156] (0) \\
\texttt{jnh215} & [33, --] (0) & [32, --] (0) & [51, --] (0) & [59, --] (0) & [52, --] (0) & [59, --] (0) & [59, 116] (0) & [59, 108] (0) \\
\texttt{mathologer4a} & \textbf{4} (5) & \textbf{4} (5) & \textbf{4} (5) & \textbf{4} (5) & \textbf{4} (5) & \textbf{4} (5) & \textbf{4} (5) & \textbf{4} (5) \\
\texttt{neos-4331170-senio} & \textbf{12} (5) & \textbf{12} (5) & \textbf{12} (5) & \textbf{12} (5) & \textbf{12} (5) & \textbf{12} (5) & \textbf{12} (5) & \textbf{12} (5) \\
\texttt{neos-4344877-stour} & [50, --] (0) & [50, --] (0) & [50, --] (0) & [66, --] (0) & [50, --] (0) & [66, --] (0) & [65, 125] (0) & [64, 145] (0) \\
\texttt{neos-942886} & [18, --] (0) & [18, --] (0) & [20, --] (0) & [33, 33] (4) & [20, --] (0) & [33, 33] (4) & [33, 33] (3) & [33, 33] (3) \\
\texttt{neos859080} & [169, --] (0) & [83, --] (0) & [175, --] (0) & [177, --] (0) & [175, --] (0) & [177, --] (0) & [176, 228] (0)$^{\ast4}$ & [176, 257] (0)$^{\ast3}$ \\
\texttt{domset-hexgrid\_03\_052} & [150, --] (0) & [47, --] (0) & [70, --] (0) & [158, --] (0) & [71, --] (0) & [158, --] (0) & [158, 231] (0) & [158, 238] (0) \\
\texttt{p0201inf} & \textbf{9} (5) & [6, --] (0) & \textbf{9} (5) & \textbf{9} (5) & \textbf{9} (5) & \textbf{9} (5) & \textbf{9} (5) & \textbf{9} (5) \\
\texttt{stein10i} & \textbf{12} (5) & \textbf{12} (5) & \textbf{12} (5) & \textbf{12} (5) & \textbf{12} (5) & \textbf{12} (5) & \textbf{12} (5) & \textbf{12} (5) \\
\texttt{stein15inf} & \textbf{42} (5) & \textbf{42} (5) & \textbf{42} (5) & \textbf{42} (5) & \textbf{42} (5) & \textbf{42} (5) & \textbf{42} (5) & \textbf{42} (5) \\
\texttt{stein45inf} & [14, --] (0) & [14, --] (0) & [20, --] (0) & [21, --] (0) & [20, --] (0) & [21, --] (0) & [21, 252] (0) & [21, 252] (0) \\
\midrule
matched & 40/100 & 35/100 & 50/100 & 59/100 & 50/100 & 58/100 & 58/100 & 58/100 \\
stopped & 0 & 0 & 0 & 0 & 0 & 1 & 5 & 4 \\
beyond limit & 0 & 0 & 0 & 0 & 0 & 1 & 6 & 5 \\
failed & 0 & 0 & 0 & 0 & 0 & 0 & 0 & 0 \\
\bottomrule
\end{tabular}}
\end{table}
\begin{table}[t]
\centering
\scriptsize
\setlength{\tabcolsep}{3pt}
\caption{Ablation with SCIP: in the upper block, the 10 instances of the ablation subset (Section~\ref{sec:setup}); in the lower block, the other 20 instances of the benchmark (Section~\ref{sec:abl30}). 5 seeds per configuration (600\,s of search, four 1-thread runs per exclusive blade). Each cell gives the median over the seeds of the final lower bound and of the verified upper bound, [LB, UB], and in parentheses the number of seeds in which the two bounds match; a bold single value means that every seed matched. A dash means no verified upper bound in most seeds; an upper bound equal to the number of atoms (the whole model) is not counted. Configurations as in Section~\ref{sec:setup}. $^{\dagger n}$ = $n$ runs failed (counted as not matched); $^{\ast n}$ = $n$ runs stopped on an undecided oracle call. Last rows of each block: matched instance-seed pairs; runs stopped on an undecided oracle call; runs whose search ended more than 5\,s after the time limit (Section~\ref{sec:safe-bc}); failed runs.}
\label{tab:v2-abl-S}
\resizebox{\textwidth}{!}{%
\begin{tabular}{@{}lrrrrrrrr@{}}
\toprule
Instance & \texttt{multitree} & \texttt{int} & \texttt{frac} & \texttt{frac+short} & \texttt{frac+tree} & \texttt{frac+short+tree} & \texttt{+init} & \texttt{full} \\
\midrule
\texttt{bell3binf} & [88, --] (0) & [30, 233] (0) & [107, 160] (0) & [117, 126] (1) & [109, 159] (0) & [117, 119] (1) & [117, 119] (2) & [118, 118] (3) \\
\texttt{bell5inf} & [73, --] (0) & [89, 103] (0) & [91, 94] (0) & \textbf{94} (5) & [91, 94] (0) & \textbf{94} (5) & \textbf{94} (5) & \textbf{94} (5) \\
\texttt{enlight11} & [18, --] (0) & [15, --] (0)$^{\dagger2}$$^{\ast3}$ & [33.5, --] (0)$^{\dagger1}$$^{\ast1}$ & [36, --] (0)$^{\dagger1}$ & [32, --] (0)$^{\dagger2}$ & [34.5, --] (0)$^{\dagger1}$ & [38, 350] (0)$^{\dagger2}$ & [38, 341] (0) \\
\texttt{g503inf} & [14, --] (0) & [21, 34] (0) & [24, 30] (0) & [27, 30] (0) & [24, 30] (0) & [28, 30] (0) & [27, 30] (0) & [27, 30] (0) \\
\texttt{hole9} & [220, --] (0) & [20, 498] (0) & [428, 495] (0) & [431, 495] (0) & [430, 484] (0) & [431, 494] (0) & [432, 490] (0) & [432, 492] (0) \\
\texttt{misc05inf} & [40, --] (0) & [28, 370] (0) & [46, 306] (0) & [82, 259] (0) & [47, 185] (0)$^{\ast1}$ & [84, 267] (0) & [81, 123] (0) & [82, 121] (0) \\
\texttt{mod008inf} & [242, --] (0) & [36, --] (0) & [292, 324.5] (0) & [323, 323] (3) & [323, 323] (4) & [323, 323] (4) & [321, 323] (2) & [316, 323] (2) \\
\texttt{neos-ythan} & [30, --] (0) & [24, 226.5] (0) & [37, --] (0) & [69, --] (0) & [52, --] (0) & [70, 206] (0) & [69, 167] (0)$^{\dagger2}$ & [69, 167] (0) \\
\texttt{grid4\_31} & [44, --] (0) & [40, 315] (0) & [73, 140] (0) & [89, 100] (0) & [73, 145] (0) & [90, 99] (0) & [89, 104] (0) & [89, 103] (0) \\
\texttt{vertexcover} & [80, --] (0) & [38, 553] (0) & [42, 565] (0) & [194, 236] (0) & [46, 567] (0) & [194, 237] (0) & [193, 200] (0) & [193, 203] (0) \\
\midrule
matched & 0/50 & 0/50 & 0/50 & 9/50 & 4/50 & 10/50 & 9/50 & 10/50 \\
stopped & 0 & 3 & 1 & 0 & 1 & 0 & 0 & 0 \\
beyond limit & 0 & 3 & 5 & 5 & 5 & 5 & 5 & 7 \\
failed & 0 & 2 & 1 & 1 & 2 & 1 & 4 & 0 \\
\midrule
\texttt{aim-100-1\_6-no-4} & \textbf{59} (5) & [46, 122] (0) & \textbf{59} (5) & \textbf{59} (5) & \textbf{59} (5) & \textbf{59} (5) & \textbf{59} (5) & \textbf{59} (5) \\
\texttt{aim-200-1\_6-no-2} & [94, --] (0) & [45, 217] (0) & [77, 373] (0) & [101, 101] (3) & [77, 268] (0) & [101, 101] (4) & [100, 101] (1) & [100, 101] (1) \\
\texttt{aim-50-2\_0-no-4} & \textbf{27} (5) & \textbf{27} (5) & \textbf{27} (5) & \textbf{27} (5) & \textbf{27} (5) & \textbf{27} (5) & \textbf{27} (5) & \textbf{27} (5) \\
\texttt{bell3ainf} & [92, --] (0) & [41, 232] (0) & [115, 119] (1) & \textbf{117} (5) & [116, 117] (1) & \textbf{117} (5) & \textbf{117} (5) & \textbf{117} (5) \\
\texttt{enlight4} & [15, --] (0) & [18, 18] (3)$^{\ast2}$ & [18, 18] (3)$^{\dagger1}$$^{\ast1}$ & \textbf{18} (5) & [18, 34.5] (1)$^{\dagger1}$$^{\ast3}$ & [18, 18] (4)$^{\ast1}$ & \textbf{18} (5) & \textbf{18} (5) \\
\texttt{flugplinf} & \textbf{23} (5) & \textbf{23} (5) & \textbf{23} (5) & \textbf{23} (5) & \textbf{23} (5) & \textbf{23} (5) & \textbf{23} (5) & \textbf{23} (5) \\
\texttt{hole10} & [15, --] (0) & [17, 667] (0)$^{\ast5}$ & [291, --] (0) & [579, --] (0)$^{\ast3}$ & [443, --] (0) & [579, --] (0)$^{\ast5}$ & [579, 659] (0)$^{\ast5}$ & [579, 659] (0)$^{\ast5}$ \\
\texttt{jnh10} & [19, --] (0) & [19, 385] (0) & [30, 464] (0) & [46, 573] (0) & [31, 273] (0) & [49, 199] (0) & [45, 88] (0) & [45, 91] (0) \\
\texttt{jnh15} & [16, --] (0) & [16, 415] (0) & [17, 382] (0) & [41, 470] (0) & [17, 380] (0) & [43, 945] (0) & [39, 163] (0) & [41, 173] (0) \\
\texttt{jnh215} & [19, --] (0) & [19, 382] (0) & [19, 513] (0) & [49, 536] (0) & [19, 475] (0) & [53, 221] (0) & [48, 115] (0) & [48, 122] (0) \\
\texttt{mathologer4a} & \textbf{4} (5) & \textbf{4} (5) & \textbf{4} (5) & \textbf{4} (5) & \textbf{4} (5) & \textbf{4} (5) & \textbf{4} (5) & \textbf{4} (5) \\
\texttt{neos-4331170-senio} & \textbf{12} (5) & \textbf{12} (5) & \textbf{12} (5) & \textbf{12} (5) & \textbf{12} (5) & \textbf{12} (5) & \textbf{12} (5) & \textbf{12} (5) \\
\texttt{neos-4344877-stour} & [3, --] (0) & [7, 824] (0) & [21, 891] (0) & --$^{\dagger5}$ & [21, 891] (0) & --$^{\dagger5}$ & --$^{\dagger5}$ & --$^{\dagger5}$ \\
\texttt{neos-942886} & [9, --] (0) & [12, 677] (0) & [13, 753] (0) & [24, 741] (0) & [12, 754] (0) & [26, 825] (0) & [24, 260] (0) & [23, 41] (0) \\
\texttt{neos859080} & [108, --] (0) & [91, 315] (0)$^{\ast2}$ & [171, 302] (0) & [174, 321] (0)$^{\ast2}$ & [172, 310.5] (0) & [174, 321] (0)$^{\ast4}$ & [174, 285] (0)$^{\ast2}$ & [174, 259] (0)$^{\ast3}$ \\
\texttt{domset-hexgrid\_03\_052} & [18, --] (0) & [9, 368] (0) & [18, 423] (0) & [158, 314] (0) & [18, 423] (0) & [158, 314] (0) & [158, 231] (0) & [158, 306] (0) \\
\texttt{p0201inf} & \textbf{9} (5) & [9, 9] (3) & \textbf{9} (5) & \textbf{9} (5) & \textbf{9} (5) & \textbf{9} (5) & \textbf{9} (5) & \textbf{9} (5) \\
\texttt{stein10i} & \textbf{12} (5) & \textbf{12} (5) & \textbf{12} (5) & \textbf{12} (5) & \textbf{12} (5) & \textbf{12} (5) & \textbf{12} (5) & \textbf{12} (5) \\
\texttt{stein15inf} & \textbf{42} (5) & \textbf{42} (5) & \textbf{42} (5) & \textbf{42} (5) & \textbf{42} (5) & \textbf{42} (5) & \textbf{42} (5) & \textbf{42} (5) \\
\texttt{stein45inf} & [12, --] (0) & [11, 262] (0) & [19, 263] (0) & [19, 264] (0) & [20, 233] (0) & [20, 233] (0) & [20, 306] (0) & [20, 306] (0) \\
\midrule
matched & 40/100 & 36/100 & 44/100 & 53/100 & 42/100 & 53/100 & 51/100 & 51/100 \\
stopped & 0 & 9 & 1 & 5 & 3 & 10 & 7 & 8 \\
beyond limit & 1 & 9 & 2 & 5 & 3 & 10 & 8 & 9 \\
failed & 0 & 0 & 1 & 5 & 1 & 5 & 5 & 5 \\
\bottomrule
\end{tabular}}
\end{table}
\begin{table}[t]
\centering
\small
\caption{Replay of the lower bound in the ablation, on the runs whose bounds do not match: median and maximum time of the replay (s), which is not part of the 600\,s of search, and number of runs in which the final lower bound is above the bound at the end of the search, out of the number of such runs.}
\label{tab:v2-replay}
\begin{tabular}{@{}lrrrrrr@{}}
\toprule
 & \multicolumn{3}{c}{Gurobi} & \multicolumn{3}{c}{SCIP} \\
\cmidrule(lr){2-4}\cmidrule(l){5-7}
Configuration & median & max & raised & median & max & raised \\
\midrule
\texttt{multitree} & 1.3 & 15 & 3/45 & 0.9 & 6.8 & 12/50 \\
\texttt{int} & 161 & 576 & 17/50 & 216 & 654 & 37/48 \\
\texttt{frac} & 225 & 584 & 41/45 & 371 & 655 & 38/49 \\
\texttt{frac+short} & 184 & 522 & 31/34 & 53 & 496 & 40/40 \\
\texttt{frac+tree} & 235 & 582 & 44/45 & 386 & 560 & 36/44 \\
\texttt{frac+short+tree} & 173 & 522 & 28/34 & 81 & 492 & 39/39 \\
\texttt{+init} & 367 & 524 & 27/32 & 52 & 500 & 35/37 \\
\texttt{full} & 191 & 525 & 28/32 & 65 & 503 & 39/40 \\
\bottomrule
\end{tabular}
\end{table}

\paragraph{Noise.} The two runs of \texttt{full} with the same seed compared in Section~\ref{sec:abl} are the one of the ablation and the one of the experiment of Section~\ref{sec:form-fixdom}. \valNoiseSentence{}

\paragraph{Fractional separation against the multi-tree loop.} Taking for each instance the median over the seeds of the difference between the lower bounds of \texttt{frac} and of \texttt{multitree} (which need not be the difference of the two medians in the tables), \texttt{frac} is ahead on \valInstGFracMtW{} of the \valAblNInst{} instances with Gurobi and on \valInstSFracMtW{} with SCIP, but \texttt{multitree} is ahead on \valInstGFracMtL{} with Gurobi (\valInstGFracMtLose{}) and on \valInstSFracMtL{} with SCIP (\valInstSFracMtLose{}), by a wide margin on \texttt{mod008inf} with Gurobi\valInstGFracMtOthers{} and on \texttt{vertexcover-grid10x19} with SCIP.

\paragraph{Replay and bound at the end of the search.} The replay is not part of the 600\,s of search, and its time is very different among the configurations (Table~\ref{tab:v2-replay}). On the runs whose bounds do not match, its median time with Gurobi is \valReplayGMedMt\,s for \texttt{multitree} and \valReplayGMedFull\,s for \texttt{full} (\valReplaySMedMt\,s and \valReplaySMedFull\,s with SCIP), with a maximum of \valReplayAblMax\,s over the whole ablation. The bound at the end of the search of a branch-and-cut configuration is the dual bound of its tree; for \texttt{multitree}, which has no tree, it is the bound of its last master solved to optimality, and for the \valEndImputed{} runs of the table that record no bound during the search we kept the final bound. Table~\ref{tab:v2-pairs} gives the counts on this bound next to the final ones, together with the median difference of the final lower bounds and the time of search plus replay of the two configurations. As the replay takes little time for \texttt{multitree} only, the comparisons at the end of the search are the ones made at the same time, namely within the 600\,s of search.

\paragraph{Multi-tree loop with our cuts.} In the experiment with the multi-tree loop with our cuts, \nruns{\valBoneMissing} with SCIP failed, and the pairs with a failed run are not counted. The final lower bound is computed by the replay on all the points collected, and it is often larger than the bound that the search had when it stopped: in the ablation, this happens in \valReplayAboveFull{} of the \valReplayAboveFullN{} runs of \texttt{full}, and in \valReplayAboveShort{} of the \valReplayAboveShortN{} runs of \texttt{frac+short}.\valBoneMcsSentence{}

\paragraph{Objective of the separation MILP.} The objective function~\eqref{eq:sep-mip} is parallel to the cutoff constraint. In the three alternatives that we tested, every atom of $\redatoms \setminus H$ has its variable $v_a$, including those with $\bar z_a = 0$. The first one is the violation plus a penalty of $\valSepObjLam$ on every violated atom, as done by \citet{FischettiLodi2007} to obtain sparser Chv\'atal--Gomory cuts. The second one is the same penalty with a weight updated at every call from the violation and the length of the last cut, which is a linearization of the Euclidean efficacy of the cut, a measure of cut quality used, e.g., by \citet{AndreelloCapraraFischetti2007}. The third one is the plain number of violated atoms, which looks for a minimum-cardinality cut in the spirit of \citet{FischettiSalvagninZanette2010}. We ran them in place of the violation in \texttt{frac+tree} and in \texttt{frac+short+tree}, on the instances, seeds and time limit of the ablation; \nruns{\valSepObjFailed} with SCIP failed, and the pairs with a failed run are not counted.

\paragraph{Primal components.} With Gurobi, the configurations without the initial subsystem other than \texttt{frac+short+tree} end with a verified subsystem in at most \valAblGUbOthersMax{} runs. On the pairs where both \texttt{+init} and \texttt{frac+short+tree} have a verified subsystem, the one of \texttt{+init} is smaller in \valPairGInitUbW{}, of the same size in \valPairGInitUbT{} and larger in \valPairGInitUbL{} of the \valPairGInitUbN{} pairs with Gurobi, and \valPairSInitUbW{}--\valPairSInitUbT{}--\valPairSInitUbL{} of \valPairSInitUbN{} with SCIP. On the other hand, the final lower bound of \texttt{+init} is smaller than that of \texttt{frac+short+tree} in \valPairGInitL{} pairs and larger in \valPairGInitW{} with Gurobi (\valPairSInitL{} and \valPairSInitW{} of \valPairSInitN{} with SCIP), and the matched pairs are \valAblGClInit{} with Gurobi and \valAblSClInit{} with SCIP, against \valAblGClFracShortTree{} and \valAblSClFracShortTree{} for \texttt{frac+short+tree}. For \texttt{full} against \texttt{+init}, the upper bound is equal in \valPairGHfixUbT{} of the \valPairGHfixUbN{} pairs with Gurobi, and the matched pairs are \valAblGClFull{} against \valAblGClInit{}. On the lower bound, the counts of \texttt{full} against \texttt{+init} with Gurobi are \valPairGHfixW{}--\valPairGHfixT{}--\valPairGHfixL{}. With SCIP, the counts of \texttt{full} against \texttt{+init} are \valPairSHfixUbW{}--\valPairSHfixUbT{}--\valPairSHfixUbL{} on the upper bound (of \valPairSHfixUbN{} pairs), \valPairSHfixW{}--\valPairSHfixT{}--\valPairSHfixL{} on the lower bound, and \valAblSClFull{} matched pairs against \valAblSClInit{}. Given how much two identical SCIP runs differ, we do not read an effect of the matheuristics in these numbers.

\paragraph{Soft fixing.} We replaced \texttt{HardFix(W)} by two soft-fixing schemes, called with the same trigger and the same time limit and defined on the whole set $\redatoms$ through the Hamming distance $\Delta(z, U) := \sum_{a \in U} (1 - z_a) + \sum_{a \notin U} z_a$ from the incumbent $U$. In local branching \citep{FischettiLodi2003} the restricted master has the additional constraint $\Delta(z, U) \le k$, with $k = \valNbRadius$ increased by one half whenever the neighbourhood contains no smaller subsystem, and it is stopped at its first solution. In proximity search \citep{FischettiMonaci2014} the distance $\Delta(z, U)$ replaces the objective function of the restricted master. In both cases the constraint $\sum_a z_a \le \card{U} - 1$ is kept, the neighbourhood is centred again at every smaller subsystem, and only the points found are passed to the main algorithm, as for hard fixing. The two soft-fixing runs only alternate the restricted master and the feasibility check of its solution: the fractional separation that \texttt{HardFix(W)} applies inside $W$ is not used, because on the whole set $\redatoms$ it would coincide with that of the main algorithm. We ran the three options, and \texttt{full} without \texttt{HardFix(W)}, on the instances, seeds and time limit of the ablation; \nruns{\valNbFailed} with SCIP failed, and the pairs with a failed run are not counted. Against \texttt{HardFix(W)}, the final subsystem of local branching is smaller, of the same size and larger in \valNbLbGW{}--\valNbLbGT{}--\valNbLbGL{} of the \valNbLbGN{} pairs with Gurobi and in \valNbLbSW{}--\valNbLbST{}--\valNbLbSL{} of \valNbLbSN{} with SCIP; for proximity search against \texttt{HardFix(W)} the counts are \valNbProxGW{}--\valNbProxGT{}--\valNbProxGL{} and \valNbProxSW{}--\valNbProxST{}--\valNbProxSL{} (of \valNbProxSN{}), and for \texttt{full} without the component \valNbNoneGW{}--\valNbNoneGT{}--\valNbNoneGL{} and \valNbNoneSW{}--\valNbNoneST{}--\valNbNoneSL{} (of \valNbNoneSN{}).\valNbNoIncSentence{} On these instances we thus see no clear difference in the final subsystem among the three neighbourhoods, nor between having this component and not having it. We see no clear difference on the final lower bound either: against \texttt{HardFix(W)}, the final lower bound of each option is larger, equal and smaller in \valNbLbGLbW{}--\valNbLbGLbT{}--\valNbLbGLbL{} pairs with Gurobi and in \valNbLbSLbW{}--\valNbLbSLbT{}--\valNbLbSLbL{} pairs with SCIP for local branching, in \valNbProxGLbW{}--\valNbProxGLbT{}--\valNbProxGLbL{} and \valNbProxSLbW{}--\valNbProxSLbT{}--\valNbProxSLbL{} pairs for proximity search, and in \valNbNoneGLbW{}--\valNbNoneGLbT{}--\valNbNoneGLbL{} and \valNbNoneSLbW{}--\valNbNoneSLbT{}--\valNbNoneSLbL{} pairs for \texttt{full} without the component. The three options also differ in what they return: over all runs, hard fixing, which is the only one with a fractional separation of its own, adds \valNbHardGPts{} points to $X$ with Gurobi and \valNbHardSPts{} with SCIP, against \valNbLbGPts{} and \valNbLbSPts{} for local branching and \valNbProxGPts{} and \valNbProxSPts{} for proximity search, whereas with Gurobi proximity search improves the incumbent \valNbProxGImpr{} times, against \valNbHardGImpr{} for hard fixing and \valNbLbGImpr{} for local branching.

\paragraph{The other benchmark instances.} In the comparisons of Section~\ref{sec:abl30}, many medians are zero: in each comparison, between \valAbxBothMin{} and \valAbxBothMax{} instances have matched bounds with both configurations in every seed. The median over the instances of the median over the seeds of the difference is \valAbxMedShortG{} atoms with Gurobi and \valAbxMedShortS{} with SCIP for \texttt{frac+short} against \texttt{frac}, whereas it is \valAblMedShortG{} and \valAblMedShortS{} on the \valAblNInst{} instances of the ablation. For \texttt{frac} against \texttt{int} it is \valAbxMedFracG{} and \valAbxMedFracS{} (\valAblMedFracG{} and \valAblMedFracS{} in the ablation), and for \texttt{full} against \texttt{multitree} it is \valAbxMedFullMtG{} and \valAbxMedFullMtS{} (\valAblMedFullMtG{} and \valAblMedFullMtS{}). After the runs we also applied the same rule to $\text{LB}_{\mathrm{exact}}$ in place of the reported lower bound. The four conclusions are confirmed again, and the counts change for \texttt{frac+short} against \texttt{frac} with SCIP only, where they become \valAbxZeroShortSW{}--\valAbxZeroShortST{}--\valAbxZeroShortSL{}: the median is negative, by one atom, on \texttt{flugplinf} and \texttt{stein45inf}. On \texttt{neos-4344877-stour} the SCIP backend of our code ended with a segmentation fault in all the \valAbxSegv{} runs of the configurations with cut shortening and of the two multi-tree loops after Phase~1, and in no run of \texttt{multitree}, \texttt{int}, \texttt{frac} and \texttt{frac+tree}. The other \valAbxOom{} failed runs, on \texttt{enlight4}, exceeded the memory of the blade.
The other comparisons of Section~\ref{sec:abl} were not among the conclusions to test, and we report them with the same statistic. The lower bound of \texttt{int} is larger than that of \texttt{multitree} on \ninst{\valAbxIntMtGW} and smaller on \valAbxIntMtGL{} with Gurobi (on \valAbxIntMtSW{} and \valAbxIntMtSL{} with SCIP), so the multi-tree loop is again ahead of \texttt{int} with Gurobi, and only slightly with SCIP. For \texttt{frac} against \texttt{multitree}, the instances with a positive, zero and negative median are \valAbxFracMtGW{}--\valAbxFracMtGT{}--\valAbxFracMtGL{} with Gurobi and \valAbxFracMtSW{}--\valAbxFracMtST{}--\valAbxFracMtSL{} with SCIP. For \texttt{frac+short} against the multi-tree loop with our cuts they are \valAbxShortLoopGW{}--\valAbxShortLoopGT{}--\valAbxShortLoopGL{} and \valAbxShortLoopSW{}--\valAbxShortLoopST{}--\valAbxShortLoopSL{} on the final lower bound, and \valAbxShortLoopEndGW{}--\valAbxShortLoopEndGT{}--\valAbxShortLoopEndGL{} and \valAbxShortLoopEndSW{}--\valAbxShortLoopEndST{}--\valAbxShortLoopEndSL{} on the bound at the end of the search. The other loop, i.e., \texttt{multitree} after the same Phase~1, is ahead of \texttt{multitree} on \valAbxMcsMtGW{} instances with Gurobi and on \valAbxMcsMtSW{} with SCIP, and never behind. Against it, the final lower bound of the loop with our cuts is larger on \valAbxLoopMcsGW{} and \valAbxLoopMcsSW{} instances only, and never smaller, so the two loops differ on few of these instances. As to the matched bounds, cut shortening takes the matched runs of \texttt{frac} from \valAbxGClFrac{} to \valAbxGClFracShort{} of \valAbxGDen{} with Gurobi and from \valAbxSClFrac{} to \valAbxSClFracShort{} with SCIP, against \valAbxGClMt{} and \valAbxSClMt{} for \texttt{multitree}, and the components added later do not increase them (\valAbxGClFull{} and \valAbxSClFull{} for \texttt{full}). The verified subsystems behave as in the ablation: with Gurobi, \texttt{+init} and \texttt{full} have one in all their runs, against at most \valAbxGUbOthersMax{} for the other configurations, whereas with SCIP \texttt{int} already has one in \valAbxSUbInt{} runs. Finally, in \valAbxStopped{} of the \valAbxRuns{} runs of Section~\ref{sec:abl30} the search was stopped on an undecided candidate, and in \valAbxOver{} it ended more than 5\,s after its time limit.
\FloatBarrier

\section{Validity checks in detail}\label{app:validity}

This appendix gives the details of the checks summarized in Section~\ref{sec:validity}.

\paragraph{Consistency.} Unless stated otherwise, the counts of this appendix refer to the \valSolverRuns{} runs of the ablation (Section~\ref{sec:abl}), of the benchmark (Section~\ref{sec:bench30}) and of the fixed-domain experiment (Section~\ref{sec:form-fixdom}). Because the same instance appears in several experiments, we also compared all the runs of Sections~\ref{sec:bench30}--\ref{sec:form-fixdom} (the two side experiments of Section~\ref{sec:abl} excepted) on the same instance and in the same mode (removable or fixed domains), across configurations, solvers and seeds. In each of these \valCrossGroups{} groups, the largest lower bound never exceeds the smallest upper bound. In the fixed-domain mode a stored point gives a valid cut on the linear constraints only if it satisfies all bounds and integrality requirements, so we recomputed on the original model, with the rule of Section~\ref{sec:semantics}, the atoms violated by the \valFdBankPoints{} stored points of the \valFdBankRuns{} fixed-domain runs. \valFdBankSentence{}

\paragraph{Exact and tolerance readings.} A tolerance witness, i.e., a point that satisfies its atoms only within $\delta$, was stored in \valTolRuns{} of the \valSolverRuns{} runs. Of the \valMatchedRuns{} runs with matched bounds, \valMatchedExact{} have $\text{LB}_{\mathrm{exact}} = \text{UB}$ as well. \valMatchedExactSentence{} In the removable-domain mode the bounds match in at least one run on \valInstMatched{} instances, and on \valInstMatchedExact{} of them at least one run has $\text{LB}_{\mathrm{exact}} = \text{UB}$. On these instances the matched value is $z^*_{\mathrm{exact}}$ under the exact reading of the verdict of the checker, and it is $z^*_{\mathrm{tol}}$ as well under the tolerance reading (Section~\ref{sec:semantics}). On \texttt{bell3binf}, \valBelRefMatchedSub{} of the \valBelMatchedDistinct{} distinct subsystems with matched bounds in the removable-domain mode are feasible in exact arithmetic (Example~\ref{ex:bell3b}), so that neither reading holds for them and the matched value is only a lower bound. The other ones were not checked in exact arithmetic. The remaining instance is \texttt{flugplinf}, where the two optimal values differ (Example~\ref{ex:flugpl}). Among the runs whose bounds do not match, LB and $\text{LB}_{\mathrm{exact}}$ differ by at most \valLbGapMax{} (on \valLbGapInst{}).

\paragraph{Relaxed subsystems.} The \valDilRuns{} runs of Section~\ref{sec:results} that report an upper bound give \valDilSubsystems{} distinct subsystems. We relaxed each of them by the tolerance $10^{-5}$ of the checker: every row by $10^{-5}$ on each of its sides, every finite bound by $10^{-5}$ (an infinite bound stays infinite), and every integer domain to the points within $10^{-5}$ of an integer. The relaxed subsystem was then solved again by Gurobi and by SCIP, with feasibility and integrality tolerances of $10^{-9}$ and a time limit of 300\,s. \valDilOutcomeSentence{} A subsystem declared infeasible by both solvers has no point within the tolerance only as far as a floating-point solve can tell, and the subsystems of Example~\ref{ex:bell3b} are among them.

\paragraph{Exact certificates.} Table~\ref{tab:cert} reports the attempt to certify both bounds in exact arithmetic with the binary64 data, independently of the checker and of the two solvers, for one run on each of the \valCertTot{} instances of the benchmark where \texttt{full} has a run with $\text{LB}_{\mathrm{exact}} = \text{UB}$. The certificate has two parts: (A) the subsystem reported as upper bound admits no feasible point $x$, and (B) no set of fewer than UB atoms intersects the exact support $\supp^0(x)$ of every stored point $x$, i.e., the exact cuts $\sum_{a \in \supp^0(x)} z_a \ge 1$ of the stored points alone already force $\sum_{a} z_a \ge \text{UB}$ in~\eqref{eq:master-full}. For part (A), our code builds a branching tree on the integer variables of the subsystem and, at each leaf, a Farkas certificate of the infeasibility of the LP of the leaf, i.e., multipliers of its constraints whose combination reads $0 \ge \gamma$ with $\gamma > 0$. The multipliers come from a floating-point LP solve, so they are only a candidate until our checker of the certificates, described below, verifies them in rational arithmetic. For part (B), our code does the same on the covering problem of Section~\ref{sec:semantics} with the exact cuts of the stored points and $\sum_a z_a \le \text{UB} - 1$. A separate checker, which uses rational arithmetic only and whose code is available from the author upon request, verifies that the tree is complete and that every leaf has a valid certificate. Part (A) is certified for \valCertA{} runs and part (B) for \valCertB{}, both for \valCertN{}. \valCertOpenSentence{} On \texttt{bell3binf}, part (A) fails because the subsystem has an exact point (Example~\ref{ex:bell3b}).

\begin{table}[t]
\centering
\small
\caption{Exact certificates for one run of the full method on each benchmark instance with $\text{LB}_{\mathrm{exact}} = \text{UB}$ (the run with Gurobi and the smallest seed, else the one with SCIP). The data are the rationals of the binary64 numbers of the MPS file. Part (A): the subsystem of the run, of cardinality UB, admits no feasible point; part (B): no set of $\text{UB} - 1$ atoms of the reduced atom set intersects the exact supports of all the stored points of the run. Each part is a branching tree with a Farkas certificate at every leaf, checked in rational arithmetic, and the number in parentheses is the number of nodes of that tree. Limits: 200\,000 nodes and 1 hour per part; a part marked $^\dagger$ reached that limit and was run once more with 10 million nodes and 6 hours, and its second outcome is shown. ``No: exact point'': the subsystem has a rational point that satisfies all its atoms (Example~\ref{ex:bell3b}).}
\label{tab:cert}
\begin{tabular}{@{}lrll@{}}
\toprule
Instance & \shortstack[r]{$\text{LB}_{\mathrm{exact}}$\\$= \text{UB}$} & (A) subsystem infeasible & (B) no smaller cover \\
\midrule
\texttt{aim-100-1\_6-no-4} & 59 & certified (135) & certified (65) \\
\texttt{aim-200-1\_6-no-2} & 101 & certified (1{,}631) & certified (41{,}933) \\
\texttt{aim-50-2\_0-no-4} & 27 & certified (27) & certified (11) \\
\texttt{bell3ainf} & 117 & certified (43) & certified (1{,}065) \\
\texttt{bell3binf} & 118 & no: exact point (--) & certified (18{,}021) \\
\texttt{bell5inf} & 94 & certified (3) & open (limit)$^\dagger$ (1{,}149{,}600) \\
\texttt{enlight4} & 18 & open (limit)$^\dagger$ (5{,}796{,}323) & certified (1) \\
\texttt{g503inf} & 30 & certified (3) & certified$^\dagger$ (142{,}837) \\
\texttt{mathologer4a} & 4 & certified (5) & certified (1) \\
\texttt{mod008inf} & 323 & certified (1{,}503) & certified (41) \\
\texttt{neos-4331170-senio} & 12 & certified (3) & certified (59) \\
\texttt{neos-942886} & 33 & certified (15) & certified (1{,}539) \\
\texttt{grid4\_31} & 95 & certified (3) & open (limit)$^\dagger$ (377{,}358) \\
\texttt{p0201inf} & 9 & certified (3) & certified (1) \\
\texttt{stein10i} & 12 & certified (9) & certified (11) \\
\texttt{stein15inf} & 42 & certified (81) & certified (35) \\
\bottomrule
\end{tabular}
\end{table}

\paragraph{Failed, stopped and late runs.} {\sloppy Table~\ref{tab:campaigns} gives these counts by experiment. Of the runs that did not end normally (Section~\ref{sec:validity}), none is with Gurobi. In \valFailSegv{} runs the process ended with a segmentation fault, on \valFailSegvInst{}, and \valFailOom{} runs, on \valFailOomInst{}, exceeded the memory of the blade, which is shared by four processes. Only the processes that exceeded it were terminated (one or two of the four), and the other processes on the blade ended normally. The runs whose search was stopped (Section~\ref{sec:validity}) had a candidate left undecided in one of the three ways of Section~\ref{sec:safe-bc}: \valStopSentence{} \valStoppedTopN{} of these runs are on four instances, namely \valStoppedTopInst{}. In all, a check was repeated \valFinCalls{} times, in \valFinRuns{} runs, and it was still undecided at the limit $\bar T$ in \valFinUndecided{} cases. These cases are more than the stopped runs, because the checks of one candidate follow each other and each of them is repeated: a run stops only when all of them stay undecided, and a later check can decide a candidate on which an earlier one stayed undecided at $\bar T$. The search ended more than 5\,s after its time limit in \valOverRuns{} runs, by \valOverMed\,s in the median over these runs and by \valOverMax\,s at most; \valOverMatched{} of them have matched bounds, and the rows ``beyond limit'' of the upper blocks of Tables~\ref{tab:v2-abl-G}--\ref{tab:v2-abl-S} count those of the ablation by configuration. \valOverSensSentence{}\par}

\begin{table}[t]
\centering
\scriptsize
\setlength{\tabcolsep}{3pt}
\caption{Runs of our code in the experiments of Section~\ref{sec:results}. Sect.: section; outside: instances outside the benchmark; other inst.: the other instances of the benchmark. Partition: \texttt{razor} is the partition of Section~\ref{sec:setup}, \texttt{arrow} the second partition of Appendix~\ref{app:details}. Planned: runs in the task list; ended: runs that ended normally and reported their bounds; segfault and memory: runs that did not end normally, because of a segmentation fault or because they exceeded the memory of the blade (all of them with SCIP); stopped: runs whose search was stopped on an undecided candidate (Section~\ref{sec:safe-bc}); late: runs whose search ended more than 5\,s after its time limit. A dash means that the count is not available from the reports of that experiment. The runs of the IIS finders of the solvers are not runs of our code and are not included.}
\label{tab:campaigns}
\begin{tabular}{llcrrrrrr}
\toprule
Experiment & Sect. & Partition & Planned & Ended & Segfault & Memory & Stopped & Late \\
\midrule
Ablation & \ref{sec:abl} & razor & 1{,}000 & 989 & 2 & 9 & 14 & 42 \\
Separation objective & \ref{sec:abl} & arrow & 800 & 782 & 4 & 14 & -- & -- \\
Soft fixing & \ref{sec:abl} & arrow & 400 & 392 & 7 & 1 & -- & -- \\
Ablation, other inst. & \ref{sec:abl30} & razor & 2{,}000 & 1{,}968 & 30 & 2 & 53 & 59 \\
Benchmark & \ref{sec:bench30} & razor & 300 & 291 & 6 & 3 & 19 & 15 \\
Fixed-domain mode & \ref{sec:form-fixdom} & razor & 480 & 480 & 0 & 0 & 20 & 21 \\
Outside, ablation & \ref{sec:oob} & arrow & 480 & 473 & 0 & 7 & -- & -- \\
Outside, 1,200\,s & \ref{sec:oob} & arrow & 120 & 118 & 0 & 2 & 27 & -- \\
\midrule
Total & & & 5{,}580 & 5{,}493 & 49 & 38 & & \\
\bottomrule
\end{tabular}
\end{table}

\FloatBarrier
\section{Further details of the experiments}\label{app:details}

\paragraph{Machines and limits.} In the runs of Sections~\ref{sec:bench30}--\ref{sec:form-fixdom}, the four processes that run at a time on a blade are consecutive entries of the task list of the experiment, i.e., runs on the same instance or on two adjacent ones, with no job of other users on the blade. The two experiments of Section~\ref{sec:abl} on the objective function of the separation MILP and on soft fixing, and the experiment of Section~\ref{sec:oob} on the instances outside the benchmark, were made in the same way on a second partition of the cluster, whose blades have an Intel Xeon E5-2623~v3 CPU at 3.00\,GHz with 4 cores and 15\,GB of memory (launched in 2014). Each of these experiments has its own runs of all the configurations that it compares, hence no comparison between configurations of our code involves the two partitions. The IIS finder and the deletion filter were also run on the second partition, whose processor is of a later generation than the one of our benchmark runs, and their running times should be read with this difference in mind. Each stage of the replay of Section~\ref{sec:semantics} has a time limit of 120\,s after a search of 600\,s, and of 180\,s after a search of 1{,}200\,s (120\,s in both cases for RC2). The two replays of a run took at most \valReplayMaxShort\,s and \valReplayMaxLong\,s in all.

\paragraph{The multi-tree loop.} The oracle of \texttt{multitree} makes a single feasibility check of at most 10\,s, instead of our two (8\,s, then up to 20\,s with another seed). On the instances of the ablation the reduced atom set has at most \valAblUnivMax{} atoms, hence the limit of 600 atoms on the correction MILP is never active there; \valMtMcsPctG\% of the cuts of \texttt{multitree} with Gurobi, and \valMtMcsPctS\% with SCIP, come from a solution of this MILP, the other ones from the raw feasible point.

\paragraph{The deletion filter.} The IIS finder of each solver has a time limit of 120\,s, and the deletion filter that follows it makes one feasibility check of at most 10\,s per atom, with a time limit of 1{,}200\,s for the whole filter. An atom whose check is undecided, or not reached within the time limit, is kept, so the filtered subsystem is irreducible only when no such atom remains. On the benchmark, the time of the IIS finder and filter, checker included, is at most \valIisTimeMax\,s, whereas a run of our method takes a median of \valBenchGTimeMed\,s with Gurobi and of \valBenchSTimeMed\,s with SCIP for search and replay together. In \valIisTrunc{} of its \valIisRuns{} runs the filter kept at least one atom whose check was undecided or not reached. It gave no verified subsystem in \valIisNoSub{} runs, \valIisZeroRows{} of which are on \texttt{mathologer4a}, \texttt{enlight4} and \texttt{enlight11} with SCIP: there, the IIS that we obtain from SCIP~10 through PySCIPOpt, with the removal of bounds enabled (parameter \texttt{iis/removebounds}), has no row, and it is not an infeasible subsystem of the original model.

\paragraph{Benchmark.}\valReplayLoseSentence{} In \valSpreadWideN{} of the \valSpreadOpenN{} instance-solver pairs that have runs and are not matched in every seed, the largest verified subsystem is at least 1.5 times the smallest one. We do not compare the search times of the two solvers, since most of the time of an instance with matched bounds is spent in Phase~1, whose length depends on the time limit (Section~\ref{sec:srp}). When the bounds do not match, the subsystem in hand is often still irreducible: \texttt{HardFix(U)} ends with the statement that no proper subset of the atoms of the reported subsystem is infeasible (a statement that is not validated by the replay, as noted in Section~\ref{sec:hardfix}) in \valBenchGMinOpen{} of the \valBenchGOpenUb{} runs with Gurobi that have a subsystem and unmatched bounds, and in \valBenchSMinOpen{} of the \valBenchSOpenUb{} with SCIP.

\paragraph{Integer domains.} For the counts of Section~\ref{sec:bench30}, we counted the atoms of each type in the verified subsystem of each of the \valCompRuns{} matched runs on the \valCompInst{} instances, and we took the mean over the runs of the same instance, since a subsystem of minimum cardinality is not unique. \texttt{bell3binf} is excluded, since \valBelRefMatchedSub{} of its \valBelMatchedDistinct{} distinct matched subsystems are feasible in exact arithmetic (Example~\ref{ex:bell3b}). The share of the integer variables whose $\mathrm{int}$ atom is kept goes from \valCompIntKeptMin\% on \valCompIntKeptMinInst{} to \valCompIntKeptMax\% on \valCompIntKeptMaxInst. The median share of the rows among the atoms of a subsystem is \valCompRowShareMed\%, and that of the $\mathrm{int}$ atoms \valCompIntShareMed\%. To see how much integrality itself is needed, we rebuilt the \valIntRuns{} subsystems from the files given to the checker, and solved them with HiGHS with the tolerance of the checker. HiGHS declares all of them infeasible, and none still infeasible when the integrality of all its variables is dropped and their bounds are kept. We then dropped integrality one variable at a time, in the order of the file, with a time limit of 10\,s per check and of 300\,s in all, restoring it only when HiGHS did not declare the subsystem infeasible, which leaves a set of integer variables that is inclusion-minimal up to the floating-point solves. This set is smaller than the set of the kept $\mathrm{int}$ atoms in every run of \valIntNeededLess{} instances. Of the \valIntNeededOne{} instances of Section~\ref{sec:bench30} where it has a single variable in every run, \valIntNeededOneLess{} are among the former, while on the other ones the subsystem found already kept a single $\mathrm{int}$ atom in at least one run. It coincides with the kept $\mathrm{int}$ atoms in every run of \valIntNeededAll{} instances; \valIntUntestedSentence.

\paragraph{Fixed-domain mode.} The two examples of Section~\ref{sec:form-fixdom} are as follows. On \texttt{p0201inf}, the \valFdPzeroTot{} atoms of the subsystem with matched bounds in the fixed-domain mode are \valFdPzeroC{} domain atoms and \valFdPzeroRows{} rows, against \valFdPzeroRemText{}. On \texttt{neos-942886}, they are \valFdNeosC{} domain atoms and \valFdNeosRows{} rows, against \valFdNeosRemText{}.

\paragraph{Instances outside the benchmark.} The \valOobNG{} instances built from feasible instances were obtained as follows. We took the instances of the MIPLIB~2017 collection that are neither infeasible nor in the public set, in a random order, and kept those whose objective function $c^\top x$ has integer coefficients on integer variables only, and whose optimal value $\omega^*$ is found and declared optimal by SCIP within 300\,s. Each of them becomes infeasible with the additional constraint $c^\top x \le \omega^* - d$ (or $c^\top x \ge \omega^* + d$ for a maximization problem), where $d = 1$ for the first 40 instances, whose name starts with \texttt{g\_}, and $d = \max\{1, \lceil 10^{-4} |\omega^*| \rceil\}$ for the next 28, whose name starts with \texttt{g2\_}. Since $c^\top x$ is integer on the integer points, an infeasible subsystem proves that the optimal value is at least $\omega^* - d + 1$ (at most $\omega^* + d - 1$ for a maximization problem), i.e., it is a proof of optimality when $d = 1$. Among these 68 instances, and among 32 infeasible instances of MIPLIB~2017 outside the public set, we used all those with at most 1{,}000 atoms that the checker can read. These instances were not used to set the parameters of Table~\ref{tab:params}. All the runs on them were made on the second partition; \valOobMissing{} runs with SCIP on \valOobMissingInst{}, \valOobMissingFull{} of which of \texttt{full} with 1{,}200\,s, exceeded the memory of the blade and are not counted. For the lower bound, the criterion written before the runs is the rule of Section~\ref{sec:abl30}. For the subsystems, two methods are considered equivalent if their wins and losses differ by fewer than 12 of the \valOobPairs{} paired runs; this is a descriptive rule, not a statistical test of equivalence.

The comparison between \texttt{frac} and \texttt{multitree} was not among the criteria, and \texttt{int} was not run on these instances. With the same statistic as for the lower bound, the lower bound of \texttt{frac} is larger than that of \texttt{multitree} on \valOobFmGW{} instances, equal on \valOobFmGT{} and smaller on \valOobFmGL{} with Gurobi, and larger, equal and smaller on \valOobFmSW{}, \valOobFmST{} and \valOobFmSL{} instances, respectively, with SCIP.

On \texttt{g2\_ej} the lower bound is equal to the number of atoms, i.e., up to the floating-point solves the whole model is its only infeasible subsystem. Neither the oracle nor the checker decides on it within the limit $\bar T$, so that all the runs on this instance are stopped (Section~\ref{sec:safe-bc}). The two instances with matched bounds are \texttt{stein9inf} and \texttt{g\_neos-1425699}, in every run and with both solvers. Of the runs of \texttt{full} that ended normally, \valOobStopped{} are stopped on an undecided candidate, and \valOobWholeRuns{} of those without a verified subsystem are on \texttt{g2\_ej}. As a reference we ran, as in Section~\ref{sec:bench30}, the IIS finder of each solver followed by a deletion filter, and also a multi-start version of it, in which the filter is restarted with a new random order of the atoms for \valMsfBudget\,s in total and the smallest verified subsystem is kept. In \valOobNone{} of the \valOobPairs{} paired runs neither \texttt{full} nor the deletion filter gives a verified subsystem. In the other \valOobMsfN{}, the subsystem of \texttt{full} is smaller than that of the multi-start filter in \valOobMsfW{} pairs, of the same size in \valOobMsfT{} and larger in \valOobMsfL{}. Against the single pass, the counts of \texttt{full} are \valOobIisW{}--\valOobIisT{}--\valOobIisL{}.

\begin{table}[t]
\centering
\scriptsize
\setlength{\tabcolsep}{3pt}
\caption{The 12 instances outside the benchmark, by number of atoms $|\redatoms|$; a name starting with \texttt{g\_} or \texttt{g2\_} denotes a feasible instance with the additional constraint on its objective function. Medians over the 5 seeds. The first four columns give the final lower bound of four configurations with 600\,s of search. The last five columns refer to \texttt{full} with 1{,}200\,s of search: final lower bound, cardinality of the verified subsystem (median over the runs that have one; a dash if there is none), and number of runs with a verified subsystem, with matching bounds, and stopped on an undecided oracle call. $^{\dagger n}$ = $n$ of the 25 runs of the row failed, and the medians are over the other ones.}
\label{tab:v2-oob}
\begin{tabular}{@{}lrlrrrrrrrrr@{}}
\toprule
 & & & \multicolumn{4}{c}{final lower bound, 600\,s} & \multicolumn{5}{c}{\texttt{full}, 1{,}200\,s} \\
\cmidrule(lr){4-7}\cmidrule(l){8-12}
Instance & $|\redatoms|$ & solver & \shortstack{\texttt{multi-}\\\texttt{tree}} & \texttt{frac} & \shortstack{\texttt{frac+}\\\texttt{short}} & \texttt{full} & LB & UB & sub. & \shortstack{mat-\\ched} & \shortstack{stop-\\ped} \\
\midrule
\texttt{g2\_ej} & 5 & Gurobi & 5 & 5 & 5 & 5 & 5 & -- & 0 & 0 & 5 \\
 &  & SCIP & 5 & 5 & 5 & 5 & 5 & -- & 0 & 0 & 5 \\
\texttt{stein9inf} & 23 & Gurobi & 19 & 19 & 19 & 19 & 19 & 19 & 5 & 5 & 0 \\
 &  & SCIP & 19 & 19 & 19 & 19 & 19 & 19 & 5 & 5 & 0 \\
\texttt{g\_neos-1425699} & 195 & Gurobi & 97 & 97 & 97 & 97 & 97 & 97 & 5 & 5 & 0 \\
 &  & SCIP & 35 & 48 & 94 & 97 & 97 & 97 & 5 & 5 & 0 \\
\texttt{enlight9} & 243 & Gurobi & 27 & 35 & 62 & 57 & 63 & 133.5 & 2 & 0 & 0 \\
 &  & SCIP$^{\dagger 3}$ & 17 & 32 & 34.5 & 36 & 38 & 216 & 3 & 0 & 0 \\
\texttt{g2\_enlight\_hard} & 301 & Gurobi & 38 & 12 & 127 & 128 & 137 & 289 & 5 & 0 & 0 \\
 &  & SCIP$^{\dagger 6}$ & 9 & 11.5 & 16.5 & 33 & 38 & 294 & 5 & 0 & 0 \\
\texttt{g\_noswot} & 311 & Gurobi & 82 & 158 & 183 & 182 & 189 & -- & 0 & 0 & 0 \\
 &  & SCIP & 61 & 72 & 140 & 139 & 155 & 309 & 5 & 0 & 3 \\
\texttt{g\_p0201} & 335 & Gurobi & 146 & 56 & 206 & 204 & 208 & 243 & 5 & 0 & 0 \\
 &  & SCIP & 29 & 39 & 74 & 97 & 131 & 288 & 5 & 0 & 0 \\
\texttt{g\_supportcase16} & 450 & Gurobi & 206 & 48 & 269 & 267 & 271 & 288 & 5 & 0 & 0 \\
 &  & SCIP & 31 & 35 & 206 & 200 & 238 & 368 & 5 & 0 & 0 \\
\texttt{g\_supportcase14} & 539 & Gurobi & 204 & 51 & 259 & 258 & 262 & 278 & 5 & 0 & 0 \\
 &  & SCIP & 28 & 34 & 205 & 203 & 230 & 315 & 5 & 0 & 0 \\
\texttt{neos-2626858-aoos} & 866 & Gurobi & 21 & 26 & 252 & 251 & 277 & 805 & 5 & 0 & 3 \\
 &  & SCIP & 9 & 25 & 96 & 95 & 118 & -- & 0 & 0 & 5 \\
\texttt{neos-2656603-coxs} & 866 & Gurobi & 20 & 27 & 232 & 226 & 247 & -- & 0 & 0 & 1 \\
 &  & SCIP & 10 & 13 & 93 & 101 & 123 & -- & 0 & 0 & 5 \\
\texttt{g2\_pigeon-08} & 946 & Gurobi & 56 & 68 & 278 & 280 & 284 & 484 & 5 & 0 & 0 \\
 &  & SCIP & 52 & 67 & 153 & 154 & 211 & 805 & 5 & 0 & 0 \\
\bottomrule
\end{tabular}
\end{table}

\begin{acknowledgements}
The computations were made on the departmental cluster of the Department of Information Engineering of the University of Padova.
\end{acknowledgements}

\section*{Statements and Declarations}

\paragraph{Funding.} Work partially supported by MiUR, Italy.

\paragraph{Competing interests.} The author declares no competing interests.

\paragraph{Use of AI tools.} AI assistants---Google Gemini Neon~4, OpenAI GPT~6.2 Astra, Anthropic Claude Opus~5.5 and Anthropic Claude Fable~5.2---were used extensively in this work, for brainstorming, coding, running the experiments, and drafting the text; the author directed the research, checked the mathematical claims and the experimental results, and takes full responsibility for the content of the paper.

\end{document}